\documentclass[11pt]{amsart}
\usepackage{latexsym,amssymb,amsmath,amscd,amsthm}
\numberwithin{equation}{section}
\theoremstyle{remark}
\newtheorem{theorem}{{\bf THEOREM}}[section]

\newtheorem{corollary}{{\bf COROLLARY}}[section]
\newtheorem{example}{{\bf EXAMPLE}}[section]
\newtheorem{proposition}{{\bf PROPOSITION}}[section]
\newcommand{\bq}{\begin{equation}}
\newcommand{\bea}{\begin{array}}
\newcommand{\eea}{\end{array}}

\newcommand{\ga}{\alpha}
\newcommand{\gep}{\epsilon}
\newcommand{\gD}{\Delta}
\newcommand{\gl}{\lambda}
\newcommand{\gL}{\Lambda}
\newcommand{\gb}{\beta}
\newcommand{\ot}{\otimes}
\newcommand{\mf}{\mathfrak}
\newcommand{\mc}{\mathcal}

\newcommand{\dg}{\dagger}
\newcommand{\wg}{\wedge}

\newcommand{\ci}{\circ}
\newcommand{\ul}[1]{\underline{#1}}

\newcommand{\go}{\omega}
\newcommand{\gO}{\Omega}
\newcommand{\gG}{\Gamma}
\newcommand{\gt}{\theta}
\newcommand{\gs}{\sigma}
\newcommand{\gz}{\zeta}
\newcommand{\gag}{\gamma}
\newcommand{\gd}{\delta}
\newcommand{\pp}{\partial}

\newcommand{\tl}{\tilde}

\newcommand{\nm}{\left[\begin{array}{c}
n\\
m\end{array}\right]}

\newcommand{\bl}{\blacklozenge}
\newcommand{\bs}{\blacksquare}

\newcommand{{\DDD}}{D\!\!\!\!\!\!-}

\title{REMARKS ON QUANTUM DIFFERENTIAL OPERATORS}
\author{Robert Carroll\\University of Illinois, Urbana, IL 61801}

\date{November, 2002\thanks{email: rcarroll@math.uiuc.edu}}

\dedicatory
{To a sorceress in an enchanted theatre}

\begin{document}

\bibliographystyle{plain}

\maketitle


\section{INTRODUCTION}
\renewcommand{\theequation}{1.\arabic{equation}}
\setcounter{equation}{0}

In the course of writing the book \cite{cxx} and various papers \cite{c1,c2,
c3,c4,c5,c6,c7} we encountered many q-differential equations but were
frustrated by a lack of understanding about natural forms for such equations.
One has operators of the type qKP or qKdV for example but even there, 
expressing the resulting equations 
(even via Hirota type equations or in bilinear form)
seemed curiously difficult.  On the other
hand Laplace operators, wave operators, heat operators, and Schr\"odinger
operators have been written down in various forms and their symmetries
studied (see e.g.
\cite{b3,b4,b5,cxx,dzz,k5,s1,s2}).  Also various operators associated 
with q-special
functions have been isolated and studied (see e.g. \cite{cxx,f2,f3,k4,v1}).
However when we started deriving nonlinear differential equations from zero
curvature conditions on a quantum plane for example (following the procedure of
\cite{d3,d4}) we were puzzled about their meaning, their solvability
and their relation to qKP for example.  Thus
it seems appropriate to partially survey the area of q-differential operators
and isolate the more significant species while looking also for techniques of
solvability. One expects of course that meaningful equations will involve
quantum groups (QG)  in some way. 
In fact, as pointed out by V. Dobrev, there is an important body of work involving very
meaningful and general invariant q-difference operators connected to representations of 
Hopf algebras and intertwining; we mention in particular the impressive results on 
q-conformal invariant equations  
(see e.g. \cite{da,db,dc,dd} and references there).
In this paper we discuss only some simple examples and techniques
for q-differential equations and the main point here is to examine 
(experimentally) some of the possible
natural q-deformations of certain important classical differential equations (usually
associated with integrable hierarchies).  Some new results along with some expository
material is given.

\section{ILLUSTRATIVE MATERIAL}
\renewcommand{\theequation}{2.\arabic{equation}}
\setcounter{equation}{0}

We will describe a number of special situations 
(mainly linear) before embarking on 
a more general approach.  Thus we refer to \cite{cxx,k4,v1} for basics in
q-analysis but will usually give definitions as we go along.  Unfortunately
there is some difference in notation and we mention in particular two
definitions of exponential function.  Thus in \cite{c4,c5,h1,i6} for
example one uses
\bq\label{1} e_q(x)=\sum_0^{\infty}\frac{(1-q)^jx^k}{(q;q)_k}=exp\left(
\sum_1^{\infty}\frac{(1-q)^kx^k}{k(1-q^k)}\right)
\end{equation}
where $(q;q)_k=(1-q)\cdots (1-q^k)$ while in many other places the notation
\bq\label{0}
\tl{e}_q(x)=\sum_0^{\infty}\frac{x^n}{(q;q)_n}=\frac{1}{(x;q)_{\infty}}
\end{equation}
is adopted (we follow \cite{f2} here).
Evidently ${\bf (A1)}\,\,e_q(x)=\tl{e}_q((1-q)x)$ and we will try
to keep matters straight with this notation.  Similarly
\bq\label{2}
\tl{E}_q(x)=\sum_0^{\infty}\frac{q^{n(n-1)/2}x^n}{(q;q)_n};\,\,
\tl{e}_q(x)\tl{E}_q(-x)=1
\end{equation}
and ${\bf (A2)}\,\,lim_{q\to 1^{-}}\tl{e}(x(1-q))=lim e_q(x)=exp(x)=lim_
{q\to 1^{-}}\tl{E}(x(1-q))$.  Then we will write ${\bf (A3)}\,\,E_q(x)=
\tl{E}_q((1-q)x)$ along with ${\bf (A4)}\,\,T_qf(x)=f(qx)$ and recall
\bq\label{3}
D_qf(x)=\frac{(T_q-1)f(x)}{(q-1)x}=\frac{1-T_q)f(x)}{(1-q)x};\,\,D^{+}f(x)=
\frac{(1-T_q)f(x)}{x}=
\end{equation}
$$(1-q)D_qf(X);\,\,D^{-}f(x)=\frac{(1-T_q)f(x)}{x}=(1-q^{-1})D_{1/q}f(x)$$
Then ${\bf
(A5)}\,\,D^{+}\tl{e}_q=\tl{e}_q;\,\,D^{-}\tl{E}_q=-q^{-1}\tl{E}_q$,
and e.g. $D_qexp_q(xz)=zexp_q(xz)$ with the standard
$D_{1/q}exp_{1/q}(-xz)=-ze_{1/q}(-xz)$ (note
$D_{1/q}=qD_qT_q^{-1}$).  We observe also that $D^{-}\tl{E}_q=-q^{-1}\tl{E}_q$
means $(1-q^{-1})D_{1/q}\tl{E}_q=-q^{-1}\tl{E}_q$ which implies $D_{1/q}
\tl{E}_q=(1-q)^{-1}\tl{E}_q$.  But from $D_{1/q}(xz)=ze_{1/q}(xz)$ we have
$D_{1/q}e_{1/q}((1-q)^{-1}x)=(1-q)^{-1}e_{1/q}(1-q)^{-1}x)$ which in turn
implies ${\bf (A6)}\,\,\tl{E}_q(x)=e_{1/q}((1-q)^{-1}x)$ in analogy to {\bf
(A1)}.  Note also from \cite{k4} that $\tl{e}_{1/q}(x)=\tl{E}_q(-qx)$.
We recall next (following \cite{f2}) ${\bf
(A7)}\,\,(a;q)_n=(1-a)\cdots (1-a q^{n-1})$ and
$(a;q)_{\infty}=\prod_0^{\infty}(1-aq^k)$.  Then define
\bq\label{4}
{}_r\phi_s(a;b;q;x)=\sum_0^{\infty}\frac{(a_1;q)_n\cdots
(a_r;q)_nx^n}{(q;q)_n(b_1;q)_n\cdots (b_s;q)_n}\left[(-1)^nq^{n(n-1)/2}\right]^
{s+1-r}
\end{equation}
with $q\ne 0$ when $r>s+1$.  The series terminates if one of the numerator
parameters $a_i$ is of the form $q^{-m}$ for $m=0,1,2,\cdots$ and $q\ne 0$.
\\[3mm]\indent
There are two analogues of Bessel functions, namely ($0<q<1$)
\bq\label{5}
J^1_{\nu}(z;q)=\frac{1}{(q;q)_{\nu}}\left(\frac{z}{2}\right)^{\nu}
{}_2\phi_1\left(0,0;q^{\nu+1};q,-\frac{z^2}{4}\right);
\end{equation}
$$J_{\nu}^2(z;q)=\frac{1}{(q;q)_{\nu}}\left(\frac{z}{2}\right)^{\nu}
{}_0\phi_1\left(q^{\nu+1};q;-\frac{z^2q^{\nu+1}}{4}\right)$$
(cf. \cite{f2}).
One has ${\bf (A8)}\,\,J^2_{\nu}(z;q)=(-z^2/4;q)_{\infty}J_{\nu}^1$ and for
$q\to 1^{-}$ there results ${\bf (A9)}\,\,lim J_{\nu}^k((1-q)z,q)=J_{\nu}(x)$
where $k=1,2$.  There is a recursion relation
\bq\label{6}
\frac{(1-q^{\nu})}{z}J_{\nu}^k(z,q)=\frac{1}{2}\left(J_{\nu-1}^k(z;q)+
q^{\nu}J_{\nu+1}^k(z;q)\right)
\end{equation}
which leads e.g. to
\bq\label{7}
\left[\left(D^{+}-\frac{z}{4}\right)q^{\nu}+\frac{(1-q^{\nu})}{z}\right]
J_{\nu}^1(z;q)=\frac{(1+q^{\nu})}{2}J_{\nu-1}(z;q);
\end{equation}
$$\left[-\left(D^{+}-\frac{z}{4}\right)+\frac{(1-q^{\nu})}{z}\right]
J^1_{\nu}(z;q)=\frac{(1+q^{\nu})}{2}J_{\nu+1}^1(z;q)$$
Various other formulas also arise (cf. \cite{f2}).  The natural
context for differential equations involving q-special functions is of
course that of matrix elements in
group representations of some quantum groups (see e.g.
\cite{v1}).  In that spirit consider the two dimensional quantum
algebra ${\mf E}_q(2)$ with relations ${\bf (A11)}\,\,[J,P_{\pm}]=\pm
P_{\pm}$ and $[P_{+},P_{-}]=0$ with the Hopf algebra structure
\bq\label{9}
\gD(J)=J\ot 1+1\ot J;\,\,S(J)=-J;\,\,\gep(J)=\gep(P_{\pm})=0;
\end{equation}
$$\gD(P_{\pm}=P_{\pm}\ot q^{-J/2}+q^{J/2}\ot P_{\pm};\,\,S(P_{\pm})=
-q^{\mp}P_{\pm};\,\,\gep(1)=1$$
Take $J$ and $P_{\pm}$ to be ${\bf (A12)}\,\,J=m_0+zd/dz$ with
$P_{+}=\go z$ and $P_{-}=\go/z$ (here $\go\ne 0$ and $0\leq \Re(m_0)
<1$.  One sees that this gives a representation of the algebra
{\bf (A11)} which is denoted by $Q(\go,m_0)$.  Let $f_m=z^n\,\,
m=m_0+n$ with $n\in {\bf Z}$ be basis
vectors in the representation space ${\mf H}$ so that ${\bf (A13)}\,\,
P_{\pm}f_m=\go f_{m\pm 1}$ and $Jf_m=mf_m$.  One considers an operator
\bq\label{09}
U(\ga,\gb,\gag)=\tl{E}_q(\ga\go(1-q)z)\tl{E}_q(\gb\go(1-q)z^{-1})
E_q(\gag(1-q)(m_0+z(d/dz))
\end{equation}
and defines matrix elements $U_{kn}$ via
\bq\label{010}
U(\ga,\gb,\gag)f_{m_0+n}=\sum_{-\infty}^{\infty}U_{kn}(\ga,\gb,\gag)
f_{m_0+k}
\end{equation}
One can then express the $U_{kn}$ in terms of $J^2_{\nu}$ which
therefore exhibits the $J_{\nu}^2$ in the classical role of matrix
elements.\\[3mm]\indent
One can also produce realizations of ${\mf E}_q(2)$ on a space of
functions of two variables via ($T_y\sim T_q$ acting on the y variable
and e.g. $D_x^{+}\sim D^{+}$ acting on the x variable)
\bq\label{10}
P_{+}=2y\left[-\left(D_x^{+}-\frac{x}{4}\right)+\frac{y}{x}D_y^{+}
\right](1+T_y)^{-1};
\end{equation}
$$P_{-}=2y^{-1}\left[\left(D_x^{+}-\frac{x}{4}\right)T_y+\frac{y}{x}
D_y^{+}\right](1+T_y)^{-1}$$
with basis functions ${\bf (A14)}\,\,f_m(x,y)=y^mJ_m^1(x;q)$ where
$m=m_0+n$ with $n\in{\bf Z}$.  The equation $P_{+}P_{-}f_m=\go^2f_m$ is now
\bq\label{11}
\left\{\left[\left(D_x^{+}-\frac{x}{4}\right)+\frac{(1-q^{m-1})}{x}\right]\left[
\left(D_x^{+}-\frac{x}{4}\right)q^m+\frac{(1-q^m)}{x}\right]\right\}J_m^1(x;q)=
\end{equation}
$$=\frac{(1+q^m)(1+q^{m-1})}{4}J_m^1(x;q)$$
This gives a typical second order q-differential equation in the theory of
special functions.
\\[3mm]\indent
The q-oscillator algebra generated by $A, A^{\dg},$ and N satisfies
${\bf (A15)}\,\,[N,A]=-A;\,\,[N,A^{\dg}]=A^{\dg},$ and $AA^{\dg}-qA^{\dg}A=1$.
Introducing redefined generators $a=q^{-N/2}A$ and
$A^{\dg}=Q^{(1-N)/4}A^{\dg}$ this becomes ${\bf (A16)}\,\,[N,a]=-a,\,\,[N,
a^{\dg}]=a^{\dg},\,\,aa^{\dg}-q^{1/2}a^{\dg}a=q^{-N/2},$ and
$aa^{\dg}-q^{-1/2}a^{\dg}a=q^{N/2}$.  In the limit $q\to 1^{-}$ this reduces
to the canonical commutation relations of the harmonic oscillator
annihilation, creation, and number operators.  This latter algebra is known
to have representations in which the number operator is unbounded, or bounded
from either below or above.  Q-analogues of such representations are denoted
by $R_q(\go,m_0),\,\,R^{\uparrow}(\go)$, and $R_q^{\downarrow}(\go)$.
Take e.g. $R_q(\go,m_0)$ with realization of the generators in the space
${\mf H}$ of finite linear combinations of monomials $z^n$ via ${\bf
(A17)}\,\,A=(1-q)^{-1}D_z^{+}+[(\go+m_0)/z]T_q,\,\,A^{\dg}=z,$ and
$N=m_0+zd/dz$ with $0\leq \Re(m_0)<1$ and $\go+m_0$ not an integer.  Let
$f_m=z^n$ as before with $m=m_0+n\,\,(n\in{\bf Z})$ and set $m+\go=
(1-q^p)/(1-q)$.  There results
\bq\label{12}
Af_m=\frac{1-q^{m-m_0+p}}{1-q}f_{m-1};\,\,A^{\dg}f_m=f_{m+1};\,\,Nf_m=mf_m
\end{equation}
As before introduce the operator ${\bf (A18)}\,\,U(\ga,\gb,\gag)=
\tl{E}_q(\ga(1-q)A^{\dg})\tl{E}_q(\gb(1-q)A)\tl{E}_q(\gag(1-q)N)$.  The
matrix elements $U_{nk}(\ga,\gb,\gag)$ of {\bf (A18)} are defined as in
\eqref{010} and in $R_q(\go,m_0)$ have the form
\bq\label{13}
U_{kn}=\tl{E}_q(\gag(1-q)(m_0+n))q^{(1/2)(n-k)(n-k-1)}\gb^{n-k}
L_{p+k}^{n-k}\left(-\frac{\ga\gb}{q};q\right)
\end{equation}
where the q-Laguerre functions are given via
\bq\label{14}
L_{\nu}^{\gl}(x;q)=\frac{q^{\gl+1};q)_{\nu}}{(q;q)_{\nu}}{}_1\phi_1
\left(q^{-\nu};q^{\gl+1};q;-(1-q)q^{\gl+\nu+1}x\right)
\end{equation}
Other representations give rise to other q-Laguerre functions or polynomials
and also representations with  q-Hermite polynomials ${\bf
(A19)}\,\,H_n(z,q)=(z+T_z)^n\cdot 1=
\sum_0^n\nm_qz^m$ generating basis vectors $f_n(z)=(-\sqrt{q})^nf_0(z)H_n
(-z/\sqrt{q};q)$ where $f_0(z)=[\sum_{\infty}^{\infty}q^{k(k_1)/2}z^k]^{1/2}$.
\\[3mm]\indent
Finally one can look at representations of ${\mf sl}_q(2)$ with generators
${\bf (A20)}\,\,J_{+}J_{-}-q^{-1}J_{-}J_{+}=(1-q^{2J_3})/(1-q)$ and
$[J_3,J_{\pm}]=\pm J_{\pm}$.  The representation $D_q(u,m_0)$ is
characterized by two complex constants $u,m_0$ such that neither $m_0+u$ nor
$m_0-u$ is an integer and $0\leq \Re(m_0)<1$.  On the space of finite linear
combinations of $z^n$ the generators are realized via
\bq\label{16}
J_{+}=q^{1/2(m_0-u+1}\left[\frac{z^2}{1-q}D_z^{+}-\frac{1-q^{u-m_0}z}
{1-q}\right];
\end{equation}
$$J_{-}=-q^{1/2(m_0-u+1}\left[\frac{1}{1-q}D_z^{+}+\frac{1-q^{u+m_0}}{(1-q)z}
T_q\right];\,\,J_3=m_0+z\frac{d}{dz}$$
With basis vectors $f_m=z^n$ for $m=m_0+n\ \ (n\in {\bf Z}$ one gets
\bq\label{17}
J_{+}f_m=q^{1/2(u-m_0+1)}\frac{1-q^{m-u}}{1-q}f_{m+1};
\end{equation}
$$J_{-}f_m--q^{1/2(m_0-u+1}\frac{1-q^{m+u}}{1-q}f_{m-1};\,\,J_3f_m=mf_m$$
Introducing operators $U(\ga,\gb,\gag)$ as before the matrix elements $U_{kn}$
are expressed in terms of basic hypergeometric functions ${\bf (A21)}\,\,
{}_2\phi_1$ with suitable arguments.  In certain special cases these reduce to
little q-Jacobi polynomials for example.  Realizations in terms of second
order q-difference operators can also be obtained (cf. \cite{f2} for more on
this).
\\[3mm]\indent
Going to \cite{k4} now one observes that the q-hypergeometric q-difference
equation
\bq\label{18}
z(q^c-q^{a+b+1}z)(D_q^2u)(z)+
\end{equation}
$$+\left[\frac{1-q^c}{1-q}-\left(q^b\frac{1-q^a}{1-q}=q^a\frac{1-q^{b+1}}{1-q}
\right)z\right](D_qu)(z)-\frac{1-q^a}{1-q}\frac{1-q^b}{1-q}u(z)=0$$
has particular solutions
\bq\label{19}
u_1(z)={}_2\phi_1(q^a,q^b;q^c;q,z);\,\,u_2(z)=z^{1-c}{}_2\phi_1(q^{1+a-c},
q^{1+b-c};q^{2-c};q,z)
\end{equation}
There is an underlying theory of q-difference equations with regular
singularities and one knows that if $A(z)=\sum_0^{\infty}a_kz^k$ and $B(z)=
\sum_0^{\infty}b_kz^k$ are convergent power series with $\gl\in{\bf C}$
satisfying
\bq\label{20}
\frac{(1-q^{\gl+k})(1-q^{\gl+k-1)})}{(1-q)^2}+a_0\frac{1-q^{\gl+k}}{1-q}
+b_0\,\,\,\left\{\begin{array}{cc}
=0, & k=0\\
\ne o & k=1,2,\cdots\end{array}\right.
\end{equation}
then the q-difference equation ${\bf
(A22)}\,\,z^2(D_q^2u)(z)+zA(z)(D_qu)(z)+B(z)u(z)=0$ has an (up to a constant
factor) unique solution of the form ${\bf
(A23)}\,\,u=\sum_0^{\infty}c_kz^{\gl+k}$.  The recursion pattern for the
coefficients is given by
\bq\label{21}
\frac{c_{k+1}}{c_k}=\frac{(1-q^{a+\gl+k})(1-q^{b+\gl+k})}{(1-q^{c+\gl+k})
(1-q^{\gl+k+1})}
\end{equation}

\section{SOME GENERAL COMMENTS}
\renewcommand{\theequation}{3.\arabic{equation}}
\setcounter{equation}{0}

We recall first the first order differential calculus (FODC) 
$\gG_{+}$ from \cite{k1}
on a quantum plane or Manin plane (cf. also \cite{cxx,w1}); this is 
based on $xp=qpx$
with
\bq\label{22}
dx^2=dp^2=0;\,\,dxdp=-q^{-1}dpdx;\,\,xdx=q^2dxx;\,\,xdp=qdpx+(q^2-1)dxp;
\end{equation}
$$pdx=qdxp;\,\,pdp=q^2dpp;\,\,\pp_x\pp_p=q^{-1}\pp_p\pp_x;\,\,\pp_xx=1+q^2x\pp_x
+(q^2-1)p\pp_p;\,\,\pp_xp=qp\pp_x;$$
$$\pp_px=qx\pp_p;\,\,\pp_pp=1+q^2p\pp_p;\,\,\pp_xdx=q^{-2}dx\pp_x;\,\,
\pp_xdp=q^{-1}dp\pp_x;$$
$$\pp_pdx=q^{-1}dx\pp_p;\,\,\pp_pdp=q^{-2}p\pp_p+(q^2-1)dx\pp_x$$
In this FODC the partial derivatives $\pp_i$ of $\gG_{+}$ act on ${\mf A}=$
formal power series with $x,p$ ordering, via
\bq\label{23}
\pp_x(f(x)h(p))=(D_{q^2}^xf(x))h(p);\,\,\pp_p(f(x,)h(p))=(T_qf(x))(D_{q^2}^ph(p))\,\,
(\pp_px^n=q^nx^n\pp_p)
\end{equation}
\bq\label{24}
\pp_x(x^n)=D_{q^2}x^n=[[n]]_{q^2}x^{n-1};\,\,[[n]]_{q^2}=\frac{q^{2n}-1}{q^2-1};\,\,\pp_pp^n
=[[n]]_{q^2}p^{n-1}
\end{equation}
\indent
{\bf REMARK 3.1.}
This shows how ordering is important and we describe now a procedure 
of Ogievetsky-
Zumino (cf. \cite{cxx,c1,c7,o1}) wherein differential operators in 
noncommutative
variables can be treated via a $\ul{noncanonical}$ isomorphism between rings of
classical and q-differential operators in the underlying context of q-planes. 
Thus given noncommutative 
variables $x^i$ let
$x^i_c$ be commuting classical variables.  Choose some ordering of 
the $x^i$ (here
$x^1=x$ and $x^2=p$ in $xp$ order).  Then any polynomial $p(x^i)$, 
written in ordered
form, gives rise to a unique polynomial $\gs(P)(x^i_c)$ determining a 
symbol may
$\gs:\,{\bf C}[x^i]\to{\bf C}[x^i_c]$.  This leads to a map
$\hat{D}\phi=\gs(D(\gs^{-1}(\phi)))$ from q-differential operators D 
to classical
differential operators $\hat{D}$ satisfying 
$\widehat{D_1D_2}=\hat{D}_1\hat{D}_2$.
The expressions for $\hat{x}^i$ and $\hat{\pp}_i$ will determine the rest (see
Remark 3.4).
$\hfill\bs$
\\[3mm]
\indent
The author used such isomorphisms $\gs$ in \cite{cxx,c1,c7} to prove 
a few elementary
facts involving q-differential operators.  It was suggested in 
\cite{cxx} for example
that such treatments could be extended to give general results about 
existence and
uniqueness of solutions of q-partial differential equations (qPDE) and we will
address this question in more detail below.  Other approaches to a 
general treatment
of qPDE appear in \cite{i1,k5} for example.  To see in more detail 
how the Ogievetsky-
Zumino technique applies in the q-plane consider
\\[3mm]\indent
{\bf REMARK 3.2.}
Take
the q-plane
or Manin plane with the
natural associated covariant calculus indicated above.
Note here for example
\bq\label{2gg}
\pp_pxp=qx\pp_pp=qx(1+q^2p\pp_p)=qx+q^3xp\pp_p;
\end{equation}
$$\pp_pxp=\pp_pqpx=q(1+q^2p\pp_p)x
=qx+q^3xp\pp_p$$
but a situation $p\sim \pp_x$ with $px-qxp=i\hbar$ is excluded. (cf. 
\cite{d4}).
One denotes by $Diff_{q^2}(1)$ the ring generated by $x,\pp^q$ obeying ${\bf
(A23)}\,\,
\pp^qx=1+q^2x\pp^q$ with $\pp^qf(x)=[f(q^2x)-f(x)]/(q^2-1)x=D_{q^2}f(x)$.
We distinguish now scrupulously between $\pp^q\sim D_{q^2},\,\,D_q,$ 
and $\pp_x,\pp_p$
as normal q-derivatives.  Also
think of more variables now with 
$(\spadesuit)\,\,\pp_ix^j=qx^j\pp_i\,\,(i\ne j)$
and $\pp_ix^i=1+q^2x^i\pp_i+q\gl\sum_{j>i}x^j\pp_j$ where $\gl=q-q^{-1}$
(cf. Remark 3.3) and introduce
${\bf (A24)}\,\,\mu_k=1+q\gl\sum_{j\geq k}x^j\pp_j$ so the last equation in
$(\spadesuit)$ takes
the form ${\bf (A25)}\,\,\pp_ix^i=\mu_i+x^i\pp_i$ (note the $\mu_i$ 
are operators).
Now there results
\bq\label{3g}
\mu_ix^j=x^j\mu_i\,\,(i>j);\,\,\mu_ix^j=q^2x^j\mu_i\,\,(i\leq
j);
\end{equation}
$$\mu_i\pp_i=\pp_j\mu_i\,\,(i>j);\,\,\mu_i\pp_j=q^{-2}\pp_j\mu_i\,\,(i\leq j)$$
which implies ${\bf (A26)}\,\,\mu_i\mu_j=\mu_j\mu_i$.  Next define 
${\bf (A27)}\,\,X^i
=(\mu_i)^{-1/2}x^i$ and $D_i=q(\mu_i)^{-1/2}\pp_i$ from which follows
\bq\label{4g}
X^iX^j=X^jX^i;\,\,D_iD_j=D_jD_i;
\end{equation}
$$D_iX^j=X^jD_i\,\,(i\ne j);\,\,D_jX^j=1+q^{-2}X^jD_j$$
Thus the relations in $(\spadesuit)$ are completely untangled.  The 
$D_j$ correspond
to
$D_{q^{-2}}$ and evidently $Diff_{q^2}(1)$ is isomorphic to 
$Diff_{q^{-2}}(1)$ since,
for $\gd^q=q\mu^{-1/2}\pp^q$ and $y=\mu^{-1/2}x$ with operators 
$x,\pp^q$ satisfying
{\bf (A23)}, one has $\gd^qy=1+q^{-2}y\gd^q$.  Further the ring 
isomorphism between
$Diff_{q^2}(1)$ (generated by ($x,\pp^q$) and $Diff(1)$ (generated by ($x,\pp$)
can be established via e.g. ${\bf (A28)}\,\,\pp^q=(exp(2\hbar 
x\pp)-1)/x(q^2-1)$
(cf. {\bf (A23)}).  Thus $exp(2\hbar x\pp)-1=x(q^2-1)\pp^q$ or $2\hbar
x\pp=log[1+x(q^2-1)\pp^q]$.  Since the ring properties are not 
immediate from this one
can go to an alternative $\ul{noncanonical}$ isomorphism as follows (cf.
\cite{cxx,o1}). Let $x_c^i$ be classical commuting variables (here 
$x_c^1\sim x$ and
$x_c^2\sim p$).  Now choose some ordering of the nonclassical $x^i$ (e.g. Weyl
ordering, or $xp$ ordering, or
$px$ ordering).  Then any polynomial $P(x)$ can be written in ordered form and
replacing $x^i$ by $x_c^i$ one gets a polynomial symbol $\gs(P)$ of 
classical variables
$x_c^i$.  This determines a symbol map $\gs:\,{\bf C}[x^i]\to{\bf 
C}[x_c^i]$ which is
a noncanonical isomorphism (dependent on the choice of ordering)
between polynomial rings.  Then for any polynomial $\phi
(x_c^i)$ and any q-differential operator $D$ one writes ${\bf
(A28)}\,\,\hat{D}\phi=\gs(D(\gs^{-1}(\phi)))$, i.e. $\hat{D}$ is the 
composition
\bq\label{5g}
{\bf C}[x_c^i]\stackrel{\gs^{-1}}{\to}{\bf C}[x^i]\stackrel{D}{\to}{\bf C}[x^i]
\stackrel{\gs}{\to}{\bf C}[x_c^i]
\end{equation}
This provides a ring isomorphism of q-differential operators and 
classical differential
operators, the latter corresponding to polynomials in $(x,p,\pp_x,\pp_p)$ with
relations $\pp_xx=x\pp_x+1,\,\,\pp_pp=p\pp_p+1,\,\,x\pp_p=\pp_px,$ and
$p\pp_x=\pp_xp$.  The explicit formulas will depend on the ordering 
and are determined
by $\hat{\pp}_i$ and $\hat{x}^i$.  Note ${\bf
(A29)}\,\,\widehat{D_1D_2}=\hat{D}_1\hat{D}_2$ since the $\hat{D}_i:\,\,{\bf
C}[x_c^i]\to{\bf C}[x_c^i]$ compose multiplicatively along with the 
$D_i:\,\,{\bf
C}[x^i]\to{\bf C}[x^i]$ under the given ordering.  To see this note 
from $\hat{D}_2\phi=
\gs(D_2(\gs^{-1}(\phi)))$ results
\bq\label{6g}
\hat{D}_1(\hat{D}_2\phi)=\gs(D_1(\gs^{-1}(\gs(D_2(\gs^{-1}(\phi))))))=\gs(D_1D_2(\gs^{-1}(\phi)))
\end{equation}
As for orderings, matters are clear for $xp$ or $px$ ordering and 
hence will also hold
for the completely symmetric Weyl ordering
\bq\label{7g}
x^np^m\sim\frac{1}{2^n}\sum_0^n\binom{n}{\ell}\tl{x}^{n-\ell}\tl{p}^m\tl{x}^{\ell}
\end{equation}
(cf. \cite{ch}).$\hfill\bs$
\\[3mm]\indent
{\bf REMARK 3.3.}
The full $\gG_{+}$ with many variables involved can be written as
\bq\label{7b}
x_i\cdot dx_j=qdx_j\cdot x_i+(q^2-1)dx_i\cdot 
x_j\,\,(i<j);\,\,x_i\cdot dx_i=q^2
dx_i\cdot x_i;
\end{equation}
$$x_j\cdot dx_i=qdx_i\cdot x_j\,\,(i<j);\,\,dx_i\wg dx_j=-q^{-1}dx_j\wg
dx_i\,\,(i<j);\,\,dx_i\wg dx_i=0$$
\bq\label{8b}
x_ix_j=qx_jx_i\,\,(i<j);\,\,\pp_i\pp_j=q^{-1}\pp_j\pp_i\,\,(i<j);\,\,\pp_ix_j=qx_j\pp_i
\,\,(i\ne j);
\end{equation}
$$\pp_ix_i-q^2x_i\pp_i=1+(q^2-1)\sum_{j>i} x_j\pp_j$$
Another FODC $\gG_{-}$ arises via $q\to q^{-1}$.
$\hfill\bs$
\\[3mm]\indent
{\bf REMARK 3.4.}
In \cite{cxx,c1,c7} we used the isomorphism $\gs$ above to pass 
between standard PDE
and qPDE.  Let us spell out the correspondence more carefully.  Thus e.g.
$\gs:\,\sum a_{nm}x^np^m\to \sum a_{nm}x_c^np_c^m$.  Following \eqref{5g} with
$D=\pp_x$ we get
\bq\label{30}
\hat{\pp}_x:\,\,\sum a_{nm}x^n_cp^m_c\stackrel{\gs^{-1}}{\to}\sum a_{nm}x^np^m
\stackrel{\pp_x}{\to}\sum a_{nm}[[n]]_{q^2}x^{n-1}p^m\stackrel{\gs}{\to}
\sum a_{nm}[[n]]_{q^2}x_c^{n-1}p_c^m
\end{equation}
This says $\hat{\pp}_x(x_c^np_c^m)=\gs(\pp_xx^np^m)=[[n]]_{q^2}x_c^{n-1}p_c^m$.
Similarly one gets
a similar equation of the form $\hat{\pp}_p
(x_c^np_c^m)=\gs(\pp_px^np^m)=q^nx_c^n[[m]]_{q^2}p_c^{m-1}$.
Since e.g. $(\pp/\pp x)x^n=nx^{n-1}$ it seems awkward however to write out 
$\hat{\pp}_x$ in
terms of $\pp/\pp x$.  One (awkward) solution would be to define an operator
$Q((\pp/\pp x)x_c^n)=[[n]]_{q^2}x_c^{n-1}$, with $Q((\pp/\pp
p)p_c^m)=[[m]]_{q^2}p^{m-1}$ and write
\bq\label{31}
\hat{\pp}_x\sum a_{nm}x_c^np_c^m=\sum a_{nm}Q[(\pp/\pp x)x_c^n]p_c^m;
\end{equation}
$$\hat{\pp}_p\sum a_{nm}x_c^np_c^m=\sum a_{nm}(T_qx^n)Q[(\pp/\pp p)p_c^m]$$
This seems rather bizarre, but appears to be consistent, and could be written
${\bf (A30)}\,\,\hat{\pp}_x=Q(\pp/\pp x)$ and 
$\hat{\pp}_p=T^x_qQ(\pp/\pp p)$.  We
emphasize that Q goes with the derivation and acts only on constants 
produced via
differentiation (i.e. $Q((\pp/\pp x)ax^n)=aQ(\pp/\pp x)x^n)$.
Therefore for $D=\sum \ga_{k\ell}\pp_x^k\pp_p^{\ell}$ one has
${\bf (A31)}\,\,\hat{D}=\sum \ga_{k\ell}[Q(\pp/\pp 
x)^k)]T^x_q[Q(\pp/\pp p)^{\ell}]$.
$\hfill\bs$
\\[3mm]\indent
Now one way of cooking up q-analogues of a classical PDE is to simply replace
$\pp/\pp x^i$ by $\pp_i$ as in $\gG_{+}$ or by $D_q$ more generally. 
Then as $q\to 1$
one recovers $\pp/\pp x^i$ and this makes the extension 
``reasonable".  However it
may not take into accound the basic nature of the original equation 
from a symmetry
point of view. Thus one would be more pleased to produce a q-version 
from q-group
or q-algebra theoretic origins which would reduce to the classical 
operator when
$q\to 1$.  Both approaches have of course been tried and we list some candidate
qPDE from various points of view.
First
one considers wave equations ${\bf (A32)}\,\,[(D_t^{+})^2-{\mc
D}_1^{-}{\mc D}_2^{-}]\phi=0$ where $D_t^{+}\sim t^{-1}(1-T_q)$ and ${\mc
D}_x^{-1}\sim x^{-1} (1-T_q^{-2})$ corresponds to the classical equation
$(\pp_t^2-4\pp_1\pp_2)\phi=0$ when
$t\to (1-q)t$ and $x_i\to (1-q^{-1})x_i/2$ with $q\to 1^{-}$. 
Solutions to {\bf
(A32)} in q-exponentials have the form ($\ga q^2=\gb\gag$)
\bq\label{q6}
\phi(t,x_1,x_2,\ga,\gb,\gag)=\tl{e}_q(\ga t)\tl{E}_q(\gb x_1)\tl{E}_q(\gag x_2)
\end{equation}
from which one can determine symmetry operators (cf. \cite{f2,f3}).
For ${\bf (A33)}\,\,[D^{-}_t-D_1^{+}D_2^{+}]\phi=0$ one has solutions
$\phi=\tl{E}_q(\ga t)\tl{e}_q(\gb x_1)\tl{e}_q(\gag x_2)$ with 
$\ga+q\gb\gag=0$.  We
recall here  that $D^{+}_z\tl{e}_q(\ga z)=\ga \tl{e}_q(\ga z)$ and
$D^{-}_z\tl{e}_q(\gb z)=-\gb q^{-1}\tl{e}_q(\gb z)$ with similar formulas for
$\tl{E}_q$.  For the Helmholz equation
${\bf (A34)}\,\,[D_1^{+}D_2^{+}-\go^2]\phi(x_1,x_2)=0$ solutions can 
be written in
terms of little q-exponentials via ($\ga\gb=\go^2$)
\bq\label{q7}
\phi(x_1,x_2,\ga,\gb)=\tl{e}_q(\ga x_1)\tl{e}_q(\gb x_2)
\end{equation}
For the heat equation in $x,t$ of the form ${\bf (A35)}\,\,[{\mc
D}_t^{-}-(D_x^{+})^2]\phi=0$ there will be some solutions ${\bf
(A36)}\,\,\phi(t,x,\ga,\gb)=\tl{E}_{q^2}(\ga t)\tl{e}_q(\gb x)$ with
$\ga+q^2\gb^2=0\,\,(\ga,\gb\in {\bf C}$).  One arrives at solutions 
to all these
equations be separating variables according to symmetry operators and their
eigenfunctions and this leads for the heat equation also to solutions
\bq\label{q8}
\phi_n(t,x)=q^{-n(n-3)/2}t^{n/2}H_n\left(\frac{x}{q\sqrt{t}};q\right)
\end{equation}
where $H_n\sim$ discrete q-Hermite polynomial.
\\[3mm]\indent
Generally one can phrase a theorem based on the $D\leftrightarrow\hat{D}$
correspondence as follows.  Let $a_{\ga}\in{\bf C}$ with ${\bf (A37)}\,\,
D=\sum a_{\ga}\pp^{\ga}=\sum a_{\ga}\pp_q^{\ga_1}\cdots\pp_n^{\ga_n}$ where the
$\pp_i\in\gG_{+}$.  Then $D\leftrightarrow\hat{D}$ gives a standard partial
differential operator (PDO) with constant coefficients of the type 
indicated in Remark
3.4.
\begin{theorem}
Let $\hat{D}$ be the classical operator corresponding to a D as in {\bf (A37)}
for some fixed ordering of the $x^i$ and assume the equation $\hat{D}\phi=0$
with say real analytic data has a unique real analytic solution. 
Then $D\psi=0$
has a unique solution.
\end{theorem}
\indent
{\bf REMARK 3.5.}
The above proposition only applies to equations with $\hat{D}$ determined 
as in Remarks
3.2-3.4 so the corresponding D may look rather funny.  To clarify this consider
some low order operators
\bq\label{32}
\hat{\pp}_x=Q(\pp/\pp x);\,\,\hat{\pp}_p=T_q^xQ(\pp/\pp 
p);\,\,\hat{\pp}_x\hat{\pp}_p=
\widehat{\pp_x\pp_p}=Q(\pp/\pp x)T^x_qQ(\pp/\pp p);
\end{equation}
$$\hat{\pp}_x^2=Q(\pp/\pp x)Q(\pp/\pp 
x);\,\,\hat{\pp}_p^2=T_q^xQ(\pp/\pp p)T_q^x
Q(\pp/\pp x)$$
This means that a q-wave type equation ${\bf 
(A38)}\,\,(\pp_x^2-\pp_p^2)\psi=0$ would
correspond to $\psi=\gs^{-1}\phi$ (recall 
$\hat{D}\phi=\gs(D(\gs^{-1}\phi))$ and set
$\phi =\sum_0^{\infty}a_{nm}x^np^m$
\bq\label{33}
(\hat{\pp}_x^2-\hat{\pp}_p^2)\phi=0\Rightarrow
\end{equation}
$$\sum_{n=2,p=0}a_{nm}[[n]]_{q^2}[[n-1]]_{q^2}x_c^{n-2}p_c^m-
\sum_{n=0,p=2}a_{nm}q^{2n}[[m]]_{q^2}[[m-1]]_{q^2}x_c^np_c^{m-2}=0$$
If one posits a Cauchy problem with say $\psi(x,0)=h(x)=\sum b_nx^n$ and $\pp_p
\psi(x,0)=0$ this translates to $\phi(x_c,0)=h(x_c)=\sum b_n x_c^n$ 
and $(\pp/\pp p)
\phi(x_c,0)=0$.  One could write out the coefficients but the 
existence of a unique
solution is not apriori known.  In this respect it is more instructive perhaps to
consider $\hat{\phi}=\gs^{-1}\phi$ with $\hat{D}\sim\pp_x^2-\pp_y^2$ so 
\bq\label{eq}
D\hat{\phi}:\,\,x^np^m\stackrel{\gs}{\to}x_c^np_c^m\stackrel{\hat{D}}{\to}
n(n-1)x_cx^{n-2}p_c^m-m(m-1)x_c^mp_c^{m-2}\stackrel{\gs^{-1}}{\to}
\end{equation}
$$\stackrel{\gs^{-1}}{\to}n(n-1)x^{n-2}p^m-(m-1)x^np^{m-2}$$
Then
\bq\label{eqq}
D\sim\frac{n(n-1)}{[[n]]_q[[n-1]]_q}(\pp_q^x)^2-\frac{m(m-1)}{[[m]]_q[[m-1]]_q}(\pp_q^p)^2
\end{equation}
acting on $x^np^m$, and the existence of unique solutions is well known.$\hfill\bs$
\\[3mm]\indent
Regarding symmetries (cf. \cite{f2,f3})
consider a wave equation in light cone coordinates ${\bf (A39)}\,\,
\pp_1\pp_2\phi=0$.  This has an infinite dimensional symmetry algebra 
generated by $v^0_m=
x_1^m\pp_1$ and $w^0_m=x_2^m\pp_2$ for ($m\in{\bf Z}$).  In fact the whole
$W_{1+\infty}\oplus W_{1+\infty}$ algebra generated by 
$v_m^k=x_1^m\pp_1^{k+1}$ and
$w_m^k=x_2^m\pp_2^{k+1}$ for $k\in{\bf Z}_{+}$ maps solutions of {\bf 
(A39)} into
solutions (definition of a symmetry) and $W_{1+\infty}$ without 
center corresponds to
$U({\mc E}(2))$.  For the q-difference version ${\bf
(A40)}\,\,D_1^{+}D_2^{+}\phi(x_1,x_2)=0$ the elements 
$V_m^k=x_1^m(D_1^{+})^{k+1}$
and
$W_m^k=x_2^m(D_2^{+})^{k+1}$ map solutions into solutions and each 
set $V_m^k$ or $W_m^k$
generates a q-deformation of $W_{1+\infty}$. However for the equation ${\bf
(A41)}\,\,[(D_t^{+})^2-(D_x^{+})^2]\phi(t,x)=0$ the situation is 
quite different.
There is still an infinite set of symmetry operators involving polynomials or
arbitrary degree in $t$ and $x$ times powers of $D_t^{+}$ and
$D^{+}_x$ but a general expression seems elusive.  This is in 
contrast to solving wave
equations $(\pp_t^2-\pp_x^2)\phi=0$ where there is conformal 
invariance in $t+x$ and
$t-x$.  One notes that $(t,x)\to (t+x,t-x)$ does not preserve the exponential
2-dimensional lattice and light cone coordinates seem more 
appropriate for q-difference
wave equations.
\\[3mm]\indent
In any case the question of what is a proper canonical q-wave equation seems to have
many answers.  We have choices
\bq\label{34}
[(D_q^t)^2-(D_q^x)^2]\phi=0\,\,(standard\,\,D_q);\,\,[(\pp_q^t)^2-(\pp_q^x)^2]\phi=0
\,\,(\pp_q\sim \pp_i\,\,in\,\,\gG_{+});
\end{equation}
$$D_q^xD_q^t=0\equiv 
D_x^{+}D_t^{+}\phi=0;\,\,\pp_q^t\pp_q^x\phi=0\,\,(\pp_q\in\gG_{+});$$
One could also take q-derivatives in one variable and ordinary 
derivatives in the other.
Now in the special function situations we see for example in 
\eqref{11} and \eqref{18}
that very special forms of second order equations are needed to 
generate q-hypergeometric or
q-Bessel functions.
Another source of information about (nonlinear) PDE arises 
theoretically via the domain of
q-integrable systems such as qKP or qKdV (cf. here in particular 
\cite{cxx,c4,c5}).  We say
theoretically since the actual form of the equations is not simple at 
all.  Thus one recalls
KP in the form ${\bf (A42)}\,\,\pp_tu=(1/4)\pp^3u+3u\pp 
u+(3/4)\pp^{-1}\pp_2^2u$ where
$\pp\sim \pp_x$ and KdV as say ${\bf (A43)}\,\,\pp_tu-6u\pp 
u+\pp^3u=0$.  One can give
relatively simple Hirota type characterizations for qKP and qKdV but 
q-equations of the form
{\bf (A42)}-{\bf (A43)} are complicated  (see here \cite{c4} and 
Section 9 below).  We 
will write out qKP and
qKdV equations later but for now we indicate another source of qPDE 
via zero curvature
conditions on a quantum plane (cf.\cite{c4}).  Thus
\begin{example}
Consider the q-plane situation
\bq\label{192}
(dx)^2=(dt)^2=0;\,\,xt-qtx=0;\,\,dxdt=-q^{-1}dtdx;\,\,xdx=q^2dxx;
\end{equation}
$$xdt=qdtx+(q^2-1)dxt;\,\,tdx=qdxt;\,\,tdt=q^2dtt$$
Then
\bq\label{193}
df=D_t^{-1}\pp_q^xfdx+\pp_q^tfdt
\end{equation}
(with $\pp_q\sim D_{q^{-2}}$ as in $\gG_{+}$ - see below).  Then $A=udx+wdt$
yields
\bq\label{194}
dA+A^2=0\leadsto -q\pp_q^tu+D_t^{-1}\pp_q^xw+uD_x^{-2}D_t^{-2}w-
\end{equation}
$$-\frac{wqt}{x}[D_x^{-1}D_t^{-1}w-D_x^{-3}D_t^{-1}w]-qwD_x^{-1}D_t^{-2}u=0$$
For $q\to q^{-1}$ evidently $\pp_qf=[f(q^{-2}x)-f(x)]/(q^{-2}-1)x\to
[f(q^2x)-f(x)]/(q^2-1)x$ and we write this latter as $\hat{\pp}_q$.  This gives
\bq\label{195}
-q^{-1}\hat{\pp}_q^tu+D_t\hat{\pp}_q^xw+uD_x^2D_t^2w-q^{-1}wD_xD_t^2u-
\frac{wt}{qx}[D_xD_tw-D_x^3D_tw]=0
\end{equation}
Note that as $q\to 1$
\eqref{195} becomes
\bq\label{196}
-\pp_tu+\pp_xw+w-wu=0
\end{equation}
Taking $w=u_x$ one has
\bq\label{197}
-u_t+u_{xx}+u_x-uu_x=0\leadsto u_t+uu_x+u_{xx}-u_x=0
\end{equation}
which is a kind of perturbed Burger's equation with perturbation $-u_x$.
One should examine here more closely the zero curvature condition
$dA+A^2=0$ and its meaning (cf. Sections 7-8).
$\hfill\bs$
\end{example}

\section{EXCURSION TO QKP AND HIROTA EQUATIONS}
\renewcommand{\theequation}{4.\arabic{equation}}
\setcounter{equation}{0}

Now referring to \cite{c4} for some details, we give a somewhat expanded version of
various features developed there.  Thus one has for qKP a Lax 
operator ${\bf (A44)}\,\,L=
D_q+a_0+\sum_1^{\infty}a_iD_q^{-i}$ with equations ${\bf 
(A45)}\,\,\pp_jL=[L^j_{+},L]$
where $\pp_j\sim\pp/\pp t_j$ for $t=(t_1,t_2,\cdots)$.
There is a gauge operator $S=1+\sum\tl{w}_jD_q^{-j}$ satisfying $L=SD_qS^{-1}$.
We use now $e_q,\,E_q$ as in \eqref{1} and {\bf (A3)} and note for 
$\xi=\sum_1^{\infty}t_iz^i$ that the
wave functions are given via ${\bf (A46)}\,\,\psi_q=Se_q(xz)exp(\xi)$ 
and $\psi^*_q=
(S^*)^{-1}_{x/q}exp_{1/q}(-xz)exp(-\xi)$.  Then there is a tau function
$\tau_q(t)=\tau(t+c(x))$  where
\bq\label{35}
c(x)=[x]_q=\left(\frac{(1-q)x}{(1-q)},\frac{(1-q)^2x^2}{2(1-q^2)},
\frac{(1-q)^3x^3}{3(1-q^3)},\cdots\right)
\end{equation}
leading to
\bq\label{36}
\psi_q=\frac{\tau_q(x,t-[z^{-1}])}{\tau_q(x,t)}exp_q(xz)exp(\xi);\,\,\psi_q^*=
\frac{\tau_q(x,t+[z^{-1}])}{\tau_q(x,t)}exp_{\frac{1}{q}}(-xz)exp(-\xi)
\end{equation}
where $[z]=(z,z^2/2,z^3/3,\cdots)$.  Vertex operators for KP are defined via
$$X(t,z)=exp(\xi)exp(-\sum_1^{\infty}(\pp_i/i)z^{-i})=exp(\xi(t,z)exp(-\xi(
\tl{\pp},z^{-1}))$$
The Schur polynomials are defined via
\bq\label{37}
\sum_{{\bf 
Z}}\tl{p}_k(x,t_1,\cdots)z^k=exp_q(xz)exp(\sum_1^{\infty}t_kz^k);\,\,
\tl{p}_k(x,t)=p_k(t+c(x));
\end{equation}
$$p_n(y)=\sum\left(\frac{y_1^{k_1}}{k_1!}\right)\left(\frac{y_2^{k_2}}
{k_2!}\right)\cdots;\,\,\,\sum_0^{\infty}p_n(y)z^n=exp(\xi(y,z))$$
It is unrealistic to compute the actual qKP equation analogous to 
{\bf (A42)} from say
{\bf (A45)} and one resorts to Hirota type equations for such 
calculations (cf. \cite{c4}
for a variety of relevant Hirota equations).  We look first at Hirota 
equations for
$f=D^n\tau_q\,\,(n\geq 1)$ in the form ($Dh(x)=T_qh(x)=h(qx)$)
\bq\label{38}
\pp_1\pp_rf\cdot f=2p_{r+1}(\tl{\pp})f\cdot f=2\sum_0^{r+1}p_j(-\tl{\pp})
fp_{r+s+1-j}(\tl{\pp})f
\end{equation}
with standard Schur polynomials.  We indicate a few calculations; 
first for $n=1,\,r=1$
\bq\label{39}
\pp_1\pp_1log(f)=\frac{1}{f^2}\sum_0^2p_j(-\tl{\pp})fp_{2-j}(\tl{\pp})f=\pp_1^2log(f)
\end{equation}
This is tautological.
We recall also $D=1+(q-1)xD_q$ and $D^2=1+(q^2-1)xD_q+q(q-1)^2x^2D_q^2$ with
\bq\label{40}
D^3=(q-1)^3q^3x^3D_q^3+q^2x^2(q-1)(q^2-1)D_q^2+(q^3-1)xD_q+1
\end{equation}
Then for $n=1,\,r=2$ one has ${\bf (A47)}\,\,
\pp_1\pp_2f\cdot f=\sum_0^2p_j(-\tl{\pp})fp_{3-j}(\tl{\pp})f$.
Better, we want $D_q^3$ in the act so try e.g. $n=3,\,r=2$.  The case $r=3,\,
n=3$ gives ${\bf (A48)}\,\,\pp_1\pp_3f\cdot f=2\sum_0^4p_j(-\tl{\pp})fp_{4-j}
(\tl{\pp})f$.  This leads to a reasonable enough expression involving 
products of
derivatives of f but it is hard to recognize its relation to KP. 
Therefore we go to a
variation of the Hirota formulas developed in \cite{c4} which leads 
directly to expressions
in terms of $log\tau$ (or $log f$ here).  We simply state the result 
here and refer to
\cite{c4} for the proof.  
\begin{theorem}
Using the differential Fay identity 
one shows (following
\cite{c9}) that the classical KP Hirota equations can be written 
directly in terms of
$u=\pp^2log\tau$ via an equation ${\bf 
(A49)}\,\,p_i(Z_j(u))=\sum_{n+\ell=\ell-1}\pp_1^{-1}\pp_nu$.
Here ${\bf 
(A50)}\,\,p_i(Z_j)=\sum_{n+\ell=i-1}p_{\ell}(-\tl{\pp})\pp_n\pp_1log\tau$ 
and $Z_j
=\sum_{k+m=j}\pp_1^2\tl{F}_{km}$ where 
$\tl{F}_{km}=p_k(-\tl{\pp})p_m(-\tl{\pp})u$.
\end{theorem}
\indent
A simple minded version of this can be phrased in a different form via
\bq\label{41}
(\sum_1^{\infty} p_k(-\tl{\pp})p_m(-\tl{\pp})\gl^{-k-m}\pp_1log\tau)(1+
\sum_1^{\infty}\gl^{-n-1}\sum_0^{\infty}p_{\ell}(-\tl{\pp})\gl^{-\ell}
\pp_n\pp_1log\tau)=
\end{equation}
$$=\sum_1^{\infty}\gl^{-n-1}\sum_0^{\infty}p_{\ell}(-\tl{\pp})\gl^{-\ell}
\pp_n\pp_1^2log\tau$$
Then equate powers of $\gl$ and write ($f=\pp_1log\tau$)
\bq\label{42}
(\sum_1^{\infty}p_{km}\gl^{-k-m}f)(1+\sum_{n=1}^{\infty}\sum_{\ell=0}^{\infty}
p_{\ell}\pp_n\gl^{-n-\ell-1}f)=\sum\sum p_{\ell}\pp_n\gl^{-n-\ell-1}\pp_1f
\end{equation}
\bq\label{43}
\sum_{k,m=1}^{\infty}p_{km}\gl^{-k-m}f+\sum_{k,m=1}^{\infty}\sum_{n=1}^{\infty}
\sum_{\ell=0}^{\infty}p_{k,m}fp_{\ell}\pp_nf\gl^{-k-m-\ell-n-1}=
\end{equation}
$$=\sum_{n=1}^{\infty}\sum_{\ell=0}^{\infty} 
p_{\ell}\pp_n\gl^{-n-\ell-1}\pp_1f$$
\bq\label{44}
\gl^{-2}:\,\,p_{11}f=\pp_1^2f;\,\,\gl^{-3}:\,\,2p_{12}f=\pp_2\pp_1f+p_1\pp_1^2f;
\end{equation}
$$\gl^{-4}:\,\,(2p_{13}+p_{22})f+p_{11}f\pp_1f=\pp_3\pp_1f+p_1\pp_2\pp_1f+
p_2\pp_1^2f$$
We can adapt this to \eqref{38} since e.g. \eqref{41}-\eqref{44} 
applies to situations with
Hirota formula $\pp_1\pp_r\tau\cdot\tau=2p_{r+1}(\tl{\pp})\tau\cdot\tau$ via
\bq\label{45}
\oint \psi^*(t,\gl)\psi(t',\gl)d\gl=0\leadsto \oint \tau(t+[\gl^{-1}])
\tau(t'-[\gl^{-1}])e^{\sum (t'_i-t_i)\gl^i}d\gl=0
\end{equation}
Let $t\to t+y,\,\,t'\to t-y$ to obtain
\bq\label{46}
0=\oint \tau(t+y+[\gl^{-1}])\tau(t-y-[\gl^{-1}])e^{-2\sum y_i\gl^i}d\gl=
\end{equation}
$$=\oint e^{\sum y_i\pp_i+\sum\tl{\pp}_i\gl^{-i}}\tau\cdot\tau e^{-2\sum
y_i\gl^i}d\gl
=\oint \sum p_n(-2y)\gl^n\sum p_{\ell}(\tl{\pp})\gl^{-\ell}e^{\sum y_i\pp_i}
\tau\cdot\tau d\gl$$
The Hirota formula for \eqref{38} arise from
\bq\label{47}
0=\oint e^{\xi(t,z)}\sum p_{\ell}(-\tl{\pp})z^{-\ell}D^n\tau_q
e^{-\xi(t',z)}\sum p_k(\tl{\pp})z^{-k}D^n\tau_q
\end{equation}
Then for the power $\gl^{-4}$ in \eqref{44} one gets
\bq\label{48}
-2\pp(-(1/3)\pp_3+(1/2)\pp\pp_2-(1/6)\pp^3)f+(1/4)(\pp^2-\pp_2)^2f+\pp^2f\pp f=
\end{equation}
$$\pp\pp_3f-\pp^2\pp_2f+(1/2)(\pp^2-\pp_2)\pp^2f$$
leading to
\bq\label{49}
\pp\pp_3f=(1/4)\pp^4f+(3/4)\pp_2^2f+3\pp^2f\pp f
\end{equation}
For $f=\pp log\tau$ we recognize this as a KP equation ${\bf (A51)}\,\,
\pp_3u=(1/4)\pp^3u+(3/4)\pp_2^2\pp^{-1}u+3u\pp u$.  For 
$f=\pp_1D^3log\tau_q$ we think of
$\pp\sim\pp_1$ and specify e.g. $\tl{u}=\pp_1^2D^3log\tau_q$.  Then
{\bf (A51)} becomes
\bq\label{50}
\pp_3\tl{u}=\frac{1}{4}\pp_1^3\tl{u}+\frac{3}{4}\pp_2^2\pp_1^{-1}\tl{u}+
3(\pp_1\tl{u})\tl{u}
\end{equation}
This says that $\tl{u}$ (and in fact any 
$\hat{u}=\pp_1^2D^nlog\tau_q$) satisfies a standard
KP equation in the t variables.  Here $t_1+x$ is the first argument of 
$\tau_q=\tau(t+c(x))$ and $x$ plays no role; 
not surprisingly
$\pp_1^2\tau(t+D^nc(x))$ satisfies KP in the t variables (note from 
\cite{a2,a3}
\\[3mm]\indent
Note now from \cite{a2,a3} that $D^nc(x)=c(x)-\sum_1^n[\gl_i^{-1}]$ 
where $D\gl_n=\gl_{n+1}$
and $\gl_i^{-1}=(1-q)xq^{i-1}$.  Generally from \cite{a2,a3} again the
operators
\bq\label{51}
Q=D+\sum_{-\infty<i\leq 0}a_i(x)D^i;\,\,Q_q=D_q+\sum_{-\infty<i\leq 
0}b_i(x)D_q^i
\end{equation}
satisfy ($\pi(k)=\pm$ for $k$ even or odd)
\bq\label{52}
a_{\ell}(x)=poly.\,\,in\,\,\left\{\begin{array}{cc}
\frac{\pp^k}{\pp t_{i_1}\cdots\pp 
t_{i_k}}log[\tau(c(x)+t)^{\pi(k)}D^{\ell+1}\tau(c(x)+t)]
& (k\geq 2)\\
\sum_1^{\ell+1}\gl_i^j(x)+\pp_jlog\frac{D^{\ell+1}\tau(c(x)+t)}{\tau(c(x)+t)} 
& (k=1)
\end{array}\right.
\end{equation}
Further
\bq\label{53}
a_i(x)=\sum_{0\leq k\leq n-i}\frac{\left[\begin{array}{c}
k+i\\
k
\end{array}\right]b_{k+i}(x)}{(-x(q-1)q^i)^k}
\end{equation}
Consider now the Hirota formula from \cite{a2,c4} ($n>m$)
\bq\label{54}
0=\oint D^ne_q(xz)e^{\xi(t,z)}\sum
p_{\ell}(-\tl{\pp})z^{-\ell}\tau(t+c(x))\cdot
\end{equation}
$$D^{m+1}e_{1/q}(-xz)e^{-\xi(t',z)}\sum
p_k(\tl{\pp})z^{-k}\tau(t'+c(x))dz$$
and recall that we need a way to handle the expression 
$D^ne_q(xz)D^{m+1}e_{1/q}(-xz)$
for $n=m+1+s$ with $s>0$.  One can write now $D^n=D^{m+1}D^s$ and use 
formulas for $D^s$ in
terms of $D_q$.  For small s these were worked out in \cite{c4} and 
we have e.g.
${\bf (A52)}\,\,D=1+(q-1)xD_q$ with
$(\clubsuit)\,\,
qxD_qf=
D_qx=qxD_q+1$.
\\[3mm]\indent
{\bf REMARK 4.1.}
A nice general formula for $D^n$ in terms of $D_q^m$ was indicated by 
J. Cigler (cf. \cite{cz}
for background).  Thus start from 
$x^n=\sum_0^n{\nm}_q(x-1)(x-q)\cdots (x-q^{m-1})$ to get
$D^n=\sum_0^n{\nm}_q(D-1)(D-q)\cdots (D-q^{m-1})$ (note $x^m(x^{-1};q)_m=(x-1)
(x-q)\cdots (x-q^{m-1})$).  Then since $(D-1)(D-q)\cdots 
D-q^{m-1})=(q-1)^mq^{m(m-1)/2}x^m
D_q^m$ one obtains the formula
\bq\label{alpha}
D^n=\sum_0^n{\nm}_q(q-1)^mx^mq^{m(m-1)/2}D_q^m
\end{equation}
We will rewrite the formulas below involving $D^s$ in terms of 
\eqref{alpha} later.$\hfill\bs$
\\[3mm]\indent
Using brute force (or \eqref{alpha}) for a few low order situations one has
${\bf (A53)}\,\,
D^2=1+(q^2-1)xD_q+q(q-1)^2x^2D_q^2$.
We note also
$$\nm_q=\frac{(q^n-1)\cdots(q^{n-m+1}-1)}{(q^m-1)\cdots(q-1)}
=\frac{[n]_q}{[m]_q[n-m]_q}$$
Similarly 
\bq\label{58}
D^3=(q-1)^3q^3x^3D_q^3+qx^2(q-1)(q^3-1)D_q^2+(q^3-1)xD_q+1
\end{equation}
Then e.g. for $s=1$ write ${\bf 
(A54)}\,\,De_q=e_q+(q-1)xze_q]$ so
$D^{m+1}De_q=D^{m+1}e_q +(q-1)q^{m+1}xzD^{m+1}e_q$.  Since
$D^{m+1}e_q(xz)D^{m+1}e_{1/q}(-xz)=1$ the integrand in \eqref{54} 
picks up a factor
$[1+(q-1)q^{m+1}xz]$.  Next for $s=2$ one has
\bq\label{59}
D^{m+1}D^2e_q=[(1+(q^2-1)q^{m+1}xz+q^{2m+3}x^2z^2(q-1)^2]D^{m+1}e_q
\end{equation}
Finally for $s=3$ we write
\bq\label{60}
D^3e_q=
[1+(q^3-1)xz+qx^2(q^3-1)z^2+(q-1)^3q^3x^3z^3]e_q;\,\,D^{m+1}D^3e_q=
\end{equation}
$$[1+(q^3-1)q^{m+1}xz+q^{2m+3}(q^3-1)x^2z^2+(q-1)^3q^{3m+6}x^3z^3]D^{m+1}e_q$$
Generally one has then ${\bf (A55)}\,\,D^se_q=\sum_0^sf_j(q,x)z^je_q$
(with the $f_j\sim f_{js}$ given via Remark 4.1) so that
\bq\label{61}
D^{m+1}D^se_q=\sum_0^sf_j(q,q^{m+1}x)z^jD^{m+1}e_q
\end{equation}
This leads to a possibly useful result for small s.
First,
using \eqref{61} with $n=m+1+s$, the Hirota formula \eqref{54} leads to
\bq\label{62}
0=\oint e^{-2\sum_1^{\infty} 
y_iz^i}\sum_0^sf_j(q,q^{m+1}x)z^je^{\sum_1^{\infty}
y_i\pp_i}\sum_0^{\infty}p_{\ell}(\tl{\pp})z^{-\ell}F_q\cdot F_q
\end{equation}
where $F_q(t)=D^{m+1}\tau_q(t)$.  Thus one has residue terms based on
\bq\label{63}
\oint \sum_0^{\infty}p_i(-2y)z^i\sum_0^sf_j^sz^j(1+\sum y_m\pp_m+O(y_iy_j))
\sum_0^{\infty}p_{\ell}(\tl{\pp})z^{-\ell}F_q\cdot F_q=0
\end{equation}
leading to the terms where $i+j-\ell=-1$, namely
\bq\label{64}
\sum_ip_i(-2y)\sum_{j=0}^sf_j^s(q,q^{m+1}x)p_{i+j+1}(\tl{\pp})(1+\sum 
y_m\pp_m)F_q\cdot
F_q=0
\end{equation}
We look for the coefficient of $y_r$ via
\bq\label{65}
-2\sum_{j=0}^sf_j^sp_{r+j+1}(\tl{\pp})F_q\cdot 
F_q+\sum_{j=0}^sf_j^sp_{j+1}(\tl{\pp})
\pp_rF_q\cdot F_q=0
\end{equation}
\begin{theorem}
Hirota equations for qKP with $n=m+1+s$ in \eqref{54} are given via 
\eqref{65} with the
$f_j^s$ for $s=1,2,3$ indicated in {\bf (A54)}, \eqref{59}, and \eqref{60}.
\end{theorem}
\indent
Then for $s=0$ we get ($f_0^0=1$)
\bq\label{66}
-2p_{r+1}(\tl{\pp})F_q\cdot F_q+\pp_1\pp_rF_q\cdot F_q=0
\end{equation}
as in \eqref{38}.  For completeness we write out a few formulas. 
Thus ${\bf (A56)}\,\,
s=1\sim \sum_0^1f_j^1z^j=1+(q-1)xz;\,\,s=2\sim 
\sum_0^2f_j^2z^2=1+(q^2-1)xz+qx^2(q-1)^2z^2;$
and $s=3\sim \sum_0^3f_j^3z^j=1+(q^3-1)xz+q^2x^2(q^2-1)z^2+(q-1)^3q^3x^3z^3$.
A typical formula then from \eqref{65} would be
\begin{corollary}
For $s=2$ the Hirota equations become 
($f_0^2=1,\,\,f_1^2=(q^2-1)x,\,\,f_2^2=qx^2(q-1)^2$)
\bq\label{67}
0=-2[p_{r+1}(\tl{\pp})+f_1^2p_{r+2}(\tl{\pp})+f_2^2p_{r+3}(\tl{\pp})]F_q\cdot 
F_q+
\end{equation}
$$+[p_1(\tl{\pp})+f_1^2p_2(\tl{\pp})+f_2^2p_3(\tl{\pp})]\pp_rF_q\cdot F_q$$
\end{corollary}
\indent
These are certainly different from the classical situation as exemplified in
\eqref{66}.  However they are
again equations in the $t_i$ variables and $D_q$ only arises via 
$F_q=D^{m+1}\tau_q$ with
D action on $\tau_q$ via
$D\tau_q=\tau(t+Dc(x))=\tau(t+c(x)-[\gl_1^{-1}])=\tau(t+c(qx))$
where $\gl_n^{-1}=(1-q)xq^{n-1}$ (note e.g. 
$x-\gl_1^{-1}=x-(1-q)x=qx$).  Thus one could
write e.g. ${\bf (A57)}\,\,Dc(x)=c(x-[\gl_1^{-1}])=c(qx)=exp(\sum 
\gl_1^{-i}\pp_i)c(x)$.
Note also ${\bf (A58)}\,\,D^nc(x)=c(x)-\sum_1^n[\gl_i^{-1}(x)]$ and 
we mention in passing
the lovely formula from \cite{a5}
\bq\label{677}
\prod_{k=m+2}^n\left(1-\frac{z}{\gl_k}\right)=\prod_{m+2}^nexp\left(\sum_1^{\infty}
\frac{1}{i}\left(\frac{z}{\gl_k}\right)^i\right)=D^ne_q(xz)D^{m+1}e_q(xz)^{-1}
\end{equation}
to supplement considerations before Theorem 3.2 (note 
$[\gl_1^{-1}]\sim ((1/i)\gl_1^{-i})$).
\\[3mm]\indent
There is some incentive now to look at Hirota formulas for qKP of the 
type mentioned in
\cite{c4,i6} written as
\bq\label{68}
Res_z D_q^n\pp^{\ga}\psi_q\psi^*_q=0;\,\,\ga=(\ga_1,\ga_2,\cdots)
\end{equation}
It is not specified whether or not $t=t'$ in $\psi_q$ and $\psi_q^*$ 
respectively and we
will consider this below.  The Dickey lemma used here is
\bq\label{69}
Res_z(Pe_q(xz)Q^*_{x/q}e_{1/q}(-xz))=Res_{D_q}(PQ)
\end{equation}
proved as follows.  Let $P=\sum p_iD_q^i$ and $Q^*=\sum q_jD^j_{1/q}$.  Then
$$Res_z[Pe_q(xz)Q^*_{x/q}e_{1/q}(-xz)]=Res[\sum p_iz^ie_q\sum 
q_j(x/q)(-q)^jz^je_{1/q}]=$$
\bq\label{70}
=\sum_{i+j=-1}(-q)^jp_i(x)q_j(x/q)
\end{equation}
Now $Q^*=\sum q_jD^j_{1/q}$ and one knows that $V=\sum v_mD_q^m$ has 
adjoint $V^*=
(D_q^*)^jq_j$ where $D_q^*=-(1/q)D_{1/q}$.  This implies that 
$V^*=\sum (-q)^{-j}D_{1/q}
q_j$ and hence $Q^*=\sum q_j(x/q)D^j_{1/q}$ comes from $Q=\sum 
(-q)^jD_q^jq_j(x/q)$ and
\bq\label{71}
Res_{D_q}(PQ)=\sum_{i+j=-1}(-q)^jp_i(x)q_j(x/q)
\end{equation}
as desired.
We note also from \cite{a5} that the Hirota bilinear identity for 
standard KP arises from
\bq\label{72}
\oint \psi(x,t,z)\psi^*(x',t',z)dz=0
\end{equation}
Hence (and this goes back to \cite{d7})
\bq\label{73}
\oint \left(\pp^{\ga}\psi(x,t,z)\right)\psi^*(x',t',z)dz=0
\end{equation}
Subsequently one can send $(x',t')\to (x,t)$ or not.  The same 
reasoning should apply for
qKP, given that e.g.
\bq\label{74}
\psi_q=(Se_q(xz))e^{\xi}=(1+\sum\tl{w}_jz^{-j})e_q(xz)e^{\xi}\frac{\tau_q(t-[z^{-1})}
{\tau_q(t)}
\end{equation}
where $S=1+\sum\tl{w}_jD_q^{-j}$ and $\tau_q(t)=\tau(t+c(x))$.  For 
$n=1=m+1$ one has also
\bq\label{75}
\oint D\psi_q(x,t,z)D\psi^*(x,t',z)dz=0\iff 
\oint\psi_q(x,t,z)\psi_q^*(x,t',z)dz=0\iff
\end{equation}
$$\iff \oint 
e^{\xi(t,z)-\xi(t',z)}\tau(t+c(x)-[z^{-1}])\tau(t'+c(x)+[z^{-1}])dz=0$$
and the latter equation is equivalent to the standard KP formula as 
in \eqref{72} but with
the stipulation $x=x'$ in \eqref{72}.  Now write 
$e_q(xz)=exp\sum_1^{\infty}c_i(x)z^i$ so
that $e_q(xz)exp(\xi)=exp[\sum_1^{\infty}(t_i+c_i(x))z^i]$ and note 
that $t_i+c_i(x)\to t_i$
and $t'_i+c_i(x)\to t_i'$ puts the second equation in \eqref{75} into 
the same pattern as
\eqref{72}, with new variables $t_i,\,t'_i$; subsequently one could set
$t_i=t_i+c(x)$ and $t_i'=t_i'+c_i(x')$, leading to
\bq\label{76}
0=\oint \psi_q(x,t,z)\psi^*_q(x',t',z)dz
\end{equation}
(alternatively one could take $t_i'=t_i''-c_i(x)+c_i(x')$).
This will enable us to state
\begin{theorem}
The Hirota relations \eqref{68} are valid in the form
\bq\label{77}
\oint \left(D^n_q\pp^{\ga}\psi_q(x,t,z)\right)\psi_q^*(x',t',z)dz=0
\end{equation}
\end{theorem}
\indent
For completeness, if nothing else, let's see what this engenders. 
Recall $D_q(fg)=(D_qf)g+
(Df)D_qg$ and $\psi_q=e_q(xz)exp(\xi)\tau(t+c(x)-[z^{-1}])$.  Also 
$Dc(x)=(c_i(x)-(1/i)
\gl_1^{-i})=(c_i(qx))$ and $\tau_q(t)=exp(\sum c_n(x)\pp_n)\tau(t)$. 
Further $D^nc(x)=
c(x)-\sum_1^n[\gl_i^{-1}]$ where again $\gl_n^{-1}=(1-q)xq^{n-1}$. 
Thus $\tau^{\pm}_q\sim
\tau_q(t+c(x)\pm [z^{-1}])$ and
\bq\label{78}
D_q\psi\sim 
e^{\xi}[ze_q(xz)\tau_q^{-}+e_q(qxz)D_q\tau^{-}_q];\,\,\pp_r\psi_q\sim
e_q(xz)[z^re^{\xi}\tau^{-}_q+e^{\xi}\pp_r\tau_q^{-}]
\end{equation}
Play this off against $e_q(xz)^{-1}exp(\xi')\tau(t'+c(x)+[z^{-1}])$ then to get
\bq\label{79}
D_q\psi_q\psi^*_q\sim e^{\xi-\xi'}[z+e_q(qxz)e_q(xz)^{-1}D_q]\tau_q^{-}(t)\tau^{+}(t')
\end{equation}
But then by \eqref{67} (cf. also \eqref{61}) we know 
$(De_q)e_q^{-1}=1-(z/\gl_1)=
1+xz(q-1)$ so ${\bf 
(A59)}\,\,(D_q\psi_q)\psi_q^*=exp(\xi-\xi')[z+(1+xz(q-1))D_q]
\tau^{-}_q(t)\tau_q^{+}(t')$.  Similarly
\bq\label{79}
\pp_r\psi_q=exp(\xi-\xi')(z^r+\pp_r)\tau^{-}_q(t)\tau_q^{+}(t')
\end{equation}
We compute also $D_q^2$ and $D_q^3$ via
\bq\label{80}
D_q^2\psi_q=e^{\xi}[z(ze_q(xz)\tau_q^{-}+e_q(qxz)D\tau_q^{-})+qze_q(qxz)D_q\tau_q^{-}
+e_q(q^2xz)D^2_q\tau_q^{-}]=
\end{equation}
$$=e^{\xi}[z^2e_q\tau_q^{-}+(q+1)ze_q(qxz)D_q\tau_q^{-}+e_q(q^2xz)D^2_q\tau_q^{-}]$$
Then note that $e_q(qxz)e_q(xz)^{-1}=1+xz(q-1)$ again and 
$e_q(q^2xz)e_q(xz)^{-1}=
[1-(z/\gl_1)][1-(z/\gl_2)]=[1-z(1-q)x][1-z(1-q)xq]=1-z(1-q^2)x+z^2(1-q)^2x^2q$
(cf. also \eqref{59}).  Consequently
\bq\label{81}
(D_q^2\psi(t))\psi^*_q(t')=e^{\xi-\xi'}\{[1+xz(q^2-1)+x^2z^2q(1-q)^2](D^2_q\tau^{-}_q(t))
\tau_q^{+}(t')+
\end{equation}
$$+z^2\tau_q^{-}(t)\tau_q^{+}(t')+[(q+1)z+(q^2-1)xz^2](D_q\tau_q^{-}(t))\tau_q^{+}(t')\}$$
One can also write now (recall $\gl_i^{-1}=(1-q)xq^{i-1}$)
\bq\label{82}
\prod_1^n\left(1-\frac{z}{\gl_k}\right)=D^ne_q(xz)e_q(xz)^{-1}=1-z\sum_1^n\frac{1}{\gl_i}
+z^2\sum_{i<j}^n\frac{1}{\gl_i\gl_j}-\cdots+\frac{(-1)^nz^n}{\gl_1\cdots\gl_n}
\end{equation}
This can be written as ${\bf 
(A60)}\,\,D^ne_q(xz)e_q(xz)^{-1}=\sum_0^n(-1)^iz^i\gL_i(q,x)$
(cf. also \eqref{61}).  For $D^3$ we write
\bq\label{83}
D^3\psi_q=
e^{\xi}\{z^2[ze_q+De_q
D_q]+(q+1)z[qze_q+D^2e_qD_q^2]+[q^2ze_q+D^3e_qD^3_q]\}\tau^{-}_q=
\end{equation}
$$=e^{\xi}\{(z^3+(q+1)qz^2+q^2z)e_q+z^2(De_q)D_q+(q+1)z(D^2e_q)D^2_q+(D^3e_q)D^3_q\}
\tau_q^{-}$$
This implies
\bq\label{84}
(D^3\psi_q)\psi^*_q=e^{\xi-\xi'}\{[z^3+(q+1)qz^2+q^2z]\tau_q^{-}(t)\tau^{+}_q(t')+
\end{equation}
$$+z^2(1-z\gL_1)(D_q\tau^{-}_q(t))\tau^{+}_q(t')+(q+1)z(1-z\gL_1+z^2\gL_2)(D^2_q\tau^{-}_q
(t)\tau_q^{+}(t')+$$
$$+(1-z\gL_1+z^2\gL_2-z^3\gL_3)(D^3_q\tau_q^{-}(t))\tau^{+}_q(t')$$
The coefficients for $D_q$ and $D_q^2$ can be expressed in the 
$\gL_i$ via (note for
$\gl_0^{-1}=(1-q)xq^{-1}$ one has $D^n\gl_0^{-1}=\gl_n^{-1}$)
\bq\label{85}
D_q(e_q\tau^{-}_q)=(ze_q+De_qD_q)\tau^{-}_q;\,\,D_q(e_q\tau_q^{-})e_q^{-1}=[z+(De_q)e_q^{-1}]
\tau_1^{-}=
\end{equation}
$$=[z+1-\gl_1z]\tau^{-}_q;\,\,D^2_q(e_q\tau_q^{-})=[(z(ze_q+De_qD_q)+qzDe_qD_q+D^2e_qD_q^2]
\tau_q^{-};$$
$$D^2_q(e_q\tau_q^{-})e_q^{-1}=[z^2+z(1-\gL_1z)D_q+qz(1-\gL_1z)D_q+(1-\gL_1z+\gL_2z^2)
D_q^2]\tau_q^{-}(t)$$
Compare this with \eqref{81} via
\bq\label{86}
z(1+q)(1-\gL_1z)=z(1+q)(1+(q-1)xz)=(q+1)z+(q^2-1)xz^2
\end{equation}
as desired.
\\[3mm]\indent
Now go to \eqref{77} and look at some Hirota formulas in the form
\bq\label{87}
\oint
[D^n_q\pp^{\ga}(e^{\xi}e_q(xz)\tau_q(t-[z^{-1}]))](e^{\xi'}e_q(xz)^{-1}\tau_q(t+[z^{-1}]))=0
\end{equation}
One can bring $\pp^{\ga}$ inside via {\bf (A59)}; for example
\bq\label{88}
\pp_a\pp_b(e^{\xi}\tau_q)=e^{\xi}(z^{a+b}+z^b\pp_a+z^a\pp_b+\pp_a\pp_b)\tau_q
\end{equation}
We will only check a few cases.  Thus for $\pp_r=\pp_3$ and $n=0$ 
with $t\to t-y$ and
$t'\to t+y$
\bq\label{89}
0=\oint e^{\xi-\xi'}[(z^3+\pp_3)\tau_q^{-}(t)]\tau_q^{+}(t')=
\end{equation}
$$=\oint e^{-2\sum y_mz^m}(z^3+\pp_3)(exp(\sum 
(-\tl{\pp}_i)z^{-i})\tau_q(t-y)exp
(\sum \tl{\pp}_jz^{-j})\tau_q(t+y)dz$$
Write this as
\bq\label{90}
0=\oint e^{-2\sum y_mz^m}z^3[e^{\sum y_n\pp_n}\sum 
p_{\ell}(\tl{\pp})z^{-\ell}\tau_q\cdot
\tau_q]dz+
\end{equation}
$$+\oint e^{-\sum2y_iz^i}e^{\sum y_m\pp_m}\sum
p_{\ell}(\tl{\pp})z^{-\ell}\tau_q\cdot\pp_3\tau_q$$
\\[3mm]\indent
{\bf REMARK 4.2.}
We note here an abuse of notation in the bilinear notation. 
Precisely one should write
\bq\label{91}
\pp_m^ja\cdot b=\frac{\pp^m}{\pp s_j^m}a(t_j+s_j)b(t_j-s_j)|_{s=0}
\end{equation}
but in KP practice one is writing e.g.
\bq\label{92}
f(t+y+[z^{-1}])g(t-y-[z^{-1}])=exp(\sum\tl{\pp}_iz^{-i})exp(\sum 
y_m\pp_m)f\cdot g
\end{equation}
and then after a residue calculation one sets the coefficient of say 
$y_r$ equal to zero.
It seems that t is playing the role of s here but it is deliberately 
not being subsequently
set equal to zero.  Thus $f\cdot g$ in KP means \eqref{92} and not \eqref{91}.
$\hfill\bs$
\\[3mm]\indent
For a residue calculation in \eqref{90} we consider first
\bq\label{93}
\oint\sum p_k(-2y)z^k\sum p_{\ell}(\tl{\pp})z^{-\ell}[z^3\tau_q\cdot\tau_q+
\tau_q\cdot \pp_3\tau_q]dz=0
\end{equation}
this seems awkward so it is better look at
\bq\label{94}
0=\oint[\sum p_k(-2y)z^kexp(-\sum y_m\pp_m)z^3\sum 
p_{\ell}(-\tl{\pp})z^{-\ell}\tau_q
exp(\sum y_m\pp_m)\sum p_n(\tl{\pp})z^{-n}\tau_q+
\end{equation}
$$+\sum p_k(-2y)z^kexp(-\sum y_m\pp_m)\pp_3\sum 
p_{\ell}(-\tl{\pp})z^{-\ell}\tau_q
exp(\sum y_m\pp_m)\sum p_n(\tl{\pp})z^{-n}\tau_q]dz$$
Now one can group powers with $k+3-\ell-n=-1$ in the first terms and 
$k-\ell-n=-1$ in the
second.  For the coefficient of $y_r$ say it is best to take 
$p_r(-2y)$ with r fixed and
vary
$\ell$ and $n$; then use $p_0(-2y)=1$ against $\pm\pp_r$ with 
$\ell+n$ varying again.
Thus pick $k=r$ and posit
\bq\label{95}
-2\sum_{\ell+n=r+4}p_{\ell}(-\tl{\pp})\tau_qp_n(\tl{\pp})\tau_q-2\sum_{\ell+n=r+1}
p_{\ell}(-\tl{\pp})\pp_3\tau_qp_n(\tl{\pp})\tau_q+
\end{equation}
$$+\sum_{\ell+n=r+4}p_{\ell}(-\tl{\pp})(-\pp_r)\tau_qp_n(\tl{\pp})\tau_q+
\sum_{\ell+n=r+1}p_{\ell}(-\tl{\pp})\pp_3\tau_qp_n(\tl{\pp})\pp_r\tau_q=0$$
as a Hirota equation.
For a simple Hirota equation involving $D_q$ one can take $\ga=0$ and 
$n=1$ in \eqref{87}.
Then via {\bf (A59)} one has
\bq\label{96}
\oint e^{\xi-\xi'}([z+xz(q-1)D_q]\tau_q^{-}(t))\tau_q^{+}(t')dz=0
\end{equation}
This leads to
\bq\label{97}
\oint e^{-2\sum 
y_mz^m}[\tau_q^{-}(t-y)\tau_q^{+}(t+y)+z([1+x(q-1)D_q]\tau_q^{-}(t-y))
\tau_q^{+}(t+y)]dz=0
\end{equation}
Written out this is
\bq\label{98}
\oint \sum p_k(-2y)z^k\{\sum p_{\ell}(-\tl{\pp})z^{-\ell}e^{-\sum 
y_j\pp_j}\tau_q
\sum p_n(\tl{\pp})z^n e^{\sum y_j\pp_j}\tau_q+
\end{equation}
$$+z[1+x(q-1)D_q]\sum p_{\ell}(-\tl{\pp})z^{-\ell}e^{-\sum 
y_j\pp_j}\tau_q\sum p_n
(\tl{\pp})z^{-n}e^{\sum y_j\pp_j}\tau_q\}dz=0$$
Picking the coefficient of $y_r$ again we get $r-\ell-n=-1$ for the 
first terms and
$r+1-\ell-n=-1$ for the second leading to
\bq\label{99}
0=-2\sum_{\ell+n=r+1}p_{\ell}(-\tl{\pp})\tau_qp_n(\tl{\pp})\tau_q-2\sum_{\ell+n=r+2}
p_{\ell}(-\tl{\pp})[1+x(q-1)D_q]\tau_q p_n(\tl{\pp})\tau_q-
\end{equation}
$$-\sum_{\ell+n=r+1}p_{\ell}(-\tl{\pp})\pp_r\tau_qp_n(\tl{\pp})\tau_q+
\sum_{\ell+n=r+2}p_{\ell}(-\tl{\pp})[1+x(q-1)D_q]\tau_q 
p_n(\tl{\pp})\pp_r\tau_q$$
This seems to provide the simplest Hirota equation involving $D_q$ directly.
\begin{theorem}
Equations \eqref{95} and \eqref{99} provide typical Hirota equations 
in the context of
Theorem 4.3.
\end{theorem}
\indent
Even for these simplest of equations the calculations leading to qKP type equations
are tedious, with little simplification arising (see Section 10 for a more
productive approach involving qKdV).

\section{EQUATIONS ON THE Q-PLANE}
\renewcommand{\theequation}{5.\arabic{equation}}
\setcounter{equation}{0}

In \cite{c4} we sketched an heuristic derivation of a q-plane 
($xt=qtx$) Burger's type equation (as in Example 3.1)
from a zero curvature condition $dA+A^2=0$ where 
$df=D_t^{-1}\pp_q^xfdx+\pp_q^tfdt$ (with e.g.
$\pp_q^xf=D_{q^{-2}}fdx$ etc.).  This was not very revealing so we 
want to examine in more
detail the zero curvature idea on a quantum plane following e.g. 
\cite{cxx,c10,c11,d8,d9,f5,m2}
(cf. also Section 6).  We note first a variation on this from \cite{c4}, 
namely
\begin{example}
Try e.g. $xdx=q^2dxx$ ($\gG_{+}$ style as in \cite{cxx,k1}) with $[dx,t]=adt,
[dx,dt]=0=(dt)^2=(dx)^2$, and $[dt,x]=adt$ for consistency.  Then $dx^n=
[[n]]_{q^{-2}}x^{n-1}dx=[[n]]_{q^2}dxx^{n-1}$; we write this here as
$\pp_qx^ndx$ (recall
$[[n]]_{q^2}=(q^{2n}-1)/(q^2-1)$).  Further take $dtt=tdt$ so 
\bq\label{186}
df=\pp_qfdx+(a+1)f_tdt
\end{equation}
Try $A=wdt+udx$ as before so
\bq\label{187}
dA=(\pp_qw+(a+1)u_t)dxdt
\end{equation}
For $A^2$ we note that $xt=tx$ is compatible with the conditions above and
e.g. $(\clubsuit)\,\,wdtwdt=w(x,t)w(x+a)(dt)^2=0$ while for $udxudx$ one
notes
$dxf(x)=D^{-2}fdx$ while $(\spadesuit)\,\,dxt^m=t^mdx+mat^{m-1}dt$ implies
that
$dxu=D_x^{-2}udx+aD_x^{-2}u_tdt$.  This leads to
\bq\label{189}
A^2=[w(x,t)w(x+a,t)+uaD_x^{-2}u_t+uD_x^{-2}w]dtdx
\end{equation}
Hence $dA+A^2=F=0$ requires
\bq\label{190}
w(x,t)w(x+a,t)+auD_x^{-2}u_t+uD_x^{-2}w+\pp_qw+(a-1)u_t=0
\end{equation}
Note if $q\to 1$ one obtains $(\bullet)\,\,w^2+uw+w_x+u_t=0$ which for
$w=u_x$ implies $(\bl)\,\,u_x^2+uu_x+u_{xx}+u_t=0$
which is a kind of perturbation of the Burger's equation 
by a $(u_x)^2$ term (constants can be adjusted by changes of variables etc.)
$\hfill\bs$
\end{example}
\indent
One can examine the zero curvature idea on a q-plane for example at some
length (cf. \cite{cxx,c10,c11,d8,d9,f5,m2} and this gets into a long
discussion of differential calculi (DC) and differential geometry (DG) over
algebras (cf. \cite{cxx} for a survey and see also \cite{w1,wq} for some
fundamental information).  In this paper we only give a brief treatment of
some essential ideas and concentrate on some examples.  First from
\cite{d9} one considers the q-plane with generators $x^i\sim (x,y)$ 
and $\xi^i\sim dx^i\sim
(\xi,\eta)$ for the algebra of forms $\gO^*=\gO^0\ot\gO^1\ot\gO^2$ 
(cf. also Remark 3.2).
The relations are
\bq\label{500}
xy=qyx;\,\,\xi^2=0;\,\,\eta^2=0;\,\,\eta\xi+\xi\eta=0;\,\,x\xi=q^2\xi x;
\end{equation}
$$x\eta=q\eta x+(q^2-1)\xi y;\,\,y\xi=q\xi y;\,\,y\eta=q^2\eta y$$
These conditions can be written in the form
\bq\label{501}
\hat{R}=\left(\begin{array}{cccc}
q & 0 & 0 & 0\\
0 & q-q^{-1} & 1 & 0\\
0 & 1 & 0 & 0\\
0 & 0 & 0 & q
\end{array}\right);\,\,\,\,
\begin{array}{c}
x^ix^j-q^{-1}\hat{R}^{ij}_{k\ell}x^kx^[{\ell}=0\\
x^i\xi^j-q\hat{R}^{ij}_{k\ell}\xi^kx^{\ell}=0\\
\xi^i\xi^j+q\hat{R}^{ij}_{k\ell}x^k\xi^{\ell}=0
\end{array}
\end{equation}
Writing the generators $(a,b,c,d)$ of $SL_q(2,{\bf C})$ in matrix 
form the invariance of $\gO^*$
under $SL_q(2{\bf C})$ action follows via
\bq\label{502}
\hat{R}^{ij}_{k\ell}a_m^ka_n^{\ell}=a^i_ka^j_{\ell}\hat{R}^{k\ell}_{mn};\,\,\,\,
a_j^i=\left(\begin{array}{cc}
a & b\\
c & d
\end{array}\right)
\end{equation}
Now introduce the trivial differential calculus (DC) on $SL_q(2,{\bf 
C})$ via ${\bf (A61)}\,\,
da_j^i=0$ so that coaction on $x^i$ and $\xi^i$ is described by ${\bf 
(A62)}\,\,x^{i'}=
a_j^i\ot x^j$ and $\xi^{i'}=a^i_j\ot \xi^j$.  From \eqref{502} it 
follows that $x^{i'}$ and $\xi^{i'}$
satisfy the same relations as $x^i$ and $\xi^i$.  Then introduce the 
1-form ${\bf (A63)}\,\,
\gt=x\eta-qy\xi$.  Evidently $\gt^2=0$ and $\gt$ is invariant under 
coaction of $SL_q(2,{\bf C})$
where ${\bf (A64)}\,\,\gt^{'}=1\ot \gt$.  In fact, up to scalar 
multiples, it is the only invariant
element of $\gO^1$ and one has ${\bf (A65)}\,\,x^i\gt=q\gt x^i$ with 
$\xi^i\gt=-q^{-3}\gt\xi^i$.
To fix the definition of a covariant derivative one must first 
introduce a suitable operator
${\bf 
(A66)}\,\,\gs:\,\gO^1\ot_A\gO^1\stackrel{\gs}{\to}\gO^1\ot_A\gO^1$ 
satisfying
$\pi(\gs+1)=0$
where $\pi$ is the projection $\gO(V)\ot_{C(V)}\gO^1(V)\to\gO^2(V)$ (do
not confuse this $\gs$ with the $\gs$ of Remark 3.1).  Here $\gs$ is 
the inverse of $q\hat{R}$
written as
\bq\label{503}
\gs(\xi\ot 
\xi)=q^{-2}\xi\ot\xi;\,\,\gs(\xi\ot\eta)=q^{-1}\eta\ot\xi;\,\,\gs(\eta\ot\eta)=
q^{-2}\eta\ot\eta;
\end{equation}
$$\gs(\eta\ot\xi)=q^{-1}\xi\ot\eta-(1-q^2)\eta\ot\xi$$
The extension to $\gO^1\ot_{\gO^0}\gO^1$ is given by the right 
$\gO^0$ linearity, and in this
situation $\gs$ is also left $\gO^0$ linear (with {\bf (A66)} 
evidently satisfied).  As a result of
linearity of finds
\bq\label{504}
\gs(\xi\ot\gt)=q^{-3}\gt\ot\xi;\,\,\gs(\gt\ot\xi)=q\xi\ot\gt-(1-q^{-2})\gt\ot\xi;
\end{equation}
$$\gs(\eta\ot\gt)=q^{-3}\gt\ot\eta;\,\,\gs(\gt\ot\eta)=q\eta\ot\gt-(1-q^{-2})\gt\ot\eta$$
along with ${\bf (A67)}\,\,\gs(\gt\ot\gt)=q^{-2}\gt\ot\gt$.  Note 
that, even though $\gt^2\ne 1$,
$\gs$ satisfies the Hecke relation ${\bf 
(A68)}\,\,(\gs+1)(\gs-q^{-2})=0$.  Now there is a unique
one parameter family of covariant derivatives compatible with the 
algebraic structure \eqref{500}.
It is given via ${\bf (A69)}\,\,D\xi=\mu^4x^i\gt\ot\gt$ where $\mu$ 
has the dimensions of inverse
length.  From the invariance of $\gt$ it follows that D is invariant 
under the coaction of
$SL_q(2,{\bf C})$ and from {\bf (A63)} one sees that the torsion 
$(d-\pi D)\gt=0$.
\\[3mm]\indent
{\bf REMARK 5.1.}
Now a few remarks about covariant derivatives may be helpful.  Thus let 
V be a differential manifold
and let $(\gO^*(V),d)$ be the ordinary DC on V.  Let H be a vector 
bundle over V associated to some
principal bundle P.  Let $C(V)$ be the algebra of smooth functions on 
V and ${\mf H}$ the left
$C(V)$ module of smooth sections of H.  A connection on P is 
equivalent to a covariant derivative on
H.  This in turn can be characterized as a linear map
\bq\label{505}
{\mf H}\stackrel{D}{\to}\gO^1(V)\ot_{C(V)}{\mf H};\,\,D(f\psi)=df\ot\psi+fD\psi
\end{equation}
for arbitrary $f\in C(V)$ and $\psi\in{\mf H}$.  There is an 
immediate extension of D to a map
${\bf (A70)}\,\,\gO^*(V)\ot_{C(V)}{\mf H}\to\gO^*(V)\ot_{C(V)}{\mf 
H}$ be requiring that it be an
antiderivation of degree 1.  This definition has a direct extension 
to NC geometry.  Thus let A be
an algebra and $(\gO^*(A),d)$ a DC over A.  One defines a covariant 
derivative on a left A-module
${\mf H}$ as a map ${\bf (A71)}\,\,{\mf 
H}\stackrel{D}{\to}\gO^1(A)\ot{\mf H}$ which satisfies
\eqref{505} with $f\in A$.  Again there is an extension of D to a map 
$\gO^*(A)\ot_A{\mf H}
\to\gO^*(A)\ot_A{\mf H}$ as an antiderivation of degree 1.  Next a 
linear connection on V can be
defined as a connection on the cotangent bundle to V, characterized via
\bq\label{506}
\gO^1(V)\stackrel{D}{\to}\gO^1(V)\ot_{C(V)}\gO^1(V);\,\, D(f\xi)=df+fD\xi
\end{equation}
for $f\in C(V)$ and $\xi\in\gO^1(V)$.
By extension a linear connection over a general NC algebra A with a DC
$(\gO^*(A),d)$ can be defined as a linear map ${\bf
(A73)}\,\,\gO^1\stackrel{D}{\to}\gO^1\ot_A\gO^1$ satisfying 
\eqref{506} for $f\in A$
and $\xi\in\gO^1=\gO^1(A)$.  The module $\gO^1$ has a natural 
structure as a right
A-module but in the NC case it is impossible in general to 
consistently impose the
condition
$D(\xi f)=D(f\xi)$ and a substitute must be found.  Here one 
postulates the existence of a map
$\gO^1\ot_A\gO^1\stackrel{\gs}{\to}\gO^1\ot_A\gO^1$ satisfying 
$\pi(\gs+1)=0$ as in {\bf (A66)}.
One defines then $D(\xi f)$ by ${\bf (A74)}\,\,D(\xi f)=\gs(\xi\ot 
df)+(D\xi)f$ and requires also
that $\gs$ be right linear.
The extension of D to $\gO^1\ot_A\gO^1$ is given by ${\bf 
(A75)}\,\,D(\xi\ot\eta)=D\xi\ot\eta
+\gs_{12}(\xi\ot D\eta)$ where $\gs_{12}=\gs\ot 1$.  Finally a metric 
$g$ on V can be defined as a
$C(V)$-linear symmetric map of $\gO^1(V)\ot_{C(V)}\gO^1$ into $C(V)$. 
This makes also if one
replaces $C(V)$ by an algebra A and $\gO^1(V)$ by $\gO^1(A)$.  By 
analogy with the commutative case
one says that the covariant derivative {\bf (A74)} is metric if the 
following digram is commutative
\[\begin{CD}
\\\gO^1\ot_A\gO^1             @>{D}>>         \gO^1\ot_A\gO^1\ot_A\gO^1\\
@V{g}VV                                      @VV{1\ot g}V\\
A                          @>>{d}>             \gO^1
\end{CD}\]
Since symmetry must in general be imposed relative to $\gs$ one 
requires of $g$ then that ${\bf
(A76)}\,\,g\gs=g$.
Going back to the q-plane covariant derivative {\bf (A69)} one can 
show that it is without torsion
but there is no solution to $g\ci\gs=g$ so the connection is not 
metric (see here also \cite{m2}).$\hfill\bs$
\\[3mm]\indent
{\bf REMARK 5.2.}
Now go to \cite{d8} where a more general picture emerges.  The 
relevant idea here for us is to work
on a generalized q-plane with algebra A generated by 
$(x,y,x^{-1},y^{-1})$ subject to $xy=qyx$
and the usual relations between elements and their inverses.  First 
on background we recall that
the largest DC over A is the differential envelope or universal DC 
(UDC) $(\gO^*(A),d_u)$
(cf. \cite{cxx}).  Every
other DC can be considered as a quotient of it.  Indeed let 
$(\gO^*(A),d)$ be another DC; then
there exists a unique $d_u$-homomorphism $\phi$ with
\[\begin{CD}
\\A            @>{d_u}>>     \gO^1_u(A)          @>{d_u}>> 
\gO^2_u(A)\cdots     \\
@V{id}VV                      @VV{\phi_1}V 
@VV{\phi_2}V       \\
A              @>>{d}>         \gO^1(A)         @>>{d}>        \gO^2(A)\cdots
\end{CD}\]
given via $(\clubsuit)\,\,\phi(d_ua)=da$.  The restriction 
$\phi_p=\phi_{\gO^p_u}$ is
defined by $(\spadesuit)\,\,\phi_p(a_0d_ua_1\cdots d_ua_p)=a_0da_1\cdots da_p$.
Consider a given algebra A and suppose that one knows how to 
construct an A-module $\gO^1(A)$ and
an application $(\bullet)\,\,A\stackrel{d}{\to}\gO^1(A)$.  Then using 
the diagram there is a method
of constructing $\gO^p(A)$ as well as the extension of the 
differential.  Since we know $\gO_u^1(A)$
and $\gO^1(A)$ we can suppose $\phi_1$ is given.  The simplest choice 
of $\gO^2(A)$ would then be
$(\bl)\,\,\gO^2(A)=\gO^2_u(A)/d_uKer(\phi_1)$; this is the largest DC 
consistent with the
constraints on $\gO^1(A)$.  The procedure defines $\phi_2$ and d is 
defined via $d(fdg)=dfdg$
(cf. \cite{cxx}).  Now to initiate the above construction one starts 
with the 1-forms using a set
of derivations.  Suppose they are interior derivations and this 
excludes the case where A is
commutative.  For each integer n let $\gl_a$ be a set of n linearly 
independent elements of A and
introduce the derivations $e_a=ad(\gl_a)$.  In general the $e_a$ do 
not form a Lie algebra but they
do satisfy commutation relations as a consequence of the commutation 
relations of A.  Define
the map $(\bullet)$ via ${\bf (A77)}\,\,df(e_a)=e_af$.  Then suppose 
there is a set of n elements
$\gt^a\in\gO^1(A)$ such that ${\bf (A78)}\,\,\gt^a(e_b)=\gd^a_b$.  In 
the examples such $\gt^a$
will exist by explicit construction and they are referred to as a 
frame which commutes with the
elements $f\in A$ via ${\bf (A79)}\,\,f\gt^a=\gt^af$.
\\[3mm]\indent
The A-bimodule $\gO^1(A)$ is generated by all elements of the form 
$fdg$ or of the form $(df)g$
and because of the Leibnitz rule these conditions are equivalent. 
Define now $\gt=-\gl_a\gt^a$
and one sees that ${\bf (A80)}\,\,df=-[\gt,f]$ so that as a bimodule 
$\gO^1(A)$ is generated by one
element.  It follows also that the 2-form $d\gt+\gt^2$ can be written 
in the form ${\bf (A81)}\,\,
d\gt+\gt^2=-(1/2)K_{ab}\gt^a\gt^b$ with coefficients $K_{ab}$ which 
lie in $Z(A)$.
By definition ${\bf (A82)}\,\,fdg(e_a)=fe_ag$ and 
$(dg)f(e_a)-(e_ag)f$ and using the frame one can
write this as ${\bf (A83)}\,\,fdg=(fe_ag)\gt^a$ and 
$(dg)f=(e_ag)f\gt^a$.  The commutation
relations of the algebra naturally constrain the relations between 
$fdg$ and $(dg)f$.  
Now as a
right module $\gO^1(A)$ is free of rank n, and because of the 
commutation relations of the algebra,
or equivalently, because of the kernel of $\phi_1$ in the quotient 
$(\bl)$ the $\gt^a$ satisfy
commutation relations which are supposed to be of the form ${\bf 
(A84)}\,\,\gt^a\gt^b+C^{ab}_{cd}
\gt^c\gt^d=0$.  If $C^{ab}_{cd}=\gd^a_c\gd^b_d$ then $\gO^2(A)=0$ and 
it follows from {\bf (A84)}
that for an arbitrary element $f\in A$ ${\bf 
(A85)}\,\,[f,C^{ab}_{cd}]\gt^c\gt^d=0$.  One supposes
now that ${\bf (A86)}\,\,C^{ab}_{ef}C^{ef}_{cd}=\gd^a_c\gd^b_d$ and 
that the relations {\bf (A84)}
are complete in the sense that if $A_{ab}\gt^a\gt^b=0$ then ${\bf 
(A87)}\,\,A_{ab}-C^{cd}_{ab}A_{cd}
=0$ (this will hold for DC on the generalized quantum plane 
considered here).  One concludes
that $C^{ab}_{cd}\in Z(A)$.  Note that in ordinary geometry one can 
choose ${\bf (A88)}\,\,
C^{ab}_{cd}=\gd^b_c\gd^a_d$ and the relation {\bf (A84)} expresses 
the anticommutativity of
1-forms.  Now let $\gL^*_C$ be the twisted exterior algebra 
determined by {\bf (A84)} and it follows
that ${\bf (A89)}\,\,\gO^*=A\ot_{{\bf C}}\gL_C^*$.  Because the 
2-forms are generated by products
of the $\gt^a$ one has ${\bf 
(A90)}\,\,d\gt^a=-(1/2)C^a_{bc}\gt^b\gt^c$ and one can assume without
loss of generality that ${\bf 
(A91)}\,\,C^a_{bc}+C^a_{de}C^{de}_{bc}=0$.   Note however that
$C^a_{bc}\ne Z(A)$ in general.  In fact from the identity 
$d(f\gt^a)=d(\gt^af)$ one has
${\bf (A92)}\,\,\{(1/2)[C^a_{bc},f]+e_{(b}f\gd^a_{c)}\}\gt^b\gt^c=0$. 
Using the definition of
derivations one can write this in the form
\bq\label{508}
\left(\frac{1}{2}C^a_{bc}+\gl_{(b}\gd^a_{c)}-\frac{1}{2}D^a_{bc}\right)\gt^b\gt^c=0
\end{equation}
with $D^a_{bc}\in Z(A)$.  One can also suppose that the $D^a_{bc}$ 
satisfy {\bf (A91)},
namely $D^a_{bc}+D^a_{de}C^{de}_{bc}=0$.  Now using this and the 
relations {\bf (A84)} and
{\bf (A91)} as well as the completeness assumption {\bf (A87)} one 
can conclude from \eqref{508}
that ${\bf 
(AA92)}\,\,C^a_{bc}-D^a_{bc}+\gl_{(b}\gd^a_{c)}-\gl_{(d}\gd^a_{e)}C^{cd}_{bc}=0$.
The equation \eqref{508} can also be written in the form ${\bf 
(A93)}\,\,d\gt^a=-[\gt,\gt^a]
-(1/2)D^a_{bc}\gt^b\gt^c$ with a graded commutator.  If $D^a_{bc}=0$ 
the form {\bf (A80)}
for the extrerior derivative is valid for all elements of $\gO^*(A)$ 
and the element $\gt$
plays the role of the phase of a generalized Dirac operator as in 
\cite{c11}.  Now from
{\bf (A93)} one finds that ${\bf 
(A94)}\,\,d\gt=-2\gt^2+(1/2)\gl_aD^a_{bc}\gt^b\gt^c$.
Comparing this with {\bf (A81)} one concludes that ${\bf 
(A95)}\,\,\gt^2=(1/2)(\gl_aD^a_{bc}
+K_{bc})\gt^b\gt^c$.  If we suppose that $K_{bc}$ satisfies {\bf 
(A91)} so ${\bf (A96)}\,\,
K_{ab}+K_{cd}C^{cd}_{ab}=0$ then ${\bf 
(A97)}\,\,(\gl_b\gl_c-C^{de}_{bc}\gl_d\gl_e-\gl_a
D^a_{bc}-K_{bc})\gt^b\gt^c=0$.  Using again {\bf (A87)} there results 
${\bf (A98)}\,\,
\gl_b\gl_c-C^{de}_{bc}\gl_d\gl_e=\gl_aD^a_{bc}+K_{bc}$.  Now 
introduce the twisted bracket
${\bf (A99)}\,\,[\gl_a,\gl_b]_C=\gl_a\gl_b-C^{cd}_{ab}\gl_c\gl_d$ and 
rewrite {\bf (A98)} in the
form ${\bf (A100)}\,\,[\gl_b,\gl_c]_C=\gl_aD^a_{bc}+K_{ab}$.  If one 
then writes out the equation
$d^2f=0$ using {\bf (A87)}-{\bf (A91)} there results ${\bf 
(A101)}\,\,[e_b,e_c]_Cf=e_afC^a_{bc}$.
This is the dual relation to the Maurer-Cartan equation {\bf (A90)}.
The construction {\bf (A100)} follows from the relation {\bf (A84)} 
assumed for the frame as well
as from the conditions imposed on the $C^{ab}_{cd}$.  In the matrix 
case the general formalism
simplifies considerably via ${\bf (A102)}\,\,\gt^a=\gl_b\gl^ad\gl^b$. 
The elements of the frame
then anticommute and one can choose $C^{ag}_{cd}=\gd^b_c\gd^a_d$. 
Further in {\bf (A93)} the first
term on the right vanishes and $D^a_{bc}=C^a_{bc}$; moreover on the 
right side of {\bf (A81)}
and {\bf (A100)} one has $K_{ab}=0$.$\hfill\bs$
\begin{example}
Finally we look at a specific example for the generalized quantum 
plane generated by $(x,y,
x^{-1},y^{-1})$ with $xy=qyx$.  Define for $q\ne 1$ ${\bf 
(A103)}\,\,\gl_1=qy/(q-1)$ and
$\gl_2=qx/(q-1)$.  The normalization has been chosen so that the 
structure elements $C^a_{12}$
contain no factors q.  Note that the $\gl_a$ are singular in the 
limit $q\to 1$ in the spirit of
the QM limit ${\bf (A104)}\,\,(1/\hbar)ad(p)\to (1/i)(\pp/\pp q)$. 
One gets ${\bf (A105)}\,\,
e_1x=-xy,\,\,e_qy=0,\,\,e_2x=0,$ and $e_2y=xy$.  These rather unusual 
derivations are extended to
arbitrary polynomials in the generators via Leibnitz and from {\bf 
(A105)} one concludes that the
commutation relations following from {\bf (A83)} are
\bq\label{509}
xdx=qdxx;\,\,ydx=q^{-1}dxy;\,\,xdy=qdyx;\,\,ydy=q^{-1}dyy
\end{equation}
Using the relations {\bf (A105)} one finds ${\bf 
(A107)}\,\,dx=-xy\gt^1$ and $dy=xy\gt^2$ so
solving for the $\gt^a$ gives ${\bf 
(A108)}\,\,\gt^1=-q^{-1}x^{-1}y^{-1}dx$ and $\gt^2=
q^{-1}x^{-1}y^{-1}dy$ satisfying ${\bf 
(A109)}\,\,(\gt^1)^2=0,\,\,(\gt^2)^2=0$, and $\gt^1\gt^2
+q\gt^2\gt^1=0$ which has the form {\bf (A84)}.  If one reorders the 
indices $(11,12,21,22)=
(1,2,3,4)$ then $C^{ab}_{cd}$ can be written as
\bq\label{510}
C=\left(\begin{array}{cccc}
1 & 0 & 0 & 0\\
0 & 0 & q & 0\\
0 & q^{-1} & 0 & 0\\
0 & 0 & 0 & 1
\end{array}\right)
\end{equation}
Thus $C^{12}_{21}=q$ and $C^{21}_{12}=q^{-1}$.  The structure 
elements $C^a_{bc}$ are given via
${\bf (A110)}\,\,C^1_{12}=-x$ and $C^2_{12}=-y$ plus {\bf (A91)}. 
Equation {\bf (AA92)} is
satisfied and for $\gt$ one finds the expression ${\bf
(A111)}\,\,\gt=(q-1)^{-1}(qx^{-1}dx-y^{-1}dy)$ which is a closed form.
\end{example}
\begin{example}
As a second example one look at the Wess-Zumino calculus determined 
for $q^4\ne 1$ via
${\bf (B1)}\,\,\gl_1=(q^4-1)^{-1}x^{-2}y^2$ and 
$\gl_2=(q^4-1)^{-1}x^{-2}y^2$.  The normalization
is chosen so that the $C^a_{12}$ contain no q factors.  Then for $q^2\ne -1$
\bq\label{513}
e_1x=\frac{-x^{-1}y^2}{q^2(q^2+1)};\,\,e_1y=\frac{-x^{-2}y^3}{(q^2+1};\,\,e_2x=0;\,\,
e_2y=\frac{-x^{-2}y}{q^2+1}
\end{equation}
 From this the commutation relations following from {\bf (A83)} are the familiar
\bq\label{514}
xdx=q^2dxx;\,\,xdy=qdyx+(q^2-1)dxy;\,\,ydx=qdxy;\,\,ydy=q^2dyy
\end{equation}
(covariant DC of Wess-Zumino).  If $q^2\ne -1$ one gets ${\bf 
(B2)}\,\,(dx)^2=0,\,\,(dy)^2=0,$ and
$dydx+qdxdy=0$.  Using the relation {\bf (A78)} one finds ${\bf
(B3)}\,\,dx=-x^{-1}y^2\gt^1/[q^2(q^2+1)]$ and 
$dy=-x^{-2}y(y^2\gt^1+\gt^2)/(q^2+1)$ and solving
for the $\gt^a$ gives ${\bf (B33)}\,\,\gt^1=-q^4(q^2+1)xy^{-2}dx$ and 
$\gt^2=-q^2(q^2+1)x
(xy^{-1}dy-dx)$.  There are commutations ${\bf 
(B5)}\,\,(\gt^1)^2=0,\,\,(\gt^2)^2=0,$
and $q^4\gt^1\gt^2+\gt^2\gt^1=0$.  This is of the form {\bf (A84)} if 
$C^{ab}_{cd}$ is given by
\bq\label{515}
C=\left(\begin{array}{cccc}
1 & 0 & 0 & 0\\
0 & 0 & q^{-4} & 0\\
0 & q^4 & 0 & 0\\
0 & 0 & 0 & 1
\end{array}\right)
\end{equation}
Thus $C^{12}_{21}=q^{-4}$ and $C^{21}_{12}=q^4$.  The structure 
elements $C^a_{bc}$ are
${\bf (B6)}\,\,C^1_{12}=-x^{-2}$ and $C^2_{12}=-x^{-2}y^2$.  Also 
{\bf (A91)} and {\bf (AA92)}
are satisfied.  For $\gt$ one has ${\bf 
(B7)}\,\,\gt=q^2y^{-1}dy/(q^2-1)$ which is a closed form.
$\hfill\bs$
\end{example}

\section{THE Q-DEFORMED LINE}
\renewcommand{\theequation}{6.\arabic{equation}}
\setcounter{equation}{0}

We note now that the search for Burger's equation or qKP type 
equations on a quantum plane is
somewhat distorted since only $D_q^x$ is involved and other variables 
$t_i$ play a traditional role.
Thus we go to \cite{c12,f6,f7} for some background information since 
ad hoc calculations are not
clearly meaningful.  We extract first from the lovely paper \cite{c12}.
Thus a DC over A is another associative algebra $\gO^*(A)$ with a 
differential d, which plays the
role of the deRham DC and must tend to this in the commutative limit. 
The DC is what gives
structure to the set of ``points" and determines the dimension.  It 
would determine the nearest
neighbors in a lattice.  Over a given A there are many possible DC 
and the choice depends on what
limit manifold is in mind.  One can start with derivations as a set 
of linear maps of the algebra
into itself satisfying a Leibnitz rule and use them as a basis for 
constructing differential forms
or one can start with a set of differential forms and construct a set 
of possibly twisted
derivations which are dual to the forms (twisted in the sense of 
obeying a modified Leibnitz rule).
One introduces derivations $e_af=[\gl_a,f]$ as before and suppose the 
algebra generated by $x^i$
with $1\leq i\leq n$.  Defining $df(e_q)=e_af$ as before one finds in 
general that
$dx^i(e_a)\ne\gd^i_a$.  Then one constructs a new basis $\gt^a$ dual 
to the basis or derivations
$\pp_i$ obeying a modified Leibnitz rule such that 
$dx^i(\pp_j)=\gd^i_j$.  In general both
approaches are equivalent.  By construction the $\gt^a$ commute with 
elements of A and define the
structure of the 1-forms as a bimodule over A.  
It is also appropriate
to add an extra generator $\gL$ called the dilatator and its inverse 
$\gL^{-1}$ chosen so that
${\bf (B8)}\,\,x^i\gL=q\gL x^i$ with $\gL$ unitary; since $r$ and 
$\gL$ do not commute the center
of the corresponding extension is trivial.  For $n=1$ now there are 
two generators $x$ and $\gL$
with $x\gL=q\gL x$ and one chooses $x$ Hermitian and $q\in {\bf 
R}^{+}$ with $q>1$.  This is a
modified version of the Weyl algebra with q real instead of having 
unit modulus.  One can represent
the algebra on a Hilbert space ${\mf R}_q$ with basis $|k>$ via ${\bf 
(B9)}\,\,x|k>=q^k|k>$ and
$\gL|k>=|k+1>$.  The spectrum of $\gL$ will be continuous (in 
contrast to the analogous Schwinger
basis).  Introduce now $y$ by the action ${\bf (B10)}\,\,y|k>=k|k>$ 
on basis elements.  Then
${\bf (B11)}\,\,\gL^{-1}y\gL=y+1$ and one can write $x=q^y$ as an 
equality within ${\bf R}_q^1$.
One occasionally renormalizes y via ${\bf (B12)}\,\,z=q^{-1}(q-1)>0$ 
with ${\bf (B13)}\,\,zy\to
y$.  With the new value of y the spacing between the spectral lines 
vanishes with z; the old
units will be referred to as Planck units and the new ones as laboratory units.
With the new value of $t$ the spacing between the spectral lines 
vanishes with $z$.
\\[3mm]\indent
One possible DC over ${\bf R}_q^1$ is constructed by setting $d\gL=0$ 
and ${\bf (B14)}\,\,
xdx=dxx$ with $dx\gL=q\gL dx$ with frame $\gt^1=x^{-1}dx$.  This 
calculus has an involution
$(dx)^*=dx^*$ but it is not based on derivations and it has no 
covariance properties with respect to $SO_q(1)$.
Hence consider another DC $\gO^*({\bf R}_q^1)$ based on 
the relations ${\bf (B15)}\,\,
xdx=qdx x$ and $dx\gL=q\gL dx$.  If one chooses $\gl_1=-z^{-1}\gL$ 
then ${\bf (B16)}\,\,
e_1x=q\gL x$ and $e_1\gL=0$ with the calculus defined by the 
condition $df(e_1)=e_1f$.  By setting
${\bf (B17)}\,\,\gl_2=z^{-1}x$ and introducing a second derivation 
${\bf (B18)}\,\,e_2\gL=
q\gL x$ and $e_2x=0$ we could extend the calculus via $df(e_2)=e_2f$ 
and find that $xd\gL=
qd\gL x$.  However this will not be done since $\gL$ is
in a certain sense an element of the
phase space associated to x and one is more interested here in position 
space geometry.
Now an adjoint derivation $e_1^{\dg}$ of $e_1$ is defined via ${\bf 
(B19)}\,\,e_1^{\dg}=(e_1f^*)^*$.
The $e_1^{\dg}$ on the left is not an adjoint of an operator $e_1$; 
it is defined uniquely in terms
of the involution whereas $e_1$ acts on the algebra as a vector 
space.  Since $\gL$ is unitary one
has $(\gl_1)^*\ne-\gl_1$ and thus $e_1$ is not a real derivation. 
Hence introduce a second DC
$\bar{\gO}^*({\bf R}_q^1)$ defined via ${\bf 
(B20)}\,\,x\bar{d}x=q^{-1}\bar{d}xx$ and
$\bar{d}x\gL=q\gL\bar{d}x$ and based on the derivation $\bar{e}_1$ formed using
$\bar{\gl}_1=-\gl_1^*$.  This calculus is defined by the condition 
$\bar{d}f(\bar{e}_1)=\bar{e}_1f$
and $\bar{e}_1$ is also not real.  However ${\bf 
(B21)}\,\,e_1^{\dg}=\bar{e}_1$ and therefore
$(df)^*=\bar{d}f^*$.  By induction ${\bf 
(B22)}\,\,e_1x^n=z^{-1}(q^n-1)\gL x^n$ and
$\bar{e}_1x^n=z^{-1}(1-q^{-n})\gL^{-1}x^n$.  One can also represent 
$\gO^*$ and $\bar{\gO}^*$
on ${\mf R}_q$ via ${\bf (B23)}\,\,dx|k>=\ga q^{k+1}|k+1>$ and 
$\bar{d}|k>=\bar{\ga}q^k
|k-1>$ with two arbitrary parameters $\ga$ and $\bar{\ga}$.  One has 
$dx^*=\bar{d}x$ if and only if
$\ga^*=\bar{\ga}$; also generally ${\bf (B24)}\,\,d=-z^{-1}\ga 
ad(\gL)$ and $\bar{d}=z^{-1}\bar{\ga}
ad(\gL^{-1})$ and the relations {\bf (B15)} and {\bf (B20)} hold. 
The frame elements $\gt^1,\,
\bar{\gt}^1$ dual to the derivations $e_1,\,\bar{e}_1$ are given by
\bq\label{516}
\gt^1=\gt_1^1dx;\,\,\gt^1_1=\gL^{-1}x^{-1};\,\,\bar{\gt}^1=\bar{\gt}^1_1\bar{d}x;\,\,
\bar{\gt}_1^1=q^{-1}\gL x^{-1}
\end{equation}
On ${\mf R}_q$ they become the operators ${\bf (B25)}\,\,\gt^1=\ga$ 
and $\bar{\gt}^1=\bar{\ga}$
proportional to the unit element.  Since $\gO^1$ and $\bar{\gO}^1$ 
are free ${\bf R}_q^1$ modules
of rank one with special bases $\gt^1$ and $\bar{\gt}^1$ one has 
${\bf (B26)}\,\,\gO^*=\gL^*
\ot{\bf R}_q^1$ and $\bar{\gO}^*=\gL^*\ot{\bf R}_q^1$ where $\gL^*$ 
is the exterior algebra over
${\bf C}$ so the extension is trivial.  One can also construct a 
real DC $\gO^*_R$ which
can be developed generally over an algebra A; since this becomes rather
complicated and not really needed we omit it here.  We also omit questions about
metrics, gauge fields, etc. and refer to Section 9 for an example.

\section{ZERO CURVATURE EXAMPLES}
\renewcommand{\theequation}{7.\arabic{equation}}
\setcounter{equation}{0}

We explore now some of the frameworks disclosed in earlier sections.
\begin{example}
Consider the generalized q-plane of Example 5.2 where 
$xdx=qdxx,\,\,ydx=q^{-1}dxy,\,\,
xdy=qdyx,$ and $ydy=q^{-1}dyy$.  Also from $qdyx=dxy$ we have 
$qdydx=dxdy$ and a little
calculation yields ${\bf (C1)}\,\,dx^n=[(1-q^{-n})/(1-q^{-1})]x^{n-1}dx$
with $dy^m=[(1-q^m)/(1-q)]y^{m-1}dy$.  Working from $f=\sum 
a_{nm}x^ny^m$ one obtains then
(note $dxy^m=q^my^mdx$)
\bq\label{534}
df=D_y\pp^x_{q^{-1}}fdx+\pp_q^yfdy
\end{equation}
Set then $A=wdy+udx$ with $dA=D_y\pp^x_{q^{-1}}wdxdy+\pp_q^yudydx$ 
and, noting that
$dyx^n=q^{-n}x^ndy,\,\, dyy^m=q^my^mdy,$ and $dxy^m=q^my^mdx$
with $dxx^n=q^{-n}x^ndx$ one gets
$dyw=D_x^{-1}D_ywdy$ and $dxu=D_x^{-1}D_yudx$ leading to
$A^2=wD_x^{-1}D_yudydx+uD_x^{-1}D_ywdxdy$; and
\bq\label{535}
dA+A^2=0=qD_y\pp_{q^{-1}}^xw+\pp_q^yu+
wD_x^{-1}D_yu+quD_x^{-1}D_yw
\end{equation}
Setting then e.g. $qw=D_y^{-1}\pp_{q^{-1}}^xu$ one gets
\bq\label{536}
\pp_q^yu+(\pp^x_{q^{-1}})^2u+q^{-1}(D_x^{-1}D_yu)(D_y^{-1}\pp_{q^{-1}}^xu)+
uD_x^{-1}\pp^x_{q^{-1}}u=0
\end{equation}
For $q\to 1$ we have ${\bf (C2)}\,\,u_y+u_{xx}+2uu_x=0$ so this 
appears to be an exact
q-form of Burger's equation.$\hfill\bs$
\end{example}
\indent
It is interesting to compare this with Example 3.1.  The difference 
arises from the term
$xdt=qdtx+(q^2-1)dxt$ in the calculus $\gG_{+}$ and this suggests 
that a KP type equation
might arise naturally in the generalized q-plane.  In this direction 
we recall first some
examples of Dimakis-M\"uller-Hoissen (cf. \cite{cxx,c3,c4,d3,d4}).
\begin{example}
Consider a calculus based on $dt^2=dx^2=dxdt+dtdx=0$ with
$[dt,t]=[dx,t]=[dt,x]=0$ and $[dx,x]=\eta dt$.  Assuming the Leibnitz rule
$d(fg)=(df)g+f(dg)$ for functions and $d^2=0$ one obtains ${\bf (C3)}\,\,
df=f_xdx+(f_t+(1/2)\eta f_{xx})dt$.  For a connection $A=wdt+udx$ the zero
curvature condition $F=dA+A^2=0$ leads to ${\bf
(C4)}\,\,(u_t-w_x+(\eta/2)u_{xx}+\eta uu_x=0$ which for $w_x=0$ is a form of
Burger's equation.
$\hfill\bs$
\end{example}
\begin{example}
Next consider $[dt,t]=[dx,t]=[dt,x]=[dy,t]=[dt,y]=[dy,y]=0$ with
$[dx,x]=2bdy$ and $[dx,y]=[dy,x]=3adt$; further assume
$dt^2=dy^2=dtdx+dxdt=dydt+dtdy=dydx+dxdy=0$.  Then ${\bf (C5)}\,\, df=f_xdx+
(f_y+bf_{xx})dy+(f_t+3af_{xy}+abf_{xxx}dt)$.  For $A= vdx+wdt+udy$ one finds
that $dA+A^2=F=0$ implies
\bq\label{01}
u_x=v_y+bv_{xx}+2bvv_x;\,\,w_x=v_t+3av_{xy}+abv_{xxx}+3auv_x+3av(v_y+bv_{xx});
\end{equation}
$$w_y+bw_{xx}=u_t+3au_{xy}+abu_{xxx}+3auu_x-v[2bw_x-3a(u_y+bu_{xx})]$$
Taking e.g. $w_x=(3a/2b)u_y+(3a/2)u_{xx}$ in the last equation to decouple one
arrives at ${\bf (C6)}\,\,\pp_x(u_t-(ab/2)u_{xxx}+3auu_x)=(3a/2b)u_{yy}$;
for suitable $a,b$ this is KP.  One would like to generalize Example 
7.3 in the generalized
q-plane of Example 8.1.
$\hfill\bs$
\end{example}
\indent
\begin{example}
We will try now a somewhat different approach.  First we take 
q-derivatives only in $x$, as in
the case of qKP for example, and use ${\bf 
(C7)}\,\,\pp_q^xf=[f(qx)-f(x)]/(q-1)x$.  We know
from Example 7.3 that {\bf (C5)} leads to interesting consequences 
so begin with an assumption
($f_y=\pp f/\pp y$, etc.)
\bq\label{1c}
df=\pp_q^xfdx+(f_y+b(\pp_q^x)^2f)dy+(f_t+3a\pp_y\pp_q^xf+ab(\pp_q^x)^3f)dt
\end{equation}
Then we can determine what elementary commutation relations between 
the variables are consistent
with \eqref{1c}.  This is rather ad hoc but we stipulate $x,y,t$ 
ordering and then 
\begin{proposition}
There are relations
\bq\label{2c}
dxx=qxdx+b[2]_qdy;\,\,dyx=q^2xdy+a[3]_qdt;\,\,dtx=q^3xdt;\,\,dxy=ydx+3adt
\end{equation}
along with
\bq\label{3c}
[dt,y]=[dy,y]=[dy,t]=[dx,t]=[dt,t]=[dx,y]=0
\end{equation}
which are determined by \eqref{1c}
\end{proposition}
\indent
{\it PROOF.}
The underlying structure for $x,y,t$ is not visible but \eqref{2c} - 
\eqref{3c} do lead to
\eqref{1c} and whatever zero curvature equations subsequently arise. 
To see how \eqref{2c}
emerges from \eqref{1c} consider the analogies to Example 7.3.  Thus
\bq\label{4c}
dx^2=xdx+dxx=(\pp_q^x)x^2dx+b(\pp_q^x)^2x^2dy=[2]_qxdx+b[2]_qdy
\end{equation}
This means ${\bf 
(C8)}\,\,dxx=([2]_q-1)xdx+b[2]_qdy=qxdx+b[(q^2-1)/(q-1)]dy$. 
Similarly
\bq\label{5c}
dx^3=x^2dx+dx^2x
=x^2dx+[2]_q[x(qxdx+b[2]_qdy)]+[2]_qbdyx
\end{equation}
Since $1+q[2]_q=[3]_q$ the $x^2dx$ terms agree.  For the $dy$ and 
$dt$ terms we have
\bq\label{6c}
b[3]_q[2]_qxdy+ab[3]_q[2]_qdt=[2]_q^2bxdy+[2]_qbdyx
\end{equation}
Consequently one takes ${\bf (C9)}\,\,dyx=q^2xdx+a[3]_qdt$ since 
$\eqref{6c}\Rightarrow
dyx=([3]_q-[2]_q)xdy+a[3]_qdt$ and $[3]_q-[2]_q=q^2$.  Next consider
\bq\label{7c}
dx^4=
x^3dx+[3]_qx^2(qxdx+b[2]_qdy)+bx[2]_q[3]_q(q^2xdy+a[3]_qdt)+ab[3]_q[2]_qdtx
\end{equation}
We note that $1+q[3]_q=[4]_q$ and $[2]_q[3]_q(1+q^2)=[3]_q[4]_q$ so 
\eqref{7c} is consistent with
\eqref{1c} provided ${\bf (C10)}\,\,dtx=q^3xt$ (since 
$[3]_q+q^3=[4]_q$).  Further checks with
$dx^5$ etc. hold and the treatment of $dy^m,\,dt^k$ is immediate.  In 
order to deal with the
mixed derivative term in \eqref{1c} one looks at $x^ny^m$ which 
yields a comparison between
\bq\label{8c}
dx^ny^m+x^ndy^m=([n]_qx^{n-1}dx+b[n]_q[n-1]_qx^{n-2}x^{n-2}dy+
\end{equation}
$$+ab[n]_q[n-1]_q[n-2]_qx^{n-3}dy)y^m+x^nmy^{m-1}dy$$
\bq\label{9c}
d(x^ny^m)=[n]_qx^{n-1}y^mdx+(x^nmy^{m-1}+b[n]_q[n-1]_qx^{n-2}y^m)dy+
\end{equation}
$$+(3a[n]_qmx^{n-1}y^{m-1}+ab[n]_q[n-1]_q[n-2]_qx^{n-3}y^m)dt$$
There is a missing $3a[n]_qmx^{n-1}y^{m-1}dt$ term in \eqref{8c} 
which can be provided upon
assuming ${\bf (C11)}\,\,dxy=ydx+3adt$ leading to 
$dxy^m=y^mdx+3my^{m-1}dt$.  Thus the origin
of \eqref{2c} is clear.$\hfill\bs$
\\[3mm]\indent
Returning now to \eqref{1c}
it remains to check now the zero curvature equation arising in the 
spirit of Example 7.3 (some
extra factors and terms will arise via noncommutativity).  Thus 
assume first $dx^2=dy^2=dt^2=0$
and take $A=vdx+wdt+udy$ so
\bq\label{10c}
dA=(v_y+b(\pp_q^x)^2v)dydx+(v_t+3a\pp_y\pp_q^xv+ab(\pp_a^x)^3v)dtdx+\pp_q^xwdxdt+
\end{equation}
$$+(w_y+b(\pp_q^x)^2w)dydt+\pp_q^xudxdy+(u_t+3a\pp_y\pp_q^xu+ab(\pp_q^x)^3u)dtdy$$
For $A^2$ we have $A^2=vdxvdx+ 
vdxwdt+vdxudy+wdtvdx+wdtwdt+wdtudy+udyvdx+udywdt+udyudy$ so one
must use \eqref{2c} in moving the interior differentials to the 
right.  Thus $dxf\sim dx\sum
f_{nmk}x^ny^mt^k$ etc. and in particular
\bq\label{11c}
dtf=\sum f_{nmk}q^{3n}x^ny^mt^kdt=D_xfdt
\end{equation}
A little calculation yields next
\bq\label{111c}
dyx^n=q^{2n}x^ndy+a[3]_qx^{n-1}q^{2(n-1)}[n]_qdt
\end{equation}
Then, given $dyy=ydy$ and $dty=ydt$ with $dyt=tdy$ and $dtt=tdt$, there results
$$dyf=dy\sum f_{nmk}x^ny^mt^k=\sum 
f_{nmk}(q^{2n}x^ndy+a[3]_q[n]_qx^{n-1}q^{2(n-1)}dt)y^mt^k=$$
\bq\label{12c}
=D_x^2fdy+a[3]_qD_x^2\pp_q^xfdt
\end{equation}
Next consider $dx f$ and some routine calculation yields
\bq\label{18c}
dxf=D_xfdx+b[2]_qD_x\pp_q^xfdy+3aD_x\pp_yfdt+ ab[3]_qD_x(\pp_q^x)^2fdt
\end{equation}
Next we put tentatively $(\clubsuit)\,\,
dxdt+dtdx=0=dxdy+dydx=dydt+dtdy$.
Then $dA+A^2=0$ involves
\bq\label{21c}
dtdx:\,\,3avD_x\pp_yv+ 
ab[3]_qvD_x(\pp_q^x)^2v-vD_xw+wD_x^3v+a[3]_quD_x^2\pp_q^xv-
\end{equation}
$$-\pp_q^xw+v_t+3a\pp_y\pp_q^xv+ab(\pp_q^x)^3v=0$$
\bq\label{22c}
dxdy:\,\,vD_xu-uD_x^2v-b[2]_qvD_x\pp_q^xv-v_y-b(\pp_q^x)^2v+\pp_q^xu=0
\end{equation}
and finally
$$dtdy:\,\,3avD_x\pp_yu+ab[3]_qvD_x(\pp_q^x)^2u+wD_x^3u-uD_x^2w-b[2]_qvD_x\pp_q^xw+
a[3]_quD_x^2\pp_q^xu-$$
\bq\label{23c}
-w_y-b(\pp_q^x)^2w+u_t+3a\pp_y\pp_q^xu+ab(\pp_q^x)^3u=0
\end{equation}
Now compare this with Example 7.3 where one has $A=vdx+wdt+udy$ and 
with a little calculation
\bq\label{24c}
dA=[v_y+bv_{xx}-u_x]dydx+[v_t+3av_{xy}+abv_{xxx}-w_x]dtdx+
\end{equation}
$$+[w_y+bw_{xx}-u_t-3au_{xy}
-abu_{xxx}]dydt$$
This is essentially the same as \eqref{10c}.  As for $A^2$ one finds
\bq\label{27c}
A^2=2bvv_xdydx+[3auv_x+3avv_y+3abvv_{xx}]dtdx+
\end{equation}
$$+[2bvw_x-3avu_y-3abvu_{xx}-3auu_x]dydt$$
There are a number of cancellations of the form $uv-vu$ which do not 
cancel in the q-situation
due to factors of $D_x$.  One has the following comparisons
\begin{itemize}
\item
$dxdy:\,\,-2bvu_x\sim vD_xu-uD_x^2v-b[2]_qvD_x\pp_q^xv$
\item
$dtdx:\,\,3auv_x+3avv_y+3abvv_{xx}\sim 
wD_x^3v-vD_xw+a[3]_quD_x^2\pp_q^xv+ 3avD_x\pp_yv
+ab[3]_qvD_x(\pp_q^x)^2v$
\item
$dtdy:\,\,3avu_y+3auu_x+3abvu_{xx}-2bvw_x\sim wD^2_xu-uD_x^3w+3avD_x\pp_yu+
a[3]_quD_x^2\pp_q^xu+ab[3]_qvD_x(\pp_q^x)^2u-b[2]_qvD_x\pp_q^xw$
\end{itemize}
These are all close in the same manner, with distortions due to 
factors of $D_x$ and q-numbers.
One expects now to be able to produce a qKP type equation directly 
from these formulas and
some manipulation as in Example 7.3.  Note that there does not 
seem to be any explicit
formula in the literature for qKP; it is defined via a hierarchy and 
calculations involving
specific equations are difficult, even in the Hirota forms (for qKdV see Section
10).
$\hfill\bs$
\end{example}
\indent
We compare now the KP type equations derivable from Examples 7.3 and 7.4.  First
\bq\label{35c}
v_y-u_x+bv_{xx}+2vvv_x=0\sim 
v_y-\pp_q^xu+b(\pp_q^x)^2v+[2]_qbvD_x\pp_q^xv+uD_x^2v-vD_xu=0
\end{equation}
\bq\label{36c}
v_t+3av_{xy}+abv_{xxx}-w_x+3auv_x+3avv_y+3abvv_{xx}=0\sim
\end{equation}
$$v_t+3a\pp_y\pp_q^xv+ab(\pp_q^x)^3v-\pp_q^xw+[3]_qauD_x^2\pp_q^xv+3avD_x\pp_yv+
[3]_qabvD_x(\pp_q^x)^2v=0$$
\bq\label{37c}
u_t+3au_{xy}+abu_{xxx}-bw_{xx}-w_y+3avu_y+3auu_x+3abvu_{xx}-2bvw_x=0\sim
\end{equation}
$$u_t+3a\pp_y\pp_q^xu+ab(\pp_q^x)^3u-b(\pp_q^x)^2w-w_y+wD_x^3u-uD_x^2w+3avD_x\pp_yu+$$
$$+[3]_qauD_x^2\pp_q^xu+[3]_qabvD_x(\pp_q^x)^2u-b[2]_qv\pp_q^xw=0$$
The similarities are obvious.  Now look at the derivation of KP in 
Example 7.3.  One takes
$w_x=(3a/2b)u_y+(3a/2)u_{xx}$ to decouple in \eqref{37c} which 
reduces \eqref{37c} to
$(\bl)\,\,u_t+3au_{xy}+abu_{xxx}-bw_{xx}-w_y=0$.  Now set
\bq\label{38c}
w=\pp^{-1}\left(\frac{3a}{2b}u_y\right)+\frac{3a}{2}u_x
\end{equation}
to obtain
\bq\label{39c}
u_t+3au_{xy}+abu_{xxx}+3auu_x-b\left[\frac{3a}{2b}u_{yx}+\left(\frac{3a}{2}\right)u_{xxx}\right]-
\end{equation}
$$-\pp^{-1}\left[\frac{3a}{2b}u_{yy}\right]-\frac{3a}{2}u_{xy}=u_t-\frac{ab}{2}u_{xxx}
-\frac{3a}{2b}\pp^{-1}u_{yy}+3auu_x=0$$
\indent
Repeating such a procedure for the q-situation 
\eqref{35c}-\eqref{37c} leads to
\begin{theorem}
For suitable $a,b$ one can derive the KP equation by a zero curvature 
condition in the form
$\eqref{39c}\,\,u_t-(ab/2)u_{xxx}-(3a/2b)\pp^{-1}u_{yy}+3auu_x=0$. 
The same procedure applied
to zero curvature equations as in Example 7.4 leads to a q-version of 
KP in the form
\bq\label{48c}
u_t+ab(\pp_q^x)^3u-\frac{[3]_qab}{[2]_q}\pp_q^xD_x(\pp_q^x)^2u+[3]_qauD_x^2\pp_q^xu
-\frac{3a}{[2]_qb}(\pp_q^x)^{-1}D_xu_{yy}=A(u,w)+B(u)
\end{equation}
where 
\bq\label{49c}
w=\frac{3a}{[2]_qb}(\pp_q^x)^{-1}D_x\pp_yu+\frac{[3]_qa}{[2]_qq}D_x\pp_q^xu;
\end{equation}
$$A(u,w)=wD_x^3u-uD_x^2w;\,\,B(u)=
\frac{[3]_qa}{[2]_qq}D_x\pp_y\pp_q^xu+\frac{3a}{[2]_q}\pp_q^x
D_x\pp_yu-3a\pp_y\pp_q^xu$$
Then $A\to 0$ and $B\to 0$ as $q\to 1$ and the equation \eqref{49c} 
goes to the standard
KP form.  Since this equation was derived from exactly the same 
procedures as the classical KP
equation (via zero curvature considerations) one expects it to be a 
good candidate for a qKP
equation.  Obtaining such a qKP equation from  qKP hierarchy 
equations or corresponding Hirota
formulas involves considerable calculation and the resulting equation will not have
a finite closed form (cf. Section 10 for qKdV).
\end{theorem}
\indent
{\bf REMARK 7.1.}
An interesting zero curvature formulation for KdV goes back to \cite{cz,dz}
and should have a quantum group adaptation.  Thus 
in the version of \cite{dz} one considers $SL(2,{\bf
R})$ with Lie algebra generators
\bq\label{1z}
T_1=ih=i\left(\begin{array}{cc}
1 & 0\\
0 & -1
\end{array}\right);\,\,T_2=if=i\left(\begin{array}{cc}
0 & 0\\
1 & 0
\end{array}\right);\,\,T_3=ie=i\left(\begin{array}{cc}
0 & 1\\
0 & 0
\end{array}\right)
\end{equation}
with ${\bf (D1)}\,\,[T_a,T_b]=iC^c_{ab}T_c$ or
$[T_1,T_2]=-2iT_2,\,\,[T_1,T_3]=2iT_3,$ and $[T_1,T_2]=iT_1$.  Writing
$g(x,t)=exp(i\gt^a(x,t)T_a)$ one sets $A_{\mu}=g^{-1}\pp_{\mu}g=
i\sum A_{\mu}^a(x,t)T_a\,\,(x^0=t,\,x^1=x)$ to obtain Maurer-Cartan equations
${\bf (D2)}\,\,\pp_{\mu}A_{\nu}^a-\pp_{\nu}A_{\mu}^a-C^a_{bc}A_{\mu}^b
A_{\nu}^c=0$ (which is equivalent to
$\pp_{\mu}A_{\nu}-\pp_{\nu}A_{\mu}+[A_{\mu},A_{\nu}]=0$).  
Note the $\gt^a(x,t)$ are ``arbitrary" functions of $x,t$ with
$i\sum\gt^aT_a\in{\mf s}\ell_2({\bf C})$.  Writing
then ${\bf (D3)}\,\,A_1^1=\sqrt{-\gl}\,\,(\gl>0)$ and $A_1^3=6$ there results
for $a=3$
\bq\label{2z}
A_0^1=\left(\frac{\sqrt{-\gl}}{6}A_0^3-\frac{1}{12}A^3_{0,x}\right)
\end{equation}
Similarly for $a=1,2$ one has (using \eqref{2z})
\bq\label{3z}
A_0^2=\left(\frac{\sqrt{-\gl}}{36}A^3_{0,x}-\frac{1}{72}A^3_{0,xx}+\frac
{1}{6}A_0^3A_1^2\right)
\end{equation}
and (using \eqref{2z}, \eqref{3z})
\bq\label{4z}
A^2_{1,t}=-\frac{1}{72}A^3_{0,xxx}+\frac{1}{3}A^3_{0,x}A_1^2+
\frac{1}{6}A_0^3A_{1,x}^2-\frac{\gl}{18}A^3_{0,x}
\end{equation}
Taking $A_1^2=-(1/36)u$ with $A^3_0=A(u(x,t))$ one obtains
\bq\label{5z} 
u_t=\frac{1}{2}A_{xxx}+\frac{1}{3}A_xu+\frac{1}{6}u_xA+2\gl A_x\equiv
\end{equation}
$$\equiv u_t=\frac{1}{2}\left[D^3+\frac{1}{3}(Du+uD)\right]A+2\gl DA$$
which is the Lenard form for KdV.  Indeed since $u_t$ is independent of 
$\gl$ one must choose $A(u)$ to be a function of $\gl$ such that 
\eqref{5z} is consistent.  It is appropriate to take
\bq\label{6z}
A(u)=2\sum_0^nA_j(u)(-4\gl)^{n-j}
\end{equation}
leading to a recursion $u_t=[D^3+(1/3)(Du+uD)A_n]$ with
\bq\label{7z} D^3+\frac{1}{3}(Du+uD)A_j=DA_{j+1}\,\,(j=0,\cdots,n-1);\,\,
A_0=1
\end{equation}
One can also set $A_j=\gd H_j/\gd u(x)$ so that
\bq\label{8z}
\left[D^3+\frac{1}{3}(Du+uD)\right]\frac{\gd H_m}{\gd u(x)}=D\frac{\gd
H_{m+1}}{\gd u(x)}
\end{equation}
There is more on this below based directly on the method of \cite{cz}.$\hfill\bs$

\section{DIFFERENTIAL CALCULI}
\renewcommand{\theequation}{8.\arabic{equation}}
\setcounter{equation}{0}

Now go to \cite{cxx,k1,m1,mz} and recall that the algebraic properties of 
$SL(2,{\bf C})$ are stored in the coordinate Hopf algebra ${\mf O}(SL(2))=
{\bf C}[u^1_1,u_2^1,u_1^2,u_2^2]/\{u^1_1u_2^2-u_2^1u_1^2=1\}$ (recall
$u^i_j(( g_{k\ell}))=g_{ij}$).  The Hopf algebras ${\mf O}(SL_q(2))$ will be
one parameter deformations of ${\mf O}(SL(2))$.  One defines first ${\mf
O}(M_q (2))$ to be the complex (associative) algebra with generators
$a,b,c,d$ satisfying
\bq\label{9z}
ab=qba;\,\,ac=qca;\,\,bd=qdb;\,\,cd=qdc;\,\,bc=cb;\,\,ad-da=(q-q^{-1})bc
\end{equation}
where $u_1^1=a,\,\,u_2^1=b,\,\,u_1^2=c$, and $u^2_2=d$.  We omit Hopf 
algebra structure momentarily.  One defines now ${\bf
(D4)}\,\,D_q=ad-qc=da-q^{-1}bc$ (quantum determinant) which belongs to the
center of ${\mf O}(M_q(2))$.  Then define
${\mf O}(SL_q(2))={\mf O}(M_q(2))/(D_q-1)$ as the coordinate algebra of
$SL_q(2))$.  There is a compact real form ${\mf O}(SU_q(2))$ for $q\in {\bf
R}$ with involution $a^*=d,\,\,b^*=-qc,\,\,c^*=-q^{-1}b,$ and $d^*=a$ so that
${\mf O}(SU_q(2))$ is the star algebra generated by two elements $a,c$ with
relations
\bq\label{10z}
ac=qca;\,\,ac^*=qc^*a;\,\,cc^*=c^*c;\,\,a^*a+c^*c=1;\,\,aa^*+q^2c^*c=1
\end{equation}
For $|q|=1$ there is also a real form ${\mf O}(SL_q(2,{\bf R}))$ of
${\mf O}(SL_q(2))$ with $a^*=a,\,\,b^*=b,\,\,c^*=c,$ and $d^*=d$.  Note that 
$SL_q(2)$ is known through its coordinate algebra but one can also envision
matrices 
\bq\label{11z}
M_m^n=\left(\begin{array}{cc}
a & b\\
c & d
\end{array}\right);\,\,(M_m^n)^{\dg}=(M_m^n)^{-1}
\end{equation}
for suitable $a,b,c,d$ (cf. \cite{mz}). 
\\[3mm]\indent
Next we recall the following universal enveloping algebras
\begin{example}
In the terminology of \cite{m1} ${\mf s}\ell_2$ has generators 
$X_{\pm},\,H$ with ${\bf (D5)}\,\,[H,X_{\pm}]=\pm 2X_{\pm}$ and
$[X_{+},X_{-}]=H$.  ${\mf U}_q({\mf s}\ell_2)$ is then defined via generators
$1,\,X_{\pm}, q^{H/2},$ and $q^{-H/2}$ with
\bq\label{12z}
q^{H/2}X_{\pm}q^{-H/2}=q^{\pm
1}X_{\pm};\,\,[X_{+},X_{-}]=\frac{q^H-q^{-H}}{q-q^{-1}}
\end{equation}
(we omit the Hopf algebra structure for now).  For real forms one has
${\mf U}_q({\mf su}_2)$ defined as ${\mf U}_q({\mf s}\ell_2)$ for 
$q\in {\bf R}$ with ${\bf (D6)}\,\,H^*=H$ and $X_{\pm}^*=X_{\mp}$ while for
${\mf U}_q({\mf s}\ell(2,{\bf R}))$ one takes $|q|=1$ with ${\bf
(D7)}\,\,H^*=-H$ and $X_{\pm}^*=-X_{\pm}$.$\hfill\bs$
\end{example}
\begin{example}
Another variation appears in \cite{k1}.  Thus for $q\ne 0$ and $q^2\ne 1$
${\mf U}_q({\mf s}\ell_2)$ is the associative algebra with unit generated by 
$E,\,F,\,K,\,K^{-1}$ satisfying
\bq\label{13z}
KK^{-1}=K^{-1}K=1;\,\,KEK^{-1}=q^2e;\,\,KFK^{-1}=q^{-2}F:\,\,[E,F]=\frac
{K-K^{-1}}{q-q^{-1}}
\end{equation}
Then the sets $\{F^{\ell}K^mE^n;\,\,m\in {\bf Z},\,\ell,n\in{\bf N}_0\}$ or
$\{E^nK^mF^{\ell};\,\,m\in{\bf Z},\,\ell,n\in{\bf N}_0\}$ are vector space
bases of ${\mf U}_q({\mf s}\ell_2)$ (again we omit temporarily the Hopf
algebra structure).  For $q\in {\bf R},\,\,{\mf U}_q({\mf su}_2)$ arises via
${\bf (D8)}\,\,E^*=FK,\,\,F^*=K^{-1}E,$ and $K^*=K$.  For $|q|=1,\,\,
{\mf U}_q({\mf s}\ell_2({\bf R}))$ is given via ${\bf (D9)}\,\,
E^*=-E,\,\,F^*=-F,$ and $K^*=K$.$\hfill\bs$
\end{example}
\indent
Now one knows that for $a\in {\mf U}({\mf g})$ and $f\in{\mf O}(G)$
there is a dual pairing determined by ${\bf (D10)}\,\,<a,f>=(\tl{a}f)(e)$
where for $a=X_1\cdots X_n$
\bq\label{14z}
(\tl{a}f)(g)=\left.\frac{\pp^n}{\pp t_1\cdots\pp t_n}\right|_{t=0}
f(gexp(t_1X_1)\cdots exp(t_nX_n))
\end{equation}
Then for $q^4\ne 0,1$ let $\breve{{\mf U}}_q({\mf s}\ell_2)$ be the algebra
over ${\bf C}$ with generators $E,\,F,\,K,\,K^{-1}$ satisfying
\bq\label{15z}
KK^{-1}=K^{-1}K=1;\,\,KEK^{-1}=qE;\,\,KFK^{-1}=q^{-1}F;\,\,[E,F]=\frac
{K^2-K^{-2}}{q-q^{-1}}
\end{equation}
One can show that there is an injective Hopf algebra homomorphism $\phi:\,\,
{\mf U}_q({\mf s}\ell_2)\to\breve{{\mf U}}_q({\mf s}\ell_2)$ determined via
$\phi(E)=EK,\,\,\phi(F)=K^{-1}F,$ and $\phi(K)=K^2$.  Under this injection
${\mf U}_q({\mf s}\ell_2)$ is a Hopf subalgebra of $\breve{{\mf U}}_q({\mf
s}\ell_2)$.  Now for ${\mf g}={\mf s}\ell(2,{\bf C})$ and $G= SL(2,{\bf
C})$ with generators $H,E,F$ for ${\mf s}\ell(2,{\bf C})$ satisfying
$[H,E]=2E,\,\,[H,F]=-2F,$ and $[E,F]=H$ the dual pairing of \eqref{14z} is
expressed via
\bq\label{16z}
-<H,a>=<H,d>=<E,c>=<F,b>=1
\end{equation}
and zero otherwise.  There is a corresponding result with bracket 
$<\,,\,>\breve{{}}$ for $\breve{{\mf U}}_q({\mf s}\ell_2)$ and ${\mf
O}(SL_q(2))$ such that
\bq\label{17z}
<K,a>\breve{{}}=q^{-1/2};\,\,<K,d>\breve{{}}=q^{1/2};\,\,
<E,c>\breve{{}}=<F,b>\breve{{}}=1
\end{equation}
and all other brackets are zero.  Similarly there is a unique dual pairing 
$<\,,\,>$ for ${\mf U}_q({\mf s}\ell_2)$ and ${\mf O}(SL_q(2))$ such that
the only nonvanishing brackets are
\bq\label{18z}
<K,a>=q^{-1};\,\,<K,d>=q;\,\,<E,c>=<F,b>=1
\end{equation}
and if q is not a root of unity both pairings are nondegenerate.
\\[3mm]\indent
We consider now the quantum tangent space for $A={\mf O}(SL_q(2))$ when
$q^2\ne 1$.  Let $a\sim u_1^1,\,\,b\sim u_2^1,\,\,c\sim u_1^2,$ and $d\sim
u_2^2$ with $E,F,K$ generators of $\breve{{\mf U}}({\mf s}\ell_2)$.
We use the dual pairing \eqref{17z} and define three linear functionals on A
via
\bq\label{19z}
X_0=q^{-1/2}FK;\,\,X_2=q^{1/2}EK;\,\,X_1=(1-q^{-2})^{-1}(\gep-K^4)
\end{equation}
where the Hopf algebra structure is given by
\bq\label{20z}
\gD(K)=K\ot K;\,\,\gD(E)=E\ot K+K^{-1}\ot E;\,\,\gD(F)=F\ot K+K^{-1}\ot F;
\end{equation}
$$S(K)=K^{-1};\,\,S(E)=-qE;\,\,S(F)=-q^{-1}F;
\,\,\gep(K)=1;\,\,\gep(E)=\gep(F)=0$$
This illustrates the advantage of $\breve{{\mf U}}({\mf s}\ell_2)$ over
${\mf U}_q({\mf s}\ell_2)$ in that $\gD(E)$ and $\gD(F)$ have the same type
of formula.  Now one gets
\bq\label{21z}
\gD(X_j)=\gep\ot  X_j+X_j\ot K^2\,\,(j=0,2);\,\,\gD(X_1)=\gep\ot X_1+X_1\ot
K^4
\end{equation}
It follows that $T=Lin\{X_0,X_1,X_2\}$ is the quantum tangent space of a 
(unique up to isomorphism) left covariant FODC $\gG$ over A (3-D calculus
of Woronowicz).
We recall here that an FODC $\gG$ over A is left covariant if and only if
there is a linear map ${\bf (D11)}\,\,\gD_L:\,\gG\to A\ot\gG$ such that
$\gD_L(adb)=\gD(a)(id\ot d)\gD(b)$.  For such $\gG$ the quantum tangent space
is ${\bf (D12)}\,\,T_{\gG}=\{X\in A';\,X(1)=0;\,X(a)=0\,\, for\,\, a\in
R_{\gG}=(ker(\gep))^2=\{f\in A;\, f(e)=(df)(e)=0\}\}$.  One writes 
${}_{inv}\gG=\{\rho\in\gG;\,\gD_L(\rho)=1\ot\rho\}$ 
(left invariant elements) and $dim{}_{inv}\gG=
dim(T_{\gG})=dim\,ker(\gep/R_{\gG})$.  For $P_L:\,\gG\to{}_{inv}\gG$
defined via $P_L(a\rho)=\gep(a)P_L(\rho)$ with $P_L(\rho)=\sum
S(\rho_{(-1)}\rho_0$ where $\gD_L\rho=\sum\rho_{(-1)}\ot\rho_0$ 
one defines ${\bf
(D13)}\,\,\go_{\gG}:\,A\to{}_{inv}\gG$ by $\go_{\gG}(a)=P_L(da)$ so
$\go(A)={}_{inv}\gG$ (see Remark 8.5 for more on $P_L$).  
In particular $R_{\gG}=\{a\in ker(\gep);\,
\go_{\gG}(a)=0\}$ and one takes $X_i$ as a basis of $T_{\gG}$ with $\go_j$
the dual basis of ${}_{inv}\gG$ (i.e. $(X_i,\go_j)=\gd_{ij}$ where the
nondegenerate pairing $(\,,\,):\,T\times\gG\to{\bf C}$ is defined via
${\bf (D14)}\,\,(X,a\cdot db)=\gep(a)X(b)$.  Now let
$\{\go_0,\go_1,\go_2\}$ be a basis of ${}_{inv}\gG$ which is dual to the
basis $\{X_0,X_1,X_2\}$ of $T_{\gG}$.  Then 
\bq\label{1x}
dx=\sum_0^2(X_j\cdot x)\go_j
\end{equation}
(recall $da=\sum (X_j\cdot a)\go_j$ and $X_i(ab)=\gep(a)X_i(b)+
\sum X_j(a)f_i^j(b)$ - $L_{X_i}(a)=X_i\cdot a\sim$ quantum analogue
of a left invariant vector field in Lie theory and there is more on
$X_j\cdot a$ below).
Since $X_0(b)=X_2(c)=X_1(a)=1$ and
$X_0(a)=X_0(c)=X_2(a)=X_2(b)=X_1(b)=X_2(c)=0$ (due to \eqref{17z}) one
obtains
\bq\label{2x}
\go_0=\go(b),\,\,\go_2=\go(c);\,\,\go_1=\go(a)=-q^{-2}\go(d);\,\,
da=b\go_2+a\go_1;
\end{equation}
$$db=a\go_0-q^2b\go_1;\,\,dc=c\go_1+d\go_2;\,\,dd=-q^2d\go_1+c\go_0$$
Comparing \eqref{21z} with the standard prescription 
${\bf (D15)}\,\,\gD(X_i)=\gep\ot X_i+X_k\ot f_i^k$ where
$\gD(f^r_s)=f^r_k\ot f^k_s$ one sees that ${\bf (D16)}\,\,f^j_j=K^2$ for
$j=0,2$ and $f_1^1=K^4$ with the other $f^r_s=0$.  This leads to
\bq\label{3x}
\go_ja=q^{-1}a\go_j;\,\,\go_jb=qb\go_j;\,\,\go_jc=q^{-1}c\go_j\,\,(j=0,2);
\end{equation}
$$\go_1a=q^{-2}a\go_1;\,\,\go_1b=q^2b\go_1;\,\,\go_1c=q^{-2}c\go_1;
\,\,\go_1d=q^2d\go_1$$
One shows now that $R=R_{\gG} = B\cdot A$ where B is the vector space
spanned by ${\bf (D17)}\,\,a+q^{-2}d-(1+q^{-2})\cdot 1,\,b^2,\,c^2,\,bc,\,
(a-1)b,$ and $(a-1)c$ so R is the right ideal of $ker(\gep)$ associated with
the FODC $\gG$.  These six generators satisfy also $S(x)^*\in R$ for ${\mf O}
(SL_q(2,{\bf R})\,\,(|q|=1)$ and ${\mf O}(SU_q(2))\,\,(q\in {\bf R})$ so
$\gG$ is a star calculus for these Hopf algebras.  It follows from the
commutativity relations in ${\mf U}_q({\mf s}\ell_2)$ that the basis
elements of $T_{\gG}$ satisfy
\bq\label{4x}
q^2X_1X_0-q^{-2}X_0X_1=(1+q^2)X_0;\,\,q^2X_2X_1-q^{-2}X_1X_2=(1+q^2)X_2;
\end{equation}
$$qX_2X_0-q^{-1}X_0X_2=-q^{-1}X_1$$
\indent
A next step could be to determine the associated bicovariant FODC and
develop corresponding q-Lie algebraic ideas (i.e. $T_{\gG}$ can be made into
a Lie algebra - cf. \cite{k1}, pp. 498-505).  First however one notes that
if $\gG^{\wg}=\oplus_0^{\infty}\gG^{\wg n}$ is a DC over A based on 
$\gG^{\wg 1}=\gG$ then
\bq\label{5x}
d\go(a)=-\sum \go(a_1)\ot \go(a_2)
\end{equation}
is called the Maurer-Cartan formula.  This can be written via $\go(a)=
\sum X_i(a)\go_i$ in the form
\bq\label{6x}
\sum X_i(a)d\go_i=-\sum X_iX_j(a)\go_i\wg \go_j
\end{equation}
Now one can determine a universal DC $\gG^{\wg}$ for any left covariant FODC
$\gG$ (see \cite{k1}, p. 509) and for the $3-D$ calculi on $SL_q(2)$ with
$A={\mf O}(SL_q(2))$ one has ($q^2\ne -1$) the following construction. 
First one constructs $S(x)$ for the six generators of R in {\bf (D17)}
and shows that $S(B)$ is spanned by elements
\bq\label{7x}
\go_0^2;\,\,\go_1^2;\,\,\go_2^2;\,\,\go_2\go_0+q^2\go_0\go_2;\,\,
\go_1\go_0+q^4\go_0\go_1;\,\,\go_2\go_1+q^4\go_1\go_2
\end{equation}
of the tensor algebra $T_A(\gG)$ over $\gG$.  This leads to (see \cite{k1}
for details)
\bq\label{8x}
\go_j\wg\go_j=0\,\,(j=0,1,2);\,\,\go_2\wg\go_0=-q^2\go_0\wg\go_2;\,\,
\go_1\wg\go_0=-q^4\go_0\wg\go_1;
\end{equation}
$$\go_2\wg\go_1=-q^4\go_1\wg\go_2$$
in $\gG^{\wg 2}$.  The 3-form $\go_0\wg\go_1\wg\go_2$ is a free left A module
basis of $\gG^{\wg 3}$ and $\gG^{\wg n}=0$ for $n\geq 4$.
\\[3mm]\indent
{\bf REMARK 8.1.}
We must spell out a little more detail regarding $X_i\cdot a$.  thus $\go$
refers to $\go_{\gG}$ and is a map $A\to {}_{inv}\gG$.  Further $\gG=
A\go(A)=\go(A)a$.  Since $\gD_L(da)=\sum a_1\ot da_2$ one has ${\bf
(D18)}\,\,\go(a)=\sum S(a_1)da_2$ and $da=\sum a_1\go(a_2)$ (see below). 
Now
$(X_i,\go(a))=X_i(a)$ via $(\go(a),X)=X(a)$ so ${\bf (D19)}\,\,
\go(a)=\sum X_i(a)\go_i$ and $da=\sum
a_1\go(a_2)=\sum\sum_ia_1X_i(a_2)\go_i=\sum_i(X_i\cdot a)\go_i$.  This
implies ${\bf (D20)}\,\,X_i\cdot a=\sum a_1X_i(a_2)$.  For bicovariant $\gG$
another variation on the Maurer-Cartan formulas exists.  Thus take
$a_j\in A$ such that $X_i(a_j)=\gd_{ij}$ and set $c_{ij}^k=X_iX_j(a_k)$.  Put
$a=a_k$ into \eqref{6x} and apply
\bq\label{9x}
[X_i,X_j]=X_iX_j-\gs^{ij}_{nm}X_nX_m=C^k_{ij}X_k
\end{equation}
where $[X,Y](a)=(X\ot Y)(Ad_R(a))$ (a left covariant FODC over A with
associated right ideal $R_{\gG}$ is bicovariant if and only if
$Ad_R(R_{\gG})\subseteq R_{\gG}\ot A$ - right coadjoint action
$Ad_R(a)=\sum a_2\ot S(a_1)a_3$).  Here
${\bf (D21)}\,\,[X_i,X_j]=X_iX_j-\gs^{ij}_{nm}X_nX_m=C^k_{ij}X_k$ where
$\gs^{nm}_{ij}=f_m^i(v_j^n)$ and $C^k_{ij}=X_j(v^i_k)$ ($\gs$ is the
braiding of $\gG\ot_A\gG$ with
$\gs(\go_i^1\ot_A\go_j^2)=\gs^{nk}_{ij}\go_n^2\ot_A\go_k^1$, $f$ is in {\bf
(D15)}, and ${\bf (D22)}\,\,\sum\gD_R(a_1X_m(a_2)\go_m)=\sum a_1X_m(a_3)\go_i\ot
a_2v^i_m$ - cf. Remark 8.3).  Consequently one obtains
\bq\label{10x}
d\go_k=-c_{ij}^k\go_i\ot\go_j;\,\,c^k_{ij}-c^k_{nm}\gs^{ij}_{nm}=C^k_{ij}
\end{equation}
We refer to \cite{k1} for perspective and details (cf. also \cite{de} for quantum Lie
algebras).
$\hfill\bs$
\\[3mm]\indent
The $3-D$ calculus constructed above is not bicovariant (cf. \cite{k1}) but
there are associated bicovariant $4-D_{\pm}$ calculi $\gG_{\pm}$ on $SL_q(2)$ which
we do not deal with here.
We now try to produce a quantum version of \eqref{1z} - \eqref{8z}.  First the
spirit of {\bf (D1)} - {\bf (D2)} lies in defining a matrix of left
invariant one forms $\gO=B^{-1}dB\sim dB=B\gO$ so that $d(dB)=0=dB\wg \gO
+Bd\gO=B(d\gO+\gO\wg\gO)$ leading to $d\gO+\gO\wg\gO=0$.  For $\gO=
\sum X_i\go_i$ one arrives at
\bq\label{14x}
\gO\wg\gO=\frac{1}{2}\sum [X_k,X_m]\go_k\wg\go_m=\frac{1}{2}\sum
C^{\ell}_{km}X_{\ell}\go_k\wg\go_m=-\sum X_{\ell}d\go_{\ell}
\end{equation}
which implies ${\bf (D23)}\,\,d\go_{\ell}=-(1/2)\sum C^{\ell}_{km}
\go_k\wg\go_m$.  We note that {\bf (D2)} is immediate via e.g.
${\bf
(D24)}\,\,\pp_{\nu}A_{\mu}=g^{-1}\pp_{\nu}gg^{-1}\pp_{\mu}g+g^{-1}\pp_{\nu}
\pp_{\mu}g=A_{\nu}A_{\mu}+g^{-1}\pp_{\nu}\pp_{\mu}g$.  The variables
$x^0,x^1\sim t,x$ are simply thrown in ad hoc; they do not themselves
represent natural Lie algebra type derivations.
In fact, recalling $g=exp(i\gt^a(x,t)T_a)$ with $A_{\mu}=g^{-1}\pp_{\mu}g$,
we see that the $\gt^a$ are simply numerical coefficients.  There is no
standard exponential map for quantum groups (cf. however \cite{m1}) and
we want an analogue for say $A_{\mu}$.  We know $T_{\gG}$ is based on
$X_1,\cdots,X_4$.  The FODC provides us with a differential calculus
according to strict rules, which we recall now.  First
$\go=\go_{\gG}:\,\,A\to {}_{inv}\gG$ is defined via $\go(a)=P_L(da)$ where
$P_L(a\rho)=\gep(a)P_L(\rho)$.  Recall for $\gD_L:\,\gG\to A\ot \gG$ one
writes $\gD_L(a)=\sum a_{(-1)}\ot a_0$ instead of $\gD(a)=\sum a_1\ot a_2$
and it is required that ${\bf (D25)}\,\,\gD_L(a\rho
b)=\gD(a)\gD_L(\rho)\gD(b)$.  Further $P_L(\rho)=\sum S(\rho_{(-1)})\rho_0$
and $\rho$ is left invariant if $\gD_L(\rho)=1\ot\rho$.  
\\[3mm]\indent
{\bf REMARK 8.2.}
It will be expedient here to review some basic ideas from Hopf algebra. 
Thus a Hopf algebra A over a field D with invertible antipode S involves
$(A,\cdot,\eta,\gD,\gep,S)$ where $\eta:\,K\to A$ is linear, $\gD:\,
A\to A\ot A$ satisfies ${\bf (D26)}\,\,(\gD\ot id)\ci\gD=(id\ot\gD)\ci\gD$
and $a=(\gep\ot id)\ci\gD(a)=(id\ot\gep)\ci\gD(a)$ ($\gep:\,A\to K$ is the
counit).  Also $\eta(\gl)a=\gl\gep(a)$ and ${\bf (D27)}\,\,\cdot(id\ot
\cdot)=\cdot(\cdot\ot id)$ with $\cdot (\eta\ot id)=id=\cdot(id\ot\eta)$.
One has also ${\bf (D28)}\,\,\gD(ab)=\gD(a)\gD(b),\,\,\gD(1)=1\ot
1,\,\,\gep(ab)=\gep(a)\gep(b),$ and $\gep(1)=1$.  Further the antipode
$S:\,A\to A$ satisfies ${\bf (D29)}\,\,\cdot (S\ot id)\ci\gD=\cdot(id\ot
S)\ci\gD=\eta\ci\gep$ so $\gep(a)=\sum a_1S(a_2)=\sum S(a_1)a_2$
and $a=\sum a_1\gep(a_2)=\sum \gep(a_1)a_2$.  One
notes also ${\bf (D30)}\,\,\gD\ci S=\tau\ci (S\ot S)\ci\gD$ where $\tau(v\ot
w)=w\ot v$ ($A\ot B$ has the natural structure $(a\ot b)(a'\ot b')=aa'\ot
bb'$).$\hfill\bs$
\\[3mm]
\indent
In working with covariant FODC in this general context one must be more
careful with the Hopf algebra constructions.  Resorting to simple $dx,\,dy$
formulas could be very risky.  Then from {\bf (D25)} one has 
${\bf (D31)}\,\,\gD_L(da)=(id\ot d)\gD(a)=\sum a_1\ot da_2$.  Recall also
$\go(a)=P_L(da),\,\,P_L(a\rho)=\gep(a)P_L(\rho),\,\,P_L(\rho)=\sum
S(\rho_{(-1)})\rho_0,$ and $P_L(\rho a)=\sum
S(a_1)P_L(\rho)a_2=ad_R(a)P_L(\rho)$.  Since $\sum \rho_{(-1)}\ot\rho_0\sim
\gD_L(da)=\sum a_1\ot da_2$ one has ${\bf (D32)}\,\,\go(a)=\sum S(a_1)da_2$
and $da=\sum a_1\go(a_2)$ (since $\rho=\sum \rho_{(-1)}P_L(\rho_0)$ implies
$da=\sum a_1P_L(da_2)=\sum a_1\go(a_2)$).
\\[3mm]\indent
{\bf REMARK 8.3.}
The other point to be developed is how to work realistically with FODC $\gG$
when they are given in the abstract context indicated.  Generally speaking
for an FODC one has an A bimodule $\gG$ over an algebra A with a linear 
$d:\,A\to\gG$ such that d satisfies Leibnitz and $\gG$ is the linear span of
elements $xdyz$.  In particular if d is a derivation of A with values in the
bimodule $\gG_0$ then $\gG=A\cdot dA\cdot A$ is an FODC over A.  Evidently
$dxy$ and $ydx$ are generally different (even when A is commutative).  One
often thinks of an FODC over a quantum space X with left coaction
$\phi:\,X\to A\ot X$ (i.e. X is a left A-comodule algebra).  Then an FODC
$\gG$ over X involves a left coaction $\Phi:\,\gG\to A\ot \gG$ such that
$\Phi(x\rho y)=\phi(x)\Phi(\rho)\phi(y)$ and $\Phi(dx)=(id\ot d)\phi(x)$.
If $\{x_i\}$ are a linearly independent set of generators of the algebra X
then there is a unique left covariant FODC $\gG$ over X such that $\{dx_i\}$
are a free right X module basis and the bimodule structure is given via
${\bf (D33)}\,\,x_i\cdot dx_j=dx_k\cdot x_{ij}^k$ (cf. \cite{k1}).  In this
situation for any $x\in X$ there are uniquely determined elements
$\pp_i(x)\in X$ such that ${\bf (D34)}\,\,dx=\sum dx_i\cdot\pp_i(x)$.  We
recall also that a bicovariant bimodule $\gG$ over A involves linear maps
$\gD_L:\,\gG\to A\ot \gG$ and $\gD_R:\,\gG\to \gG\ot A$ such that
$(id\ot \gD_R)\ci\gD_L=(\gD_L\ot id)\ci\gD_R$.  If $\gG$ is a bicovariant
bimodule over A and $\go_i$ is a basis of ${}_{inv}\gG$ there exist
pointwise finite matrices $v=(v^i_j)$ and $f=(f^i_j)$ with $v_j^i\in A$ and
$f^i_j\in A'$ such that 
\begin{itemize}
\item
$\go_ia=(f^i_k\cdot a)\go_k$
\item
$\gD_R(\go_i)=\go_k\ot v^k_i$
\item
$f^i_j(ab)=f^i_k(a)f^k_j(b)$ and $f^i_j(1)=\gd_{ij}$
\item
$\gD(v^i_j)=v^i_k\ot v^k_j$ and $\gep(v^i_j)=\gd_{ij}$
\item
$\sum_kv_i^k(a\cdot f_j^k)=\sum_k(f^j_k\cdot a)v^j_k$
\end{itemize}
Moreover $\eta_i=\go_jS(v^j_i)$ is a basis
of $\gG_{inv}$ (right invariant forms from $\gG$) and both sets $\go_i$ and
$\eta_j$ are free left (resp. right) A module basis of $\gG$.$\hfill\bs$
\\[3mm]\indent
Now go back to $A={\mf O}(SL_q(2))$ with $X_i\in A'\,\,(i=1,\cdots,4)$ 
and $T_{\pm}\subset A'$.  Recall $E,F,K$ generate $\breve{{\mf U}}_q({\mf
s}\ell_2)$ and the bracket $(\,,\,):\,T\times \gG\to {\bf C}$ is defined via
$(X,adb)=\gep(a)X(b)$ (this gives a nondegenerate dual pairing between
${}_{inv}\gG=\go(A)$ and $T$ with $\gG\sim AdA$ and $(\go(a),X)=X(a)$).
For $\gz=\sum a_idb^i\in \gG$ and $X\in T$ one has
$(X,\gz)=X(\sum\gep(a_i)b^i)$ and this is independent of the representation
of $\gz$.  A few other notions to recall involve (cf. Remark 8.1)
$\go(a)=\sum S(a_1)da_2$ and 
$(X_i,\go(a))=X_i(a)$ with $\go(a)=\sum X_i(a)\go_i$
and $da=\sum a_1\go(a_2)$; recall also $P_L(da)=\go(a)=\sum S(a_1)da_2$
and $X_i\cdot a=\sum a_1X_i(a_2)$.  In bicovariant situations the right
invariant forms satisfy $\eta(a)=\sum da_1\cdot S(a_2)$ and $da=\sum
\eta(a_1)a_2$.
\\[3mm]\indent
{\bf REMARK 8.4.}
It will be worthwhile to list still more structural and computational
features of ${\mf O}(SL_Q(2))$ and $\breve{{\mf U}}_q({\mf s}\ell_2)$
before going to calculations for qKdV type equations.  Thus following
\cite{k1} we return to $\breve{{\mf U}}_q({\mf s}\ell_2)$ (cf. also 
\eqref{19z} - \eqref{21z}).  Note first for a dual pairing $<\,,\,>$ of
bialgebras U and A via $<\,,\,>:\,U\times A\to {\bf C}$ one requires
\bq\label{15x}
<\gD_Uf,a_1\ot a_2>=<f,a_1a_2>;\,\,<f_1f_2,a>=<f_1\ot f_2,\gD_Aa>;
\end{equation}
$$<f,1_A>=\gep_U)f);\,\,<1_U,a>=\gep_A(a)$$
Now in order to calculate $X_0(b)=\cdots$ after \eqref{1x} we need the
operators on ${\mf O}(M_q(2))$ (extended to ${\mf O}(SL_q(2))$) given via
\bq\label{16x}
\gD(a)=a\ot a+b\ot c;\,\,\gD(b)=a\ot b+b\ot d;\,\,\gD(c)=c\ot a+d\ot c;\,\,
\gD(d)=c\ot b+d\ot d
\end{equation}
with ${\bf (D35)}\,\,\gep(a)=\gep(d)=1$ and $\gep(b)=\gep(c)=0$.  This is
equivalent to ${\bf (D36)}\,\,\gD(u^i_j)=\sum u_k^i\ot u^k_j$ and
$\gep(u^i_j)=\gd_{ij}$.  One can then calculate e.g. ($<\,,\,>\,\sim\,<\,,\,>
\breve{{}}\,$)
\bq\label{17x}
X_0(b)=q^{-1/2}FK(b)=q^{-1/2}(<F\ot K,a\ot b+b\ot d>)=1
\end{equation}
In general for q not a root of unity the dual pairing is given via
\bq\label{18x}
<K^mE^nF^{\ell},d^sc^rb^t>=q^{(n-2)^2}\left[\begin{array}{c}
s\\
n-r
\end{array}\right]_{q^2}\gag^{srt}_{mn\ell}
\end{equation}
for $0\leq n-r=\ell-t\leq s$ (otherwise $0$) and
\bq\label{19x}
<K^mE^nF^{\ell},a^sc^rb^t>=\gd_{rn}\gd_{t\ell}\gag^{-srt}_{mn\ell}
\end{equation}
where in all cases
\bq\label{20x}
\gag^{srt}_{mn\ell}=\frac{q^{m(s+r-t)/2}q^{-s(n+\ell)/2}(q^2;q^2)_{\ell}
(q^2;q^2)_n}{q^{n(n-1)/2}q^{\ell(\ell-1)/2}(1-q^2)^{\ell+n}}
\end{equation}
(cf. \cite{kx} for proof).  Thus e.g. for $X_1(a)$ the term in $K^4$ involves
$n=\ell=0$ and $m=4$ with $s=1$ and $r=t=0$ so $n-r=0=\ell-t$ and, since
$\gep(a)=1$, we obtain
\bq\label{21x}
\gag_{400}^{-100}=q^{-2}=<K^4,a>;\,\,
X_1(a)=(1-q^{-2})^{-1}(1-q^{-2})=1
\end{equation}
as desired.  For the 3-D calculations we shift notation $a,b,c,d\to
A,B,C,D$ and write (cf. \eqref{2x})
\bq\label{22x}
\go_0=\go(B));\,\,\go_2=\go(C);\,\,\go_1=\go(A)=-q^{-2}\go(D);\,\,
dA=B\go_2+A\go_1;
\end{equation}
$$dB=A\go_0-q^2B\go_1;\,\,dC=C\go_1+D\go_2;\,\,dD=-q^2
D\go_1+C\go_0$$
One can show (since $\gD(X_i)=\gep\ot X_i+K_k\ot f^k_i$ with $\gD
f^r_s=f^r_k\ot f^k_s$) that $f_j^j=K^2$ for $j=0,2$ and $f^1_1=K^4$ with
$f^r_s=0$ otherwise.  This leads to commutation relations (cf. \eqref{3x})
\bq\label{23x}
\go_jA=q^{-1}A\go_j;\,\,\go_jB=qB\go_j;\,\,\go_jC=q^{-1}C\go_j;\,\,
\go_jD=qD\go_j\,\,(j=0,2);
\end{equation}
$$\go_1A=q^{-2}A\go_1;\,\,\go_1B=q^2\go_1;\,\,\go_1C=q^{-2}C\go_1;
\,\,\go_1D=q^2D\go_1$$
Further one has (cf. \eqref{4x})
\bq\label{24x}
q^2X_1X_0-q^{-2}X_0X_1=(1+q^2)X_0;\,\,q^2X_2X_1-q^{-2}X_1X_2=(1+q^2)X_2;
\end{equation}
$$qX_2X_0-q^{-1}X_0X_2=-q^{-1}X_1$$
These will correspond to Maurer-Cartan equations (see below).
$\hfill\bs$
\\[3mm]\indent
Now look at the Maurer-Cartan (MC) equations for the 3-D calculus in
the form
\eqref{5x} or \eqref{6x}, namely
\bq\label{25x}
d\go(a)=-\sum\go(a_1)\ot\go(a_2);\,\,\sum X_i(a)d\go_i=-\sum
X_iX_j(a)\go_i\wg\go_j
\end{equation}
We know formulas for $\go_i$ via \eqref{22x}-\eqref{23x} and $X_i$ via
\eqref{24x} and one would like to put variables $x,t$ in the $\go_i$ or the $X_i$ to
give equations of the form \eqref{2x}-\eqref{4x}, etc.  Note also \eqref{9z} for
other relations between $A,B,C,D$.
Now from \eqref{22x} and \eqref{9z} we have
$dA=A\go_1+B\go_2$ and $dC=C\go_1+D\go_2$.  Then, recalling e.g. $BD=qDB$,
there results $qD(dA)=qDA\go_1+qDB\go_2$ while $B(dC)=BC\go_1+BD\go_2$ which
implies
\bq\label{26x}
qDdA-qDA\go_1=BdC-BC\go_1\Rightarrow \go_1=\frac{BdC-qDdA}{BC-qDA}=
DdA-q^{-1}BdC
\end{equation}
since $AD-qBC=DA-q^{-1}BC=1\Rightarrow (qDA-BC)q^{-1}=1$.  Similarly
$AC=qCA$ leads to $AdC=AC\go_1+AD\go_2$ and $qCdA=qCA\go_1+qCB\go_2$ which
implies
\bq\label{27x}
AdC-AD\go_2=qCdA-qCB\go_2\Rightarrow\go_2=\frac{AdC-qCdA}{qCB-AD}=qCdA-AdC
\end{equation}
Next from $dB=A\go_0-q^2B\go_1$ and $dD=-q^2D\go_1+C\go_0$ we get via
$BD=qDB$, the formulas $qDdB=qDA\go_0-q^3DB\go_1$ and
$BdD=-q^2BD\go_1+BC\go_0$ which implies $qDdB-qDA\go_0=BdD-BC\go_0$ or
\bq\label{28x}
\go_0=\frac{qDdB-BdD}{qDA-BC}=DdB-q^{-1}BdD
\end{equation}
\begin{proposition}
The MC equations can be expressed in terms of $A,B,C,D,dA,dB,dC,$ and $dD$
via \eqref{26x} - \eqref{28x} in the form
\bq\label{29x}
\go_1=DdA-q^{-1}BdC;\,\,\go_2=qCdA-AdC;\,\,\go_0=DdB-q^{-1}BdD
\end{equation}
\end{proposition}
\indent
{\bf REMARK 8.5.}
We sketch a few ideas regarding $P_L$ from \cite{k1,sp,sq,wq}.
The presentation in \cite{wq} is definitive and we will a more or less
complete treatment here, proving a few results which are only sketched in
\cite{k1} for example.  Recall first the standard formulas $\gep(a)=\sum
S(a_1)a_2=\sum a_1S(a_2),\,\,a=\sum a_1\gep(a_2)=\sum \gep(a_1)a_2,$ etc.
Now one shows that there exists a unique projection $P_L:\,\gG\to{}_{inv}\gG$
such that ${\bf (D37)}\,\,P_L(a\rho)=\gep(a)P_L(\rho)$ with $\rho=\sum
a_kP_L(\rho_k)$ for $\gD_L(\rho)=\sum a_k\ot \rho_k$.  To see this write for
$\rho\in\gG,\,\,P_L(\rho)=\sum S(a_k)\rho_k$.  Then $P_L$ is clearly linear
and to show $P_L(\rho)$ is left invariant consider ${\bf (D38)}\,\,
\gD(a_k)=\sum a_{k\ell}\ot b_{k\ell}$ with $\gD_L(\rho_k)=\sum
c_{km}\ot\rho_{km}$ (a more strict Sweedler notation is adopted below).
For left covariant bimodules the diagram 
\[\begin{CD}
\\\gG            @>{\gD_L}>>         A\ot\gG\\
@V{\gD_L}VV                                      @VV{\gD\ot id}V\\
A\ot\gG                          @>>{id\ot\gD_L}>             A\ot A\ot \gG
\end{CD}\]
is commutative so ${\bf (D38)}\,\,\sum a_{k\ell}\ot b_{k\ell}\ot\rho_k=\sum
a_k\ot c_{km}\ot \rho_{km}$ and
\bq\label{30x}
\gD_LP_L(\rho)=\sum\gD_LS(a_k)\rho_k=\sum\gD
S(a_k)\gD_L\rho_k=\sum\gD(S(a_k))(c_{km}\ot\rho_{km})
\end{equation}
Now (recall ${\bf (D30)}:\,\gD\ci S=\tau\ci(S\ot S)\ci\gD$)
\bq\label{31x}
\cdot(\gD S\ot1\ot 1)(a_k\ot c_{km}\ot\rho_{km})=\cdot(\gD S\ot 1\ot 1)
(a_{k\ell}\ot b_{k\ell}\ot\rho_k)=
\end{equation}
$$=\gD S(a_{k\ell})(b_{k\ell}\ot \rho_k)=\sum(S(b_{k\ell m})\ot S(a_{k\ell
m}))(b_{k\ell}\ot\rho_k)=\sum S(b_{k\ell m})b_{k\ell}\ot S(a_{k\ell
m})\rho_k$$
This is better expressed in a more strict Sweedler notation as
(writing $\gD_L(\rho)=\sum
\rho_{(-1)}\ot\rho_0,\,\,\gD(\rho_{(-1)})=\rho_{(-2)}\ot\rho_{(-3)},\,\,
\gD_L\rho_0=\sum \rho_{(-4)}\ot \rho_1,$ and $\gD(\rho_{(-2)})=\rho_{(-5)}\ot
\rho_{(-6)}$)
\bq\label{40x}
\cdot(\gD S\ot 1\ot 1)(\rho_{-1)}\ot\rho_{(-4)}\ot\rho_1)=\cdot(\gD S\ot 1
\ot 1)(\rho_{(-2)}\ot\rho_{(-3)}\ot\rho_0)=
\end{equation}
$$=\gD S(\rho_{(-1)})(\rho_{(-4)}\ot \rho_1)=\gD
S(\rho_{(-2)})(\rho_{(-3)}\ot\rho_0)=\sum S(\rho_{(-6)})\rho_{(-3)}
\ot S(\rho_{(-5)})\rho_0=$$
$$=\sum S(\rho_{(-7)})1\ot S(\rho_{(-5)})\rho_0=1\ot\sum S(\rho_{(-1)})\rho_0
=1\ot P_L(\rho)$$
The last step involves
$\gep(\rho_{(-7)})S)\rho_{(-5)})=S(\gep(\rho_{(-7)})\rho_{(-5)})=
S(\rho_{(-6)})$.  The proof of $P_L(a\rho)=\sum \gep(a)P_L(\rho)$ is
straightforward (cf. \cite{k1,wq}) and one notes also that for $\rho=
\sum \rho_{(-1)}P_L(\rho_0)$ we have $P_L(\rho)=\sum \gep(\rho_{(-1)})
P_L(\rho_0)=P_L(\sum \gep(\rho_{(-1)})\rho_0)=P_L(\rho)$.  Note 
{\bf (D37)} says that any $\rho\in\gG$ has the form $\rho=\sum
\rho_{(-1)}P_L(\rho_0)$ so ${}_{inv}\gG$ serves as a basis for $\gG$ (i.e.
$\gG=A\cdot{}_{inv}\gG$).
Next from $\gD_L(adb)=\gD(a)(id\ot d)\gD(b)$ one gets $\gD_L(da)=\sum a_i\ot
da_2$ so $\rho=\sum a_1P_L(da_2)$ and $\go(a)=P_L(da)=\sum
S(\rho_{(-1)})\rho_0=\sum S(a_1)da_2$.  One has then ${\bf (D39)}\,\,
\go(a)=\sum S(a_1)da_2$ and $da=\sum a_1\go(a_2)=\sum a_1P_L(da_2)$. 
Finally in deriving the MC equations one writes $\gep(a)=\sum S(a_1)a_2$ to
obtain
\bq\label{41x}
0=\sum dS(a_1)a_2+\sum S(a_1)da_2\Rightarrow \sum S(a_1)da_2S(a_3)=-\sum
dS(a_1)a_2S(a_3)=
\end{equation}
$$=-\sum dS(a_1)\gep(a_4)=-\sum dS(a)$$
(since $\sum a_1\gep(a_4)\sim a$).$\hfill\bs$
\\[3mm]\indent
Now given \eqref{41x} one computes the MC equations \eqref{5x} as in
\cite{k1} via
\bq\label{42x}
d\go(a)=\sum d(S(a_1)da_2)=\sum dS(a_1)\wg da_2=-\sum S(a_1)da_2S(a_3)\wg
da_4=
\end{equation}
$$=-\sum S(a_1)da_2\wg S(a_3)da_4=-\sum \go(a_1)\wg\go(a_2)$$
For the 3-D calculus we have now \eqref{9z}, \eqref{19z}-\eqref{21z},
\eqref{2x}-\eqref{5x}, \eqref{8x}, and \eqref{29x} which we collect in
\begin{proposition}
Properties of the 3-D calculus include:
\begin{enumerate}
\item
$AB=qBA,\,\,AC=qCA,\,\,BD=qDB,\,\,CD=qDC,\,\,BC=CB,\,\,AD-DA=(q-q^{-1})BC$
\item
$X_0=q^{-1/2}FK,\,\,X_2=q^{1/2}EK,\,\,X_1=(1-q^{-2})^{-1}(\gep-K^4)$
\item
$\gD(K)=K\ot K,\,\,\gD(E)=E\ot K+K^{-1}\ot E,\,\,\gD(F)=F\ot K+K^{-1}\ot
F,\,\,S(K)=K^{-1},\,\,S(E)=-qE,\,\,S(F)=-q^{-1}F,\,\,\gep(K)=1,\,\,\gep(E)=
\gep(F)=0$
\item
$\gD(X_j)=\gep\ot X_j+X_j\ot K^2\,\,(j=0,2),\,\,\gD(X_1)=\gep\ot X_1+X_1\ot
K^4$
\item
$\go_0=\go(B),\,\,\go_2=\go(C),\,\,\go_1=\go(A)-q^{-2}\go(D),
\,\,dA=B\go_2+A\go_1,\,\,dB=A\go_0-q^2
B\go_1,\,\,dC=C\go_1+D\go_2,\,\,dD=-q^2D\go_1+C\go_0$
\item
$\go_jA=q^{-1}A\go_j,\,\,\go_jB=qB\go_j,\,\,
\go_jC=q^{-1}C\go_j\,\,(j=0,2),\,\,\go_1A=q^{-2}A\go_1,
\,\,\go_1B=q^2B\go_2,\,\,\go_1C=q^{-2}C\go_1,\,\,\go_1 D=q^2
D\go_1$
\item
$q^2X_1X_0-q^{-2}X_0X_1=(1+q^2)X_0,\,\,q^2X_2X_1-q^{-2}X_1X_2=(1+q^2)X_2,
\,\,qX_2X_0-q^{-1}X_0X_2=
-q^{-1}X_1$
\item
$\go_j\wg\go_j=0\,\,(j=0,1,2),\,\,\go_2\wg\go_0=-q^2\go_0\wg\go_2,
\,\,\go_1\wg\go_0=-q^4\go_0\wg\go_1,\,\,\go_2\wg\go_1=-q^4\go_1\wg\go_2$
\item
$\go_1=DdA-q^{-1}BdC,\,\,\go_2=qCdA-AdC,\,\,\go_0=DdB-q^{-1}BdD$
\end{enumerate}
\end{proposition}
\indent
The equations for KdV in \eqref{1x}-\eqref{8x} involve
$A_{\mu}=g^{-1}\pp_{\mu}g=i\sum A_{\mu}^a(x,t)T_a$ where $T_a\sim X_a\sim$ 
Lie algebra elements, and one puts the $x,t$ in $g(x,t)=exp(i\sum
\gt^a(x,t)T_a)$.  We don't really have an exponential map now so we should
phrase the MC equations in terms of the $X_i$.  Given (8) above the MC
equations \eqref{6x} take the form
\bq\label{43x}
\sum
X_i(a)d\go_i=-(X_0X_1-q^4X_1X_0)(a)\go_0\wg\go_1-(X_0X_2-q^2X_2X_0)(a)
\go_0\wg\go_2-
\end{equation}
$$-(X_1X_2-q^4X_2X_1)(a)\go_1\wg\go_2$$
On the other hand from \cite{k1}, p. 510 one gets, using (5), (8), and \eqref
{42x},
\bq\label{44x}
d\go_0=(q^2+q^4)\go_0\wg\go_1;\,\,d\go_1=-\go_0\wg\go_2;
\,\,d\go_2=(q^2+q^4)\go_1\wg\go_2
\end{equation}
This should all be compatiable with $\eqref{5x}\sim (7)$ and we put \eqref
{44x} into \eqref{43x} to get
\bq\label{45x}
X_0(a)(q^2+1)\go_0\wg\go_1+X_1(a)(-\go_0\wg\go_2)+X_2(a)q^2(q^2+1)\go_1\wg
\go_2=
\end{equation}
$$=(q^4X_1X_0-X_0X_1)(a)\go_0\wg\go_1+(q^2(X_2X_0-X_0X_2)(a)\go_0\wg\go_2
+(q^4X_2X_1-X_1X_2)(a)\go_1\wg\go_2$$
which becomes e.g.
\bq\label{46x}
X_0(a)q^2(q^2+1)=(q^4X_1X_0-X_0X_1)(a)\equiv (1+q^2)X_0(a)=
(q^2X_1X_0-q^{-2}X_0X_1)(a)
\end{equation}
Thus (7) corresponds to the MC equations and now one could envision elements
$\ga(x,t)=\sum \gt_a(x,t)X_a$ in working toward a qKdV equation.  We will
think of $\pp_x$ as some sort of q-derivative and $\pp_t\sim \pp/\pp t$.  
First let us recall that $A'$ is a dual Hopf algebra via
\bq\label{52x}
fg(a)=(f\ot g)\gD(a)=f(a_1)\ot f(a_2);\,\,\gD f(a\ot b)=(f\ci m)(a\ot
b)=f(ab);
\end{equation}
$$Sf(a)=f(S(a));\,\,\gep_{A'}(f)=f(1);\,\,1_{A'}(a)=\gep(a)$$
and as before $(X,adb)=\gep(a)X(b)\in{\bf C}$ expressses the duality
of $T_{\gG}$ and $\gG$.  One should be able to use (7) as rules of
multiplication in $A'$ so consider elements
\bq\label{53x}
{\mf B}_x=\ga^{-1}\pp_x\ga;\,\,{\mf B}_t=\ga^{-1}\pp_t\ga
\end{equation}
leading to ($\pp_x\ga^{-1}=-\ga^{-1}\pp_x\ga\ga^{-1}$)
\bq\label{54x}
{\mf B}_x{\mf B}_t=\ga^{-1}\pp_x\ga\ga^{-1}\pp_t\ga;\,\,{\mf B}_t{\mf
B}_x=\ga^{-1}\pp_t\ga\ga^{-1}\pp_x\ga
\end{equation}
\bq\label{55x}
\pp_t{\mf B}_x=\pp_t(\ga^{-1}\pp_x\ga)=-\ga^{-1}\pp_t\ga\ga^{-1}\pp_x\ga+
\ga^{-1}\pp_t\pp_x\ga;
\end{equation}
$$\pp_x{\mf B}_t=-\ga^{-1}\pp_x\ga\ga^{-1}\pp_t\ga+\ga^{-1}\pp_t\pp_x\ga$$
Assume now $\pp_t\pp_x=\pp_x\pp_t$ in which case
\bq\label{56x}
\pp_t{\mf B}_x-\pp_x{\mf B}_t={\mf B}_x{\mf B}_t-{\mf B}_t{\mf B}_x
\end{equation}
which corresponds to {\bf (D2)}.  Now ${\mf B}_x=\ga^{-1}\sum
\pp_x\gt^aX_a$, etc. so one can write ${\mf B}_x{\mf B}_t=\ga^{-1}\sum
\pp_x\gt^aX_a\ga^{-1}\sum\pp_t\gt^bX_b$.  Computation of $\ga^{-1}$ or $X_a\ga^{-1}$
seems excessive here since \eqref{56x} is tautological and we simply state in
summary
\begin{theorem}
Maurer-Cartan equations in $T_{\gG}$ correspond simply to (7) of Proposition 8.2 and
in principle
\eqref{56x}, corresponding to {\bf (D2)} should lead to qKdV type equations.
\end{theorem}
\indent
This should work but we can also devise a simpler procedure directly from
\cite{cz} (which was the source from which  \cite{dz} extracted the more general
\eqref{1x}-\eqref{8x}).  Thus look at
$SL(2,{\bf R})$ with matrices
\bq\label{58x}
X=\left(\begin{array}{cc}
a & b\\
c & d
\end{array}\right);\,\,ad-bc=1
\end{equation}
The right invariant MC form is
\bq\label{59x}
\go=dXX^{-1}=\left(\begin{array}{cc}
\go^1_1 & \go^2_1\\
\go_2^1 & \go_2^2
\end{array}\right)
\end{equation}
where $\go_1^1+\go^2_2=0$.  The structure equation of $SL(2,{\bf R})$ or 
MC equation is ${\bf (D40)}\,\,d\go=\go\wg\go$ or explicitly
\bq\label{60x}
d\go^1_1=\go_1^2\wg \go_2^1;\,\,d\go_1^2=2\go^1_1\wg\go_1^2;\,\,d\go_2^1
=2\go_2^1\wg\go_1^1
\end{equation}
Now let U be a neighborhood in the $(x,t)$ plane and consider a smooth map
$f:\,U\to SL(2,{\bf R})$.  The pullback of the MC form can be written as
\bq\label{61x}
\go_1^1\sim \eta dx+Adt;\,\,\go_1^2\sim Qdx + Bdt;\,\,\go_2^1\sim rdx
+Cdt
\end{equation}
with coefficient functions of $x,t$.  The equations \eqref{60x} become
\begin{enumerate}
\item
$-\eta_t+A_x-QC+rB=0$
\item
$-Q_t+B_x-2\eta B+2QA=0$
\item
$-r_t+C_x-2rA+2\eta C=0$
\end{enumerate}
Take $r=1$ with $\eta$ independent of $(x,t)$ and set $Q=u(x,t)$.  Then from 
(1) and (3) one gets
\bq\label{63x}
A=\eta C+\frac{1}{2}C_x;\,\,B=uC-\eta C_x-\frac{1}{2}C_{xx}
\end{equation}
Putting this in the (2) above yields $u_t=K(u)$ where
\bq\label{64x}
K(u)=u_xC+2uC_x+2\eta^2C_x-\frac{1}{2}C_{xxx}
\end{equation}
In the special case $C=\eta^2-(1/2)u$ one gets the KdV equation
\bq\label{65x}
u_t=\frac{1}{4}u_{xxx}-\frac{3}{2}uu_x
\end{equation}
\indent
Now a map $f:\,U\to SL_q(2)$ means simply expressing $A,B,C,D$ as functions
of $(x,t)$ and the pullback of forms from $\gG$ built over ${\mf O}(SL_q(2))$
should be modelable on the standard procedure from differential geometry.
Thus for a manifold map $f:\,M\to N$ one has maps $f_*:\,TM\to TN$ and
$f^*:\,T^*N\to T^*M$ defined via $f_*v(g)=v(g\ci f)\,\,(v\in T_p(M))$ and
$f^*(dg)=d(g\ci f)\,\,(dg\in T^*_{f(p)}(N))$.
The analogue here would work from formal power series $g$ in terms of
$A^sC^rB^n$ or $D^sC^rB^n$ and formal power series $f(x,t)$ in
$x^{\ga}t^{\gb}$ (i.e. $f\in {\bf C}[[x,t]]$).  The $X_a\in T_{\gG}$
correspond to elements in $TN$ and $\gG\sim T^*N$ so in some sense e.g.
${\bf (D41)}\,\,f^*(da)=d(a(x,t))$ etc.  As an appropriate FODC for
$f^*(\gG)=\tl{\gG}$ it may be
necessary to have a noncommutative formulation for $\tl{\gG}$ and we proceed in a
somewhat ad hoc manner.  Thus one could use \eqref{29x} to express the $\go_i$ in
terms of $dA,dB,dC,dD$ and \eqref{44x} for the MC equations.  First however let us
simply write d'apr\`es \eqref{61x}
\bq\label{66x}
\go_1^1\sim\go_1=\eta dx+\mu dt;\,\,\go_2^1\sim w_0=\ga dx+\nu dt;\,\,
\go_1^2\sim \go_2=\gag dx+\gb dt
\end{equation}
(so $\eta\sim\eta,\,\,A\sim\mu,\,\,r\sim\ga,\,\,B\sim\gb,\,\,u\sim\gag,\,\,
C\sim\nu$).  Then 
\bq\label{67x}
d\go_0=(q^2+q^4)\go_0\wg\go_1;\,\,d\go_1=-\go_0\wg\go_2;\,\,d\go_2=
(q^2+q^4)\go_1\wg\go_2
\end{equation}
Pick again $\eta$ constant, $\gag=u(x,t)$, and $\ga=1$ so
$(\spadesuit)\,\,
\go_1=\eta dx+\mu dt;\,\,\go_0=dx+\nu dt;\,\,\go_2=udx+\gb dt$
Assume first $dh=\pp_xhdx+h_tdt$ for $\pp_x$ possibly a q-derivative.  Then
\eqref{67x} for $d\go_1$ implies $(\bullet)\,\,
\pp_x\mu-\nu u+\gb=0$
which corresponds to (1) above
and for $d\go_0$ we get
$\pp_x\nu=(q^2+q^4)(\mu-\nu\eta)$
analogous to (3) above.  Finally from $d\go_2$ there arises $(\bl)\,\,
\pp_x\gb-u_t=(q^2+q^4)(\eta\gb-\mu u)$
analogous to (2) above.  Summarizing (with $q^2+q^4={\mf Q}$)
\begin{enumerate}
\item
$\pp_x\mu -\nu u+\gb=0$
\item
$-u_t+\pp_x\gb={\mf Q}(\eta\gb-\mu u)$
\item
$\pp_x\nu={\mf Q}(\mu-\nu\eta)$
\end{enumerate}
To eliminate as in \eqref{63x}-\eqref{65x} one has from (1) and (3)
$(\clubsuit)\,\,
\pp_x\nu={\mf Q}(\mu-\nu\eta);\,\,\pp_x\mu=u\nu-\gb$
Hence ${\bf (D42)}\,\,\gb=\nu u-\pp_x(\nu\eta+(1/{\mf Q})\pp_x\nu)=\nu
u-\eta\pp_x\nu-(1/{\mf Q})\pp_x^2\nu$ as in \eqref{63x}.  Now put this in
(2) to get
\bq\label{73x}
u_t=\pp_x\gb-{\mf Q}(\eta\gb-\mu u)=(\pp_x\nu)u+\nu\pp_xu-\eta\pp_x^2\nu
-(1/{\mf Q})\pp_x^3\nu-
\end{equation}
$$-{\mf Q}\eta(\nu u-\eta\pp_x\nu-(1/{\mf Q})\pp_x^2\nu)+{\mf Q}
u((1/{\mf Q})\pp_x\nu+\nu\eta)$$
and modeled on \eqref{63x}-\eqref{65x} one would try $\nu=\eta^2-Pu$
and $P=1/{\mf Q}$ to get
\bq\label{75x}
u_t=\frac{1}{{\mf Q}^2}\pp_x^3u-\frac{3}{2}u\pp_xu
\end{equation}
which actually can be rescaled to become a standard KdV equation with variables
depending on q.  This assumes $\pp_x$ is a normal derivative.  If e.g. $\pp_x\sim
\pp_q^x$ then we need only check first $\pp_t\pp_q^x=\pp_q^x\pp_t$ which is
immediate and note that one is assuming $dxdt=-dtdx$.  However for $\pp_x\sim
\pp_q^x$ the term $\pp_q^x(\nu u)$ becomes ${\bf (D43)}\,\,
\pp_q^x(\nu u)=\pp_q^x\nu D_xu+\nu\pp_q^xu$ where $D_xu=u(qx,t)$.  
This changes some terms in \eqref{73x} et suite and leads to ($\pp_q\sim\pp_q^x$)
\bq\label{79x}
u_t=\frac{1}{{\mf Q}^2}\pp_q^3u-\frac{2}{{\mf Q}}u\pp_qu-D_xu\pp_qu
\end{equation}
\begin{theorem}
Assume we can write \eqref{66x}, with $dxdt=-dtdx$, and $\pp_x=\pp_q^x$;
then \eqref{79x} is a kind of qKdV equation arising naturally from 
MC equations.
\end{theorem}

\section{MAURER-CARTAN AND THE QUANTUM LINE}
\renewcommand{\theequation}{9.\arabic{equation}}
\setcounter{equation}{0}

We go now to a quantum line ${\bf R}_q^1$ as in Section 6 coupled with a
time variable (e.g. $A_q=C({\bf R})\ot {\bf
R}_q^1$).  Recall one possibility based on {\bf (B58)}-{\bf (B61)} involves
\bq\label{1y}
x\gL=q\gL x;\,\,xdx=qdxx;\,\,dx\gL=q\gL dx;\,\,xd\gL=qd\gL x;\,\,e_1x=q\gL
x;\,\,
\end{equation}
$$e_1\gL=0;\,\,df(e_1)=e_1f;\,\,e_2\gL=q\gL x;\,\, e_2x=0;\,\,df(e_2)=e_2f$$
This was rejected before since $\gL$ is in a sense an element of the phase
space associated to $x$ and one was interested in position space geometry.
For our purposes however we could consider this situation, treating
$\gL$ as a parameter to give structure to the picture.
Another formulation can be based on the real DC via {\bf (B84)}-{\bf (B94)} 
but we omit this here.  Thus we mimic \eqref{66x} et suite in taking 
\bq\label{2y}
\go_1^1\sim \go_1=\eta dx+\mu dt+\phi d\gL;\,\,\go_2^1\sim \go_0=\ga dx+\nu
dt+\psi d\gL;
\end{equation}
$$\go_1^2\sim \go_2=\gag dx+\gb dt+\chi d\gL$$
The MC equations \eqref{44x} are
\bq\label{3y}
d\go_0={\mf
Q}(\go_0\wg\go_1);\,\,d\go_1=-\go_0\wg\go_2;\,\,d\go_2={\mf Q}(\go_1\wg\go_2)
\end{equation}
and we must find expressions for $df$ and $d\go$.  First (recalling e.g.
Example 7.4) consider from \eqref{1y}
\bq\label{4y}
dx^2=xdx+dxx=xdx+q^{-1}xdx=(1+q^{-1})xdx;\,\,\cdots,\,\,dx^n=[n]_{q^{-1}}
x^{n-1}dx
\end{equation}
Since $\gL d\gL$ does not arise in \eqref{1y} we go back to Example 8.1 of
the generalized q-plane which has relations
\bq\label{5y}
xdx=qdxx;\,\,dx\gL=q\gL dx;\,\,xd\gL=qd\gL x;
\end{equation}
$$d\gL\gL=q\gL
d\gL;\,\,dx^2=d\gL^2=0;\,\,dxd\gL+qd\gL dx=0$$
We add the additional hypothesis ${\bf (E1)}\,\,d\gL\gL=q\gL d\gL$ and then
from {\bf (C74)}
\bq\label{6y}
d\gL^m=[m]_q\gL^{m-1};\,\,df=D_{\gL}\pp_{q^{-1}}^xfdx+\pp_q^{\gL}fd\gL
\end{equation}
We are after some form of qKdV now and not qKP - the variable $\gL$ is purely an
artifice to give ``quantum" meaning to ${\bf R}$.  Note here ${\bf
(E2)}\,\,dx\gL^m=q^m\gL^mdx$ plays a role in \eqref{6y}.  Now from
\eqref{3y} we write
\bq\label{7y}
d\go_0=(D_{\gL}\pp_{q^{-1}}^x\ga dx+\pp_q^{\gL}\ga
d\gL+\ga_tdt)dx+(D_{\gL}\pp_{q^{-1}}^x\nu dx+\pp_q^{\gL}\nu d\gL+\nu_tdt)dt+
\end{equation}
$$+(D_{\gL}\pp_{q^{-1}}\psi dx+\pp_q^{\gL}\psi d\gL+\psi_tdt)d\gL=
{\mf Q}(\go_0\wg\go_1)={\mf Q}(\ga dx+\nu dt+\psi d\gL)\wg (\eta dx+\mu
dt+\phi d\gL)$$
This leads to
\bq\label{8y}
\pp_q^{\gL}\ga d\gL dx+\ga_tdtdx+D_{\gL}\pp_{q^{-1}}^x\nu dxdt+
\pp_q^{\gL}\nu d\gL dt+D_{\gL}\pp_{q^{-1}}^x\psi dx d\gL+\psi_tdtd\gL=
\end{equation}
$$={\mf Q}\{\ga dx\eta dx+\ga dx\mu dt+\ga dx\phi d\gL+\nu dt\eta dx+\nu
dt\mu dt+\nu dt\phi d\gL+\psi d\gL\eta dx+\psi d\gL\mu dt+\psi d\gL\phi
d\gL\}$$
Now e.g. $dx\eta=dx\sum \eta_{nmk}x^n\gL^mt^k=\sum \eta_{nmk}q^{-n}x^n
q^mx^n\gL^mt^kdx=D_x^{-1}D_{\gL}\eta dx$ (since $dxx^n=q^{-n}x^ndx$ and {\bf
(E2)} holds).  Similarly $d\gL x=q^{-1}xd\gL\Rightarrow d\gL
x^n=q^{-n}x^nd\gL$ and $d\gL\gL=q\gL d\gL\Rightarrow d\gL\gL^m=q^m\gL^md\gL$
so $d\gL\eta=D_x^{-1}D_{\gL}\eta d\gL$ and one has (assuming $dxdt+dtdx=0$
and $d\gL dt+dtd\gL=0$ and omitting some calculations)
\bq\label{10y}
D_{\gL}\pp_{q^{-1}}^x\psi-q^{-1}\pp_q^{\gL}\ga={\mf Q}(\ga D_x^{-1}D_{\gL}
\phi-q^{-1}\phi D_x^{-1}D_{\gL}\eta);
\end{equation}
$$D_{\gL}\pp_{q^{-1}}^x\nu-\ga_t={\mf Q}[\ga D_x^{-1}D_{\gL}\mu-\nu\eta];\,\,
\pp_q^{\gL}\nu-\psi_t={\mf Q}[\psi D_x^{-1}D_{\gL}\mu-\nu\phi]$$
with the $dxdt$ term in the middle.  Similarly computing for $d\go_1$
leads to
\bq\label{14y}
D_{\gL}\pp^x_{q^{-1}}\phi-q^{-1}\pp_q^{\gL}\eta=q^{-1}D_x^{-1}D_{\gL}\chi-
\psi D_x^{-1}D_{\gL}\psi;
\end{equation}
$$D_{\gL}\pp^x_{q^{-1}}\mu-\eta_t=\nu\gag-\ga D_x^{-1}D_{\gL}\gb;\,\,
\pp_q^{\gL}\mu-\phi_t=\nu\chi-\psi D_x^{-1}D_{\gL}\gb$$
Finally from $d\go_2$
\bq\label{18y}
D_{\gL}\pp^x_{q^{-1}}\chi-q^{-1}\pp_q^{\gL}\gag={\mf Q}(\eta
D_x^{-1}D_{\gL}\chi-q^{-1}\phi D_x^{-1}D_{\gL}\gag);
\end{equation}
$$D_{\gL}\pp^x_{q^{-1}}\gb-\gag_t={\mf Q}(\eta D_x^{-1}D_{\gL}\gb-\mu\gag);\,\,
\pp_q^{\gL}\gb-\chi_t={\mf Q}(\phi D_x^{-1}D_{\gL}\gb-\chi\mu)$$
We look first at the $dxdt$ equations in \eqref{10y}, \eqref{14y}, and
\eqref{18y} to get
\bq\label{19y}
D_{\gL}\pp^x_{q^{-1}}\nu-\ga_t={\mf Q}[\ga D_x^{-1}D_{\gL}\mu-\nu\eta];
\end{equation}
$$D_{\gL}\pp_{q^{-1}}^x\mu-\eta_t=\nu\gag-\ga D_x^{-1}D_{\gL}\gb;\,\,
D_{\gL}\pp^x_{q^{-1}}\gb-\gag_t={\mf Q}(\eta D_x^{-1}D_{\gL}\gb-\mu\gag)$$
Compare this with the enumeration in Section 8, namely
\bq\label{20y}
\pp_x\nu={\mf Q}(\mu-\nu\eta);\,\,-u_t+\pp_x\gb={\mf Q}(\eta\gb-\mu u);\,\,
\pp_x\mu-\nu\mu+\gb=0
\end{equation}
where $u\sim\gag$ and $r\sim \ga$.  If we take $\eta=$ constant and $\ga=1$
again in \eqref{19y} there results
\bq\label{21y}
D_{\gL}\pp^x_{q^{-1}}\nu={\mf Q}(D_x^{-1}D_{\gL}\mu-\eta\nu);
\end{equation}
$$D_{\gL}\pp^x_{q^{-1}}\gb-u_t={\mf Q}(\eta D_x^{-1}D_{\gL}\gb-\mu u);\,\,
D_{\gL}\pp^x_{q^{-1}}\mu=\nu u-D_x^{-1}D_{\gL}\gb$$
This is quite parallel, modulo shifts $D_x$ and $D_{\gL}$.  Thus
\begin{enumerate}
\item
$\pp_x\mu-\nu u+\gb=0\sim D_{\gL}\pp^x_{q^{-1}}\mu=\nu u-D_x^{-1}D_{\gL}\gb$
\item
$\pp_x\gb-u_t={\mf Q}(\eta\gb -\mu u)\sim D_{\gL}\pp^x_{q^{-1}}\gb-u_t={\mf
Q}(\eta D_x^{-1}D_{\gL}\gb-\mu u);$
\item
$\pp_x\nu={\mf Q}(\mu-\nu\eta)\sim D_{\gL}\pp^x_{q^{-1}}\nu={\mf
Q}(D_x^{-1}D_{\gL}\mu-\eta \nu)$
\end{enumerate}
It is interesting that there are no $\pp_q^{\gL}$ terms here, just shifts
$D_{\gL}$.  Moreover only $D_{\gL}\gb$ appears but pairs
$(D_{\gL}\mu,\,\mu)$ and $(D_{\gL}\nu,\,\nu)$ both appear.  
We reduce matters as in Section 8.  Thus set
$D_{\gL}\gb=\tl{\gb}$ and then from (1) and (3) we have
$(\hat{\pp}_x\sim\pp_{q^{-1}}^x$)
\bq\label{22y}
\hat{\pp}_xD_{\gL}\nu={\mf
Q}(D_x^{-1}D_{\gL}\mu-\eta\nu);\,\,\hat{\pp}_xD_{\gL}\mu=\nu
u-D_x^{-1}\tl{\gb}
\end{equation}
This means $D_{\gL}\mu=D_x[(1/{\mf Q})\hat{\pp}_xD_{\gL}\nu+\eta\nu]$ and 
$D_x^{-1}\tl{\gb}=\nu u-\hat{\pp}_xD_x\left(\frac{1}{{\mf
Q}}\hat{\pp}_xD_{\gL}\nu+\eta\nu\right)$.
Note now ${\bf (E3)}\,\,D_x\hat{\pp}_xf=D_x\pp_{q^{-1}}^xf=\pp^x_qf$ so
\bq\label{24y}
D_x^{-1}\tl{\gb}=\nu u-\hat{\pp}_x\left[\frac{1}{{\mf
Q}}\pp_q^xD_{\gL}\nu+\eta D_x\nu\right]\Rightarrow
\tl{\gb}=D_x(\nu u)-\pp_q^x\left[\frac{1}{{\mf Q}}\pp_q^xD_{\gL}\nu+\eta
D_x\nu\right]
\end{equation}
Note also $\pp_q^xD_xf=qD_x\pp_q^xf$ so
$\tl{\gb}=D_x(\nu u)-\frac{D_{\gL}}{{\mf Q}}\left(\pp_q^x\right)^2\nu-\eta q
D_x\pp_q^x\nu$
(since $D_{\gL}\pp_q^x=\pp_q^xD_{\gL}$).  Putting this in (2) now involves
computing
\bq\label{26y}
D_{\gL}\pp^x_{q^{-1}}\tl{\gb}=D_{\gL}q\pp_q^x(\nu u)-\frac{1}{{\mf Q}}
D_{\gL}^2\pp^x_{q^{-1}}(\pp_q^x)^2\nu-\eta q^2(\pp_q^x)^2\nu
\end{equation}
(note $\pp_{q^{-1}}^xD_xf=q\pp_q^xf$).  However
\bq\label{27y}
\pp^x_{q^{-1}}\pp_q^xf=\frac{1}{(q^{-1}-1)x}\left[\pp^x_{q^{-1}}f-
\pp_q^xf\right];
\,\,\pp_q^x\pp_{q^{-1}}^xf=\frac{1}{(q-1)x}
\left[\pp_q^xf-\pp^x_{q^{-1}}f\right]
\end{equation}
Then from (2), after some calculation, using
$\pp_q^x(\nu u)=\pp_q^x\nu D_xu+\nu\pp_q^xu$ and ${\bf (E4)}\,\,
{\mf Q}D_x^{-1}D_{\gL}\mu=\pp^x_{q^{-1}}D_{\gL}\nu+\eta\nu{\mf Q}\Rightarrow
D_{\gL}\mu=D_x(\eta\nu)+(1/{\mf Q}D_x\pp^x_{q^{-1}}D_{\gL}\nu=(1/{\mf
Q})\pp_q^xD_{\gL}\nu+D_x\eta\nu$, one gets (puting $\nu=\eta^2-Pu)$
\bq\label{30y}
u_t=\frac{1}{{\mf
Q}}D_{\gL}^2\pp_{q^{-1}}^x(\pp_q^x)^2Pu+\eta q^2(\pp_q^x)^2Pu
-\eta D_x^{-1}D_{\gL}(\pp_q^x)^2Pu-u\pp_q^xPu-
\end{equation}
$$-D_{\gL}q(\pp_q^xPu)D_xu-D_{\gL}qPu\pp_q^xu+D_{\gL}q\eta^2\pp_q^xu-
{\mf Q}\eta^2q\pp_q^xPu-{\mf Q}\eta^3u+{\mf Q}\eta uPu+$$
$$+{\mf Q}uD_{\gL}^{-1}D_x\eta^3-{\mf Q}uD_{\gL}^{-1}D_x\eta Pu$$
The terms in $u,\,u^2,\,\pp_xu$ which cancelled before now have the form
\bq\label{31y}
\eta^2q(D_{\gL}\pp_q^xu-{\mf Q}\pp_q^xPu);\,\,-{\mf Q}\eta^3u+{\mf
Q}uD_{\gL}^{-1}D_x\eta^3;\,\,{\mf Q}\eta u(1-D_{\gL}^{-1}D_x)Pu
\end{equation}
Then $P=1/{\mf Q}$ gives terms $(\bullet\bullet)\,\,
\eta^2q(D_{\gL}\pp_q^xu-\pp_q^xu);\,\,\eta u(1-D_{\gL}^{-1}D_x)u$,
neither of which vanish.  However from \cite{c12} one can assume $\gL\to 1$ as $q\to
1$ so both terms in $(\bullet\bullet)$ vanish in the limit as desired.  Moreover
\eqref{30y} tends to (${\mf Q}\to 2$ and $P\to 1/2$)
$u_t=\frac{1}{4}\pp_q^xu-\frac{3}{2}u\pp_xu$
as desired.
This leads one to think of a qKdV type equation based on the quantum
line to have the form 
\bq\label{34y}
u_t=\frac{1}{{\mf Q}^2}D_{\gL}^2\pp^x_{q^{-1}}(\pp_q^x)^2u+\frac{\eta}{{\mf
Q}}\left((\pp^x_{q^{-1}}D_x)^2-D_x^{-1}D_{\gL}(\pp_q^x)^2\right)u-
\end{equation}
$$-\frac{1}{{\mf
Q}}\left[D_{\gL}q(\pp_q^x)D_xu+(1+D_{\gL}q)(u\pp_q^xu)\right]+
q\eta^2(D_{\gL}-1)\pp_q^xu+\eta u(u-D_{\gL}^{-1}D_xu)$$ (since
$\pp_{q^{-1}}^xD_xf=q\pp_q^xf\sim (\pp_{q^{-1}}^xD_x)^2f=
q^2(\pp_q^x)^2f$).  If now $u=u(x,t)$ does not depend on $\gL$ one obtains
\begin{theorem}
A qKdV type equation based on the quantum line can be formulated as
\bq\label{35y}
u_t=\frac{1}{{\mf Q^2}}\pp_{q^{-1}}^x(\pp_q^x)^2u+\frac{\eta}{{\mf
Q}}\left((\pp^x_{q^{-1}}D_x)^2-D_x^{-1}(\pp_q^x)^2\right)u-
\end{equation}
$$-\frac{1}{{\mf Q}}[q(\pp_q^x)D_xu+(1+q)u\pp_q^xu]+\eta u(1-D_x)u$$
\end{theorem}

\section{STANDARD FORMS FOR QKDV}
\renewcommand{\theequation}{10.\arabic{equation}}
\setcounter{equation}{0}

There are a number of standard forms treated in e.g.
\cite{a2,a3,c4,c5,f1,h1,i6,kz,t2} and here we mainly want to exhibit specific
examples of qKdV based on $\pp_tu=\ga\pp_x^3u+\gb u\pp_xu$ for $q=1$.  That
is, we want to explicitly write out the coefficients and this does not seem to be
available in the literature.   We recall
($\tau\sim\tau_q$ in {\bf (F4)} below)
\bq\label{1w}
L^2=\pp_q^2+(q-1)xu\pp_q+u;\,\,u=\pp_q\pp_1log[\tau(x,t)D\tau(x,t)]
\end{equation}
where $L^2\sim\pp_q^2+u_1\pp_q+u_0$ with (cf. \cite{c4})
\bq\label{2w}
L=\pp_q+s_0+s_1\pp_q^{-1}+\cdots;\,\,u_1=(q-1)xu=s_0+Ds_0;\,\,\pp_{t_1}u=
\pp_qu-\pp_q^2s_0-\pp_qs_0^2
\end{equation}
Recall here ${\bf (F1)}\,\,s_0=(q-1)x\pp_1\pp_qlog\tau$ so $s_0\to 0$ as 
$q\to 1$ and $\pp_{t_1}u\to\pp_xu$ as desired. 
Some further
information in terms of tau functions is given in \cite{a3} and we will
return to that later.  Now write out $L^2$ via $(\clubsuit)\,\,
L^2=(\pp_q+s_0+s_1\pp_q^{-1}+\cdots)(\pp_q+s_0+s_1\pp_q^{-1}+\cdots)$
bearing in mind ($n\in {\bf Z}$)
\bq\label{4w}
\pp_q^nf=\sum_{k\geq 0}{\nm}_qD^{n-m}(\pp_q^mf)\pp_q^{n-m}
\end{equation}
Thus $\pp_qf=(\pp_qf)+(Df)\pp_q$ and ${\bf
(F2)}\,\,\pp_q^{-1}f=\sum_{k\geq
0}(-1)^kq^{-k(k+1)/2}D^{-k-1}(\pp_q^kf)\pp_q^{-k-1}$ which means in
particular that one has ${\bf
(F3)}\,\,\pp_q^{-1}f=(D^{-1}f)\pp_q^{-1}-q^{-1}D^{-2}\pp_qf
\pp_q^{-2}+\cdots$.  
Now the $+$ parts of $L^2$ are known to specify
\bq\label{6w}
(q-1)xu=u_1=(Ds_0)+s_0;\,\,u_0=u=(\pp_qs_0)+s_0^2+(Ds_1)+s_1;
\end{equation}
and we will write out the first $-$ parts as follows (cf. \eqref{4w}
and \eqref{41w}-\eqref{42w} below)
\bq\label{5w}
L^2_{1,2,3}=
[(\pp_qs_1)+(Ds_2)+s_0s_1+s_2+s_1D^{-1}s_0]\pp_q^{-1}+[(\pp_qs_2)+
(Ds_3)+s_3+s_0s_2-
\end{equation}
$$-q^{-1}(D^{-2}\pp_qs_0)+
s_1(D^{-1}s_1)+s_2(D^{-2}s_0]\pp_q^{-2}+[s_1q^{-3}(D^{-2}\pp_q^2s_0)-
q^{-1}s_1(D^{-2}\pp_qs_1)+s_0s_3+$$
$$\pp_qs_3+(Ds_4)
+s_1(D^{-1}s_2)+s_4+s_2D^{-2}s_1-s_2[2]_qq^{-2}(D^{-3}\pp_qs_0)+
s_3(D^{-3}s_0)]\pp_q^{-3}$$
(cf. here also \eqref{39w}).  Note that we
determined $s_0=(q-1)x\pp_1\pp_qlog\tau$ in \cite{c4} (via the expression
for $u$ in \eqref{1w}) so in principle $(Ds_1+s_1)$ is determined via
\eqref{6w} (as in \cite{c4}).  This would give $Ds_2$ if $s_1$ could be
isolated and eventually coefficients $s_i$ could then be calculated. 
However the expression for $s_1+(Ds_1)$ is not trivial since
\bq\label{7w}
s_0^2=(q-1)^2x^2(\pp_1\pp_qlog\tau)^2;\,\,\pp_qs_0=(q-1)\pp_1
[qx\pp_q^2log\tau+\pp_qlog\tau]
\end{equation}
We recall here ${\bf (F4)}\,\,\tau\sim \tau_q=\tau(c(x)+t)$ where 
$c(x)=\left(\frac{(1-q)^nx^n}{n(1-q^n)}\right)$ and $\tau$ is an ordinary
tau function for KdV.  Let us gather together some formulas now in writing
$\pp_1\pp_qlog\tau=\ga$; thus
\bq\label{9w}
u=(1+qD)\ga;\,\,s_0=(q-1)x\ga;\,\,\pp_qs_0=(q-1)[xq\pp_q\ga+\ga];\,\,
s_0^2=(q-1)^2x^2\ga^2;
\end{equation}
$$\pp_qu=\pp_q\ga+q\pp_qD\ga=(1+q^2D)\pp_q\ga$$
Note that a term $\pp_q^3u$ for example, or some variation on this involving
$D$ operations, should arise in the eventual $u_t$ equation ($t\sim t_3$)
and we may get some ideas by computing $\pp_q^3u$ from \eqref{9w}.
Thus (recall $qD\pp_q=\pp_qD$)
\bq\label{10w}
\pp_q^2u=\pp_q^2\ga+q^2\pp_qD\pp_q\ga=(1+q^3D)\pp_q^2\ga;\,\,
\pp_q^3u=(1+q^4D)\pp_q^3\ga
\end{equation}
Consider also
\bq\label{11w}
\pp_q^2s_0=(q-1)[q^2x\pp_q^2\ga+q\pp_q\ga+\pp_q\ga]=
(q-1)\{q^2x\pp_q^2\ga+[2]_q\pp_q\ga\}
\end{equation}
\bq\label{12w}
\pp_qs_0^2=\pp_q(q-1)^2(x^2\ga^2)=(q-1)^2\{q^2x^2\pp_q\ga^2+[2]_qx\ga^2\}
\end{equation}
We note now
\bq\label{13w}
\pp_q(fh)=\frac{(fh)(qx)-(fh)(x)}{(q-1)x}=Df(\pp_qh)+(\pp_qf)h=(\pp_qf)Dh+
f(\pp_qh)
\end{equation}
We could use the second form to generate D action on $\ga$ and in that
spirit consider
\bq\label{15w}
\pp_qs_0=\pp_q(q-1)x\ga=(q-1)\{D\ga+x\pp_q\ga\};\,\,\pp_q^2s_0=
(q-1)\{(q+1)D\pp_q\ga+x\pp_q^2\ga\}
\end{equation}
These have a different flavor that \eqref{11w}.  Similarly
\bq\label{16w}
\pp_qs_0^2=(q-1)^2\pp_q(x^2\ga^2)=(q-1)^2\{[2]_qxD\ga^2+x^2[\ga+
D\ga]\pp_q\ga\}
\end{equation}
\indent
We note that given (a suitable) $f$ one can define a split $f=(D+1)h$ via a
formal series
\bq\label{30w}
(1+D)^{-1}f(x)=\sum_0^{\infty}(-1)^nD^nf(x)=\sum_0^{\infty}(-1)^nf(q^nx)
\end{equation}
If e.g. $f=x^p$ this requires convergence of $(\sum_0^{\infty}(-1)^n
q^{np})x^p$ so e.g. $|q|<1$ would do.  In any case we have a formal solution
for $s_1$ via
\begin{theorem}
Given $s_1+Ds_1=u-\pp_qs_0-s_0^2=f$ there is a formal solution
\bq\label{31}
s_1(x)=\sum_0^{\infty}(-1)^nf(q^nx)
\end{equation}
Going back then to \eqref{6w} we see that $s_1$ can be determined and 
and this leads eventually to a determination of all the $s_k$.
\end{theorem}
\indent
With the coefficients now accounted for let us look for an explicit form for
qKdV based on ${\bf (F5)}\,\,L_t=[L_{+}^3,L]$.  Thus consider ($u$ and 
$u_1$ given via \eqref{6w})
$L^3=L^2L=(\pp_q^2+u_1\pp_q+u_0)(\pp_q+s_0+s_1\pp_q^{-1}+s_2\pp_q^{-2}+\cdots)$
We recall \eqref{4w}, $\pp_qf=(\pp_qf)+(Df)\pp_q$, {\bf (F2)}, {\bf (F3)},
and note ${\bf (F6)}\,\,\pp_q^2f=\pp_q[(\pp_qf)+(Df)\pp_q]=(D^2f)\pp_q^2+
\left[\begin{array}{c}
2\\
1\end{array}\right]_q
(D\pp_qf)\pp_q+(\pp_q^2f)=(D^2f)\pp_q^2+
[2]_q(D\pp_qf)\pp_q+(\pp_q^2f)=(D^2f)\pp_q^2+(q+1)(D\pp_qf)\pp_q+(\pp_q^2f)$.
Further
\bq\label{36w}
\pp_q^3f=(D^3f)\pp_q^3+\left[\begin{array}{c}
3\\
1\end{array}\right]_q(D^2\pp_qf)\pp_q^2+\left[
\begin{array}{c}
3\\
2\end{array}\right]_q(D\pp_q^2f)\pp_q+(\pp_q^3f)
\end{equation}
(where the q-brackets are equal to $(q^3-1)/(q-1)=[3]_q$.
In addition let us recall
\bq\label{41w}
\pp_q^{-1}=\sum\left[\begin{array}{c}
-1\\
m\end{array}\right]_qD^{-1-m}(\pp_q^mf)\pp_q^{-1-m};\,\,
\pp_q^{-2}f=\sum\left[\begin{array}{c}
-2\\
m\end{array}\right]_qD^{-2-m}(\pp_q^mf)\pp_q^{-2-m}
\end{equation}
etc. and a little calculation gives
\bq\label{42w}
\left[\begin{array}{c}
-n\\
k\end{array}\right]_q=\frac{[-n]_q\cdots[-n-k+1]_q}
{[k]_q[k-1]_q\cdots[1]_q};\,\,
\left[\begin{array}{c}
-1\\
m\end{array}\right]_q=(-1)^mq^{-m(m+1)/2};
\end{equation}
$$\left[\begin{array}{c}
-2\\
m\end{array}\right]_q=(-1)^m[m+1]_qq^{-(m(m+3)/2};\,\,
\left[\begin{array}{c}
-3\\
m\end{array}\right]_q=(-1)^m[m+2]_q[m+1]_qq^{-m(m+5)/2}$$
\indent
Now we calculate
\bq\label{46w}
L^3_{+}=(\pp_q^2+u_1\pp_q+u_0)(\pp_q+s_0+s_1\pp_q^{-1}+s_2\pp_q^{-2}+\cdots)=
\pp_q^3+w_2\pp_q^2+w_1\pp_q+w_0
\end{equation}
where
(recall $u_1=(q-1)xu=(1+D)s_0$ and $u=s_1+Ds_1+s_0^2+\pp_qs_0$)
\bq\label{39w}
w_2=D^2s_0+u_1=D^2s_0+Ds_0+s_0;
\end{equation}
$$w_1=(q+1)(D\pp_qs_0)+D^2s_1+[(Ds_0)+s_0](Ds_0)+
u;$$
$$w_0=\pp_q^2s_0
+(q+1)(D\pp_qs_1)+u_1\pp_qs_0+u_1(Ds_1)+us_0+D^2s_2$$
Now we note from \eqref{5w} that since $L^2_{+}=L^2$ the coefficients of 
$\pp_q^{-m}$ in \eqref{5w} must be zero and this leads to equations for $s_2,\,s_3,\,s_4$
via
\bq\label{5ww}
[(\pp_qs_1)+(Ds_2)+s_0s_1+s_2+s_1D^{-1}s_0]=0;
\end{equation}
$$[(\pp_qs_2)+
(Ds_3)+s_3+s_0s_2-q^{-1}(D^{-3}\pp_qs_0)+
s_1(D^{-1}s_1)+s_2(D^{-2}s_0]=0;$$
$$[s_1q^{-3}(D^{-3}\pp_q^2s_0)-
q^{-1}s_1(D^{-2}\pp_qs_1)+s_0s_3+\pp_qs_3+(Ds_4)+
s_1(D^{-1}s_2)+s_4-$$
$$-s_2[2]_qq^{-3}(D^{-3}\pp_qs_0)+
s_3(D^{-3}s_0)]=0$$ 
\begin{proposition}
The coefficients $s_i$ can be determined as in \eqref{5ww} using Theorem 10.1.
\end{proposition}
\indent
We note that $\pp_tL=[L^3_{+},L]$  where
$L_t=\pp_ts_0+\pp_ts_1\pp_q^{-1}+\cdots$ and since $s_0\to 0$ as $q\to 1$ we
concentrate on the $\pp_q^{-1}$ term involving $s_1$.  Thus we compute first
(omitting some calculation)
\bq\label{xxx}
[L^3_{+},L]_{-1}=
\{\pp^3s_1+[3]_q(D\pp_q^2s_2)+[3]_q(D^2\pp_qs_3)+D^3s_4
+w_2[\pp_q^2s_1+
\end{equation}
$$+[2]_q(D\pp_qs_2)]+w_2(D^3s_3)
+w_1[\pp_qs_1+(Ds_2)]+w_0s_1
-s_1q^{-3}D^{-3}\pp_q^2w_2+s_1q^{-1}D^{-2}\pp_qw_1-$$
$$-s_1D^{-1}w_0+
s_2[2]_qq^2D^{-3}\pp_q w_2-s_2D^{-2}w_1-s_3D^{-3}w_2-s_4\}\pp_q^{-1}$$
To have now $\pp_tL=[L^3_{+},L]$ requires of course some compatibility conditions
to guarantee that the coefficients of $\pp_q^3,\,\pp_q^2,\,\pp_q$ vanish, but these
are in fact automatic (see below).
Then, omitting $\pp_ts_0$ (since $s_0\to 0$ as $q\to 1$ and this gives us no
check), one can write
\bq\label{zzz}
\pp_ts_1=\pp_q^3s_1+[3]_qD\pp_q^2s_2+[3]_qD^2\pp_qs_3+D^3s_4+w_2[(\pp_q^2s_1)+
[2]_q(D\pp_qs_2]+w_2D^3s_3+
\end{equation}
$$w_1[\pp_qs_1+Ds_2]+w_0s_1-s_1q^{-3}D^{-3}w_2+ s_1q^{-1}D^{-2}\pp_qw_1-$$
$$-s_1D^{-1}w_0+s_2
[2]_qq^2D^{-3}w_2-s_2D^{-2}w_1-s_3D^{-3}w_2-s_4$$
\indent
Let us first check some of this for $q\to 1$ where $s_0\to 0$ and $2s_1=u$
with $w_2\to 0,\,\,w_1\to 3s_1,$ and $w_0\to s_2=-(1/2)\pp s_1$ (cf. \eqref
{39w}, \eqref{5ww}, etc.).  Further from \eqref{5ww} one has $2s_3+\pp s_2
+s^2_1=0\Rightarrow s_3=-(1/2)s_1^2+(1/4)\pp^2s_1$ and $s_1\pp s_1+\pp
s_3+2s_4+s_1s_2=0\Rightarrow s_4=-(1/8)\pp^3s_1+(1/2)s_1\pp s_1$.  One can then
write \eqref{zzz} as
\bq\label{xx}
\pp_ts_1\to \pp^3s_1+3\pp^2s_2+3\pp s_3+s_4+w_1(\pp s_1+s_2)+w_0s_1+s_1\pp
w_1-s_1w_0-
\end{equation}
$$-s_2w_1-s_4=
\pp^3s_1+3\pp^2\left(-\frac{1}{2}\pp
s_1\right)+3\pp\left(-\frac{1}{2}s_1^2+\frac{1}{4}\pp^2s_1\right)+3s_1(\pp
s_1-(1/2)\pp s_1)+s_1(\pp 3s_1)-$$
$$-3s_1(-(1/2)\pp s_1)
=\frac{1}{4}\pp^3s_1+3s_1\pp s_1$$
which has the correct KdV form.  Hence one can suggest
\begin{theorem}
The qKdV equation can be extracted from \eqref{zzz}.
\end{theorem}
\indent
However \eqref{zzz} is awkward at best and there is a better procedure which
expresses matters directly in terms of $u$.  Thus
first note from \cite{kz} that for for PSDO $L$ in ${\mf D}_q=x\pp_q$ with 
$\pp_mL=[L^m_{+},L]$ we have
$[L^m,L]=[L^m_{+},L]+ [L^m_{-},L]=0$  
(recall that in \cite{kz} $D_q$ is replaced by
$xD_q=(q-1)^{-1}\gD_q={\mf D}_q$ with greater ease in calculation since
$D{\mf D}_q={\mf D}_qD$).  
This means that
$deg\,[L^m_{+},L]=deg\,[L^m_{-},L]=deg\,L^m_{-}+det\,L=-1+1=0$ and this argument
clearly carries over to PSDO $D_q\sim\pp_q$.  
Consequently
$deg\,\pp_mL$ is zero and     
hence the compatibility conditions mentioned above 
(before \eqref{zzz}) are automatically satisfied (i.e. $\pp_tL=\pp_ts_0+\cdots$ so
$[L^3_{+},L]$ must be of degree zero).  
Now for qKP equations one might
imagine working with 
\bq\label{xyz}
\pp_2\psi=L^2_{+}\psi;\,\,\pp_3\psi=L^3_{+}\psi;\,\,\pp_3\pp_2\psi=\pp_3L^2_{+}\psi
+L^2_{+}\pp_3\psi=
\end{equation}
$$=\pp_2\pp_3\psi=\pp_2L_{+}^3\psi+L_{+}^3\pp_2\psi\Rightarrow
\pp_3L^2_{+}=\pp_2L^3_{+}+[L_{+}^3,L_{+}^2]$$
Thus for qKdV where $L$ does not depend on $t_2$ and $L^2_{+}=L^2$ one has for 
$t\sim t_3$ the result ${\bf (F9)}\,\,\pp_tL^2=[L_{+}^3,L^2]$.  This can be
written then as 
\bq\label{xzy}
\pp_tu=[L^3_{+},L^2]_0=[\pp_q^3+w_2\pp_q^2+w_1\pp_q+w_0,\pp_q^2+u_1\pp_q+u]_0=
\end{equation}
$$(\pp_q^3u)+w_2(\pp_q^2u)+w_1(\pp_qu)-[(\pp_q^2w_0)+u_1(\pp_qw_0)]$$
\indent
{\bf REMARK 10.1.}
We have not seen equations like \eqref{zzz} or \eqref{xzy} 
written out before (evidently for good reason) and the
motivation here was to compare this with qKdV type equations derived by
other means in other contexts (e.g. zero
curvature derivations).  Noncommutative KdV (using Moyal type brackets) is
also of interest here (cf. \cite{cxx,d3,fyy,hzz,pyz,yzz}) and the equations are much
simpler. Other quantum versions of KdV are indicated in \cite{fz} and references
there; we refer to \cite{kq,sz} for general discrete situations.
One can also consider the Frenkel form (cf. \cite{a2,a3,f1,kz} and see
\cite{jzz,jzx,lzz} for other Lax forms)
\bq\label{1ww} 
Q=D+\sum_0^{\infty}u_iD^{-i};\,\,Q_q=D_q+\sum_0^{\infty}v_iD_q^{-i}
\end{equation}
\bq\label{2ww}
\pp_nQ=[A^n_{+},Q];\,\,\pp_nQ_q=[(Q_q^n)_{+},Q_q]
\end{equation}
The calculations in \cite{a2,a3,f1,dz} (and in \cite{c4} and Sections 1-4 of
this paper) never actually produce an explicit formula for qKdV.  There are
equations for the coefficients in terms of tau functions and relations
between $u_i$ and $v_i$ but nothing explicitly in a form based on
$\pp_tu=a\pp_x^3u+bu\pp_xu$.  In \eqref{zzz}, \eqref{xzy}, and calculations below
we hope to remedy this situation.  
$\hfill\bs$
\\[3mm]\indent
We want to write \eqref{xzy} now in terms of $u$ and $u_1=u_1(u)$ so recall
$u_1=(q-1)xu=s_0+Ds_0$ and $u=s_1+Ds_1+s_0^2+\pp_qs_0$ with (from \eqref{39w})
$w_2=D^2s_0+u_1,\,\,w_1=(q+1)D\pp_qs_0+D^2s_1+u_1Ds_0+u,$ and
$w_0=\pp_q^2s_0+(q+1)D\pp_qs_1+u_1\pp_qs_0+u_1Ds_1+us_0+D^2s_2$.  Here
$s_2+Ds_2=-s_1D^{-1}s_0-\pp_qs_1-s_0s_1$ from \eqref{5ww}.
Now before going further we check \eqref{xzy} for $q\to 1$.  Then
\bq\label{s0s}
2s_0=0;\,\,2s_1=u;\,\,w_2=u_1=0;\,\,w_1=s_1+u;\,\,2s_2=-\pp s_1;\,\,w_0=
(3/2)\pp s_1
\end{equation}
Hence \eqref{xzy} becomes ($s_1=u/2$)
\bq\label{s1s}
\pp_tu=\frac{1}{4}\pp^3u+\frac{3}{2}u\pp u
\end{equation}
which is fine (and agrees with \eqref{xx}).  To directly express \eqref{xzy} in
terms of u alone requires an expression for the $w_i$ as functions of u but here we
see that an explicit finite form is unlikely.  This is seen immediately in
$u_1=(q-1)xu=s_0+Ds_0$ for example which has a solution
\bq\label{s2s}
s_0=(1+D)^{-1}u_1=\sum_0^{\infty}(-1)^nu_1(q^nx)=(q-1)x\sum_0^{\infty}
(-1)^nq^nu(q^nx)
\end{equation}
Similarly from Theorem 9.1 $s_1=\sum_0^{\infty}(-1)^nf(q^nx)$ for
$f=u-\pp_qs_0-s_0^2$.  Thereafter we can write e.g.
\bq\label{s3s}
w_2=(q-1)(D^2+D+1)x\sum_0^{\infty}(-1)^nq^nu(q^nx)=
(q-1)x\left[\sum_2^{\infty}(-1)^pq^pu(q^px)-\right.
\end{equation}
$$\left.-\sum_1^{\infty}(-1)^sq^su(q^su(q^sx)+\sum_0^{\infty}
(-1)^nq^nu(q^nx)\right]=(q-1)x\left[\sum_0^{\infty}(-1)^nq^nu(q^nx)+qu(qx)\right]$$
One can also compute $w_1$ and $w_0$ in terms of infinite series in u and one has
\begin{theorem}
The qKdV equation has the form \eqref{xzy} where $u_1=(q-1)xu$ and $w_i=w_i(u)$.
\end{theorem}

\end{document}